\documentclass{article}
\usepackage[T1]{fontenc}  
\usepackage[utf8]{inputenc}
\usepackage{hyperref}
\usepackage{textalpha}
\usepackage{amsmath,amssymb,graphicx,subfig,mathrsfs,amsthm}

\begin{document}
\title{Book I of Euclid's {\em Elements} and application of areas} 
\author{Jordan Bell\\ \texttt{jordan.bell@gmail.com}}
\date{\today}

\maketitle

\begin{abstract}
We work through Book I of Euclid's Elements with our focus on application of areas (I.42, I.44, I.45). We summarize alternate constructions from medieval editions of Euclid's elements and ancient and medieval commentaries. We remark that Euclid's proof of I.44 involves a seldom commented on use of superposition, but that several medieval editions of Euclid give constructions that avoid the use of superposition. This use of superposition is also avoided in Ralph Abraham's ``VCE: The Visual Constructions of Euclid'' C\#12, C\#12B at \url{http://www.visual-euclid.org/vce/contents.html}

We collate the figures with the digitized editions of Euclid at
(P) Biblioteca Apostolica Vaticana (BAV), Vat. gr. 190,
(F) Biblioteca Medicea Laurenziana (BML), Plut. 28.03,
(B) Bodleian, MS. D’Orville 301,
(V) Österreichische Nationalbibliothek, Cod. Phil. gr. 31,
(b) Biblioteca Comunale dell'Archiginnasio, Collocazione A 19,
(p) Bibliothèque nationale de France, Grec 2466.
\end{abstract}

\tableofcontents

\section{Introduction}
Taisbak \cite[pp.~28--29]{data}:

\begin{quote}
It may be appropriate to introduce the {\em The Helping Hand}, a well-known factotum in Greek
geometry, who takes care that lines are drawn, points are taken, circles described,
perpendiculars dropped, etc. The perfect imperative passive is its verbal mask: `Let a
circle have been described with centre $A$ and radius $AB$'; `let it lie given' {\em keistho
dedomenon}. No one who has done the {\em Elements} in Greek will have missed it; never is
there any of the commands or exhortations so familiar from our own class-rooms:
`Draw the median from vertex $A$', or `If we cut the circle by that secant', or `Let us add
those squares together'. Always {\em The Helping Hand} is there first to see that things are
done, and to keep the operations free from contamination by our mortal fingers.
\end{quote}

On the perfect passive imperative, cf.
Priscian, {\em Institutiones Grammaticae}: Book 8, Keil, 1870, II: 406, 15--27; 407, 1--9;
Book 18, Keil, 1870, III: 238, 12--26.

Netz \cite[p.~xvii]{netz}:

\begin{quote}
As is explained in chapter 1, most of the diagrams in Greek mathematical
works have not yet been edited from manuscripts. The figures
in modern editions are reconstructions made by modern editors, based
on their modern understanding of what a diagram should look like.
However, as will be argued below, such an understanding is culturally
variable. It is therefore better to keep, as far as possible, to the diagrams
as they are found in Greek manuscripts (that is, generally
speaking, in Byzantine manuscripts).
\end{quote}

Netz \cite[p.~16]{netz}:

\begin{quote}
Diagrams, as a rule, were not drawn on site.
The limitations of the media available suggest, rather, the preparation
of the diagram prior to the communicative act -- a consequence of the
inability to erase.
\end{quote}

Netz \cite[p.~25]{netz}:

\begin{quote}
What we see, in short is that while the text is being worked through,
the diagram is assumed to exist. The text takes the diagram for granted.
This reflects the material implementation discussed above. This, in
fact, is the simple explanation for the use of {\em perfect} imperatives in the
references to the setting out -- `let the point $A$ have been taken'. It
reflects nothing more than the fact that, by the time one comes to
discuss the diagram, it has already been drawn.
\end{quote}

Netz \cite[pp.~94--95]{netz}:

\begin{quote}
That numbers are absent from the original is not just an accident,
the absence of a tool we find useful but the Greeks did not require.
The absence signifies a different approach to definitions. The text of
the definitions appears as a continuous piece of prose, not as a discrete
juxtaposition of so many definitions. \dots. So the principle
is this. Mathematical texts start, most commonly, with some piece of
prose preceding the sequence of proved results. Often, this is developed
into a full `introduction', usually in the form of a letter (prime
examples: Archimedes or Apollonius). Elsewhere, the prose is very
terse, and supplies no more than some reflections on the mathematical
objects (prime example: Euclid).

I suggest that we see the shorter, Euclid-type introduction as an
extremely abbreviated, impersonal variation upon the theme offered
more richly in Archimedes or Apollonius. Then it becomes possible to 
understand such baffling `definitions' as, e.g., {\em Elements} I.3: `and the
limits of a line are points'. This `definition' is not a definition of any
of the three nouns it contains (lines and points are defined elsewhere,
and no definition of limits is required here). It is a brief second-order
commentary, following the definitions of `line' and `point'. Greek 
mathematical works do not start with definitions. They start with
second-order statements, in which the goals and the means of the work are
settled. Often, this includes material we identify as `definitions'. In
counting definitions, snatches of text must be taken out of context, and
the decision concerning where they start is somewhat arbitrary. (Bear 
in mind of course that the text was written -- even in late manuscripts -- 
as a continuous, practically unparagraphed whole.)

\dots.

Most definitions do not prescribe equivalences
between expressions (which can then serve to abbeviate, no more).
They specify the situations under which properties are considered to
belong to objects. Now that we see that most definitions are simply
part of the introductory prose, this makes sense. There is no metamathematical
theory of definition at work here. Before getting down to
work, the mathematician describes what he is doing -- that's all.
\end{quote}

Netz \cite[p.~238]{netz}:

\begin{quote}
The first floor of Greek mathematics is the general tool-box,
Euclid's results. To master it, even superficially, is to become a passive
mathematician, an initiate.

The second floor is made up of such works as the first four books of
Apollonius' {\em Conics}; other, comparable works are, e.g. works on trigonometry.
Such results are understood passively even by the passive
mathematician, the one who knows no more than Euclid's {\em Elements}.
They are mastered by a creative mathematician in the same way in
which the first floor is mastered by the passive mathematician.

The third floor relies on the results of the second floor. This is
where the professionals converse. They have little incentive and occasion
to master this level in the way in which the second floor is mastered.
What is more important, the cognitive situation excludes such
mastery, for oral storing and retrieval have their limitations. Already
the second floor is invoked only through some very specific results.
When Archimedes entered the scene, the tool-box was already full.
Mathematics would explode exponentially only when storing and retrieval
became much more written, and when the construction of the
tool-box was done methodically rather than through sheer exposure.
\end{quote}

Proclus 83--84 \cite[pp.~68--69]{proclus}:

\begin{quote}
The book is divided into three major parts. The first reveals
the construction of triangles and the special properties of their
angles and their sides, comparing triangles with one another
as well as studying each by itself. Thus it takes a single triangle
and examines now the angles from the standpoint of the sides
and now the sids from the standpoint of the angles, with 
respect to their equality or inequality; and then, assuming two
triangles, it investigates the same properties in various ways.
The second part develops the theory of parallelograms, beginning
with the special characteristics of parallel lines and
the method of constructing the parallelograms and then
demonstrating the properties of parallelograms. The third part
reveals the kinship between triangles and parallelograms both
in their properties and in their relations to one another. Thus
it proves that triangles or parallelograms on the same or equal
bases have identical properties; it shows [what is the relation
between] a triangle and a parallelogram on the same base,
how to construct a parallelogram equal to a triangle, and
finally, with respect to the squares on the sides of a right-angled
triangle, what is the relation of the square on the side
that subtends the right angle to the squares on the two sides
that contain it. Something like this may be said to be the
purpose of the first book of the {\em Elements} and the division 
of its contents.
\end{quote}

First part: 1--26, second part: 27--34, third part: 35--48; cf. Proclus 352--353 \cite[p.~275]{proclus} and 395 \cite[p.~311]{proclus}.

Proclus 203 \cite[p.~159]{proclus}:

\begin{quote}
Every problem and every theorem that is furnished with all its parts
should contain the following elements: an enunciation, an
exposition, a specification, a construction, a proof, and a
conclusion. Of these the enunciation states what is given and
what is being sought from it, for a perfect enunciation consists
of both these parts. The exposition takes separately what is
given and prepares it in advance for use in the investigation.
The specification takes separately the thing that is sought and
makes clear precisely what it is. The construction adds what
is lacking in the given for finding what is sought. The proof
draws the proposed inference by reasoning scientifically from
the propositions that have been admitted. The conclusion
reverts to the enunciation, confirming what has been proved.
\end{quote}

Netz \cite[pp.~10--11]{netz}:

\begin{itemize}
\item[] enunciation={\em protasis}
\item[] exposition=setting out={\em ekthesis}
\item[] specification=definition of goal={\em diorismos}
\item[] construction={\em kataskeu\={e}}
\item[] proof={\em apodeixis}
\item[] conclusion={\em sumperasma}
\end{itemize}

Mueller \cite[p.~11]{mueller}: 
{\em protasis}, {\em ekthesis}, {\em diorismos}, {\em kataskeu\={e}}, {\em apodeixis}, {\em sumperasma} 

Al-Nayrizi \cite[p.~102]{alnayriziI}:

\begin{quote}
The figures, all of them, theorems and constructions, have
been named with a common name, and each one of them, namely,
theorem and construction (and locating too, if it is something else
apart from the two of them), is divided into six divisions, namely,
proposition, exemplification, separation, construction, proof, and
conclusion.
\end{quote}

Al-Nayrizi \cite[p.~102]{alnayriziI}:

\begin{quote}
The separation is what separates what is requested in the
proposition, what is set down in the exemplification, from its
common genus and requests that it be constructed and proved.
\end{quote}

Al-Nayrizi \cite[p.~103]{alnayriziI}:

\begin{quote}
The conclusion is what teaches the proposition, as when
you say, ``We have now proven that in every triangle, the three angles
are truly equal to two right angles.'' We say it with confidence since it
has been proven, and for that reason we do not add anything at all to
it except ``therefore.'' 
\end{quote}

Proclus 207 \cite[p.~162]{proclus}:

\begin{quote}
Furthermore, mathematicians are accustomed to draw what
is in a way a double conclusion. For when they have shown
something to be true of the given figure, they infer that it is
true in general, going from the particular to the universal
conclusion. Because they do not make use of the particular
qualities of the subjects but draw the angle or the straight
line in order to place what is given before our eyes, they
consider that what they infer about the given angle or straight
line can be identically asserted for every similar case. They
pass therefore to the universal conclusion in order that we
may not suppose that the result is confined to the particular
instance. This procedure is justified, since for the demonstration
they use the objects set out in the diagram not as these
particular figures, but as figures resembling others of the same
sort. It is not as having such-and-such a size that the
angle before me is bisected, but as being rectilineal and 
nothing more. Its particular size is a character of the given angle,
but its having rectilineal sides is a common feature of all rectilineal
angles. Suppose the given angle is a right angle. If
I used its rightness for my demonstration, I should not be
able to infer anything about the whole class of rectilineal
angles; but if I make no use of its rightness and consider 
only its rectilineal character, the proposition will apply
equally to all angles with rectilineal sides.
\end{quote}

Proclus 208 \cite[pp.~162--163]{proclus}:

\begin{quote}
Let us view the things that have been said by applying
them to this our first problem. Clearly it is a problem, for it
bids us  devise a way of constructing an equilateral triangle.
In this case the enunciation consists of both what is given and
what is sought. What is given is a finite straight line, and what
is sought is how to construct an equilateral triangle on it. The
statement of the given precedes and the statement of what is
sought follows, so that we may weave them together as ``If
there is a finite straight line, it is possible to construct an
equilateral triangle on it.''
\end{quote}

Proclus 208 \cite[p.~163]{proclus}:

\begin{quote}
Next after the enunciation is the exposition: ``Let this be
the given finite line.'' You see that the exposition itself mentions
only the given, without reference to what is sought.
Upon this follows the specification: ``It is required to construct
an equilateral triangle on the designated finite straight
line.''  In a sense the purpose of the specification is to fix our
attention; it makes us more attentive to the proof by announcing
what is to be proved, just as the exposition puts us in
a better position for learning by producing the given element
before our eyes.
\end{quote}

Book I, Definitions \cite[pp.~153--154]{euclidI}:

\begin{itemize}
\item[8] A \textbf{plane angle} is the inclination to one another of
two lines in a plane which meet one another and do not lie in
a straight line.

\item[9] And when the lines containing the angle are straight,
the angle is called \textbf{rectilineal}.

\item[10] When a straight line set up on a straight line makes
the adjacent angles equal to one another, each of the equal
angles is \textbf{right}, and the straight line standing on the other is
called a \textbf{perpendicular} to that on which it stands.

\item[13] A \textbf{boundary} is that which is an extremity of anything.

\item[14] A \textbf{figure} is that which is contained by any boundary
or boundaries.

\item[15] A \textbf{circle} is a plane figure contained by one line
such that all the straight lines falling upon it it from one point among
those lying within the figure are equal to one another;

\item[16] And the point is called the \textbf{centre} of the circle.

\item[22] Of quadrilateral figures, a \textbf{square} is that which is
both equilateral and right-angled; an \textbf{oblong} that which is
right-angled but not equilateral; a \textbf{rhombus} that which is
equilateral but not right-angled; and a \textbf{rhomboid} that which
has its opposite sides and angles equal to one another but is
neither equilateral nor right-angled. And let quadrilaterals
other than these be called \textbf{trapezia}.

\item[23] \textbf{Parallel} straight lines are straight lines which,
being in the same plane and being produced indefinitely in
both directions, do not meet one another in either direction.
\end{itemize}

Book I, Postulates \cite[pp.~154--155]{euclidI}:

\begin{enumerate}
\item To draw a straight line from any point to any point.
\item To produce a finite straight line continuously in a
straight line.
\item To describe a circle with any centre and distance.
\item That all right angles are equal to one another.
\item That, if a straight line falling on two straight lines
make the interior angles on the same side less than two right
angles, the two straight lines, if produced indefinitely, meet
on that side on which are the angles less than the two right
angles.
\end{enumerate}

If $A,B,C,D$ are angles, the statement that $A$ and $B$ are less than $C$ and $D$ means that
$A+B<C+D$, not that $A<C$ and $B<D$. The statement that $A$ and $B$ are less than $C$ and $D$ respectively means that
$A<C$ and $B<D$. Likewise, the statement that $A$ and $B$ are equal to $C$ and $D$ means that $A+B=C+D$, and the statement that $A$ and $B$ are equal to
 $C$ and $D$ respectively means that $A=C$ and $B=D$. 

Plural predication in Plato, {\em Hippias Major} 300d--e \cite[p.~259]{hippias}:

\begin{quote}
Hippias: You'll have  finer knowledge than anyone whether or not I'm
playing games, Socrates, when you try to describe these notions of yours
and are shown to be talking nonsense. It's quite impossible -- you'll never
find an attribute which neither you nor I have, but which both of us have.

Socrates: Are you sure, Hippias? I suppose you've got a point, but
{\em I} don't understand. Let me explain more clearly what I'm getting at: it
seems to me that both of us together may possess as an attribute something
which neither I have as an attribute nor am (and neither are you); and, to
put it the other way round, that neither of us, as individuals, may be
something which both of us together have as an attribute.
\end{quote}

Socrates ironically says the following, 301d--e \cite[p.~260]{hippias}:

\begin{quote}
You see, before you spoke, my friend,
we were so inane as to believe that {\em each} of us -- you and I -- is one, but that
both of us together, being two not one, are not what each individual is.
See how stupid we were! But now we know better: you've explained that
if both together are two, then each individual must be two as well; and if
each individual is one, both must be one as well.
\end{quote}

Heath \cite[p.~201]{euclidI}:

\begin{quote}
As to the {\em raison d'\^etre} and the place of Post. 4 one thing is quite certain.
It was essential from Euclid's point of view that it should come before Post. 5,
since the condition in the latter that a certain pair of angles are together less
than two right angles would be useless unless it were first made clear that
right angles are angles of determinate and invariable magnitude.
\end{quote}

Common Notions \cite[p.~155]{euclidI}:

\begin{enumerate}
\item Things which are equal to the same thing are also equal to one another.
\item If equals be added to equals, the wholes are equal.
\item If equals be subtracted from equals, the remainders
are equal.
\item Things which coincide with one another are equal to
one another.
\item The whole is greater than the part.
\end{enumerate}

\section{Book I}
Conspectus siglorum: 

\begin{itemize}
\item[P] Biblioteca Apostolica Vaticana (BAV), Vat. gr. 190
\item[F] Biblioteca Medicea Laurenziana (BML), Plut. 28.03
\item[B] Bodleian, MS. D'Orville 301
\item[V] \"Osterreichische Nationalbibliothek, Cod. Phil. gr. 31
\item[b] Biblioteca Comunale dell'Archiginnasio, Collocazione A 19
\item[p] Bibliothèque nationale de France, Grec 2466
\end{itemize}

I.1: ``On a given finite straight line to construct an equilateral
triangle.''

I.2: ``To place at a given point (as an extremity) a straight line equal to a given straight line.''

\begin{figure}
\begin{center}
\includegraphics[width=\textwidth]{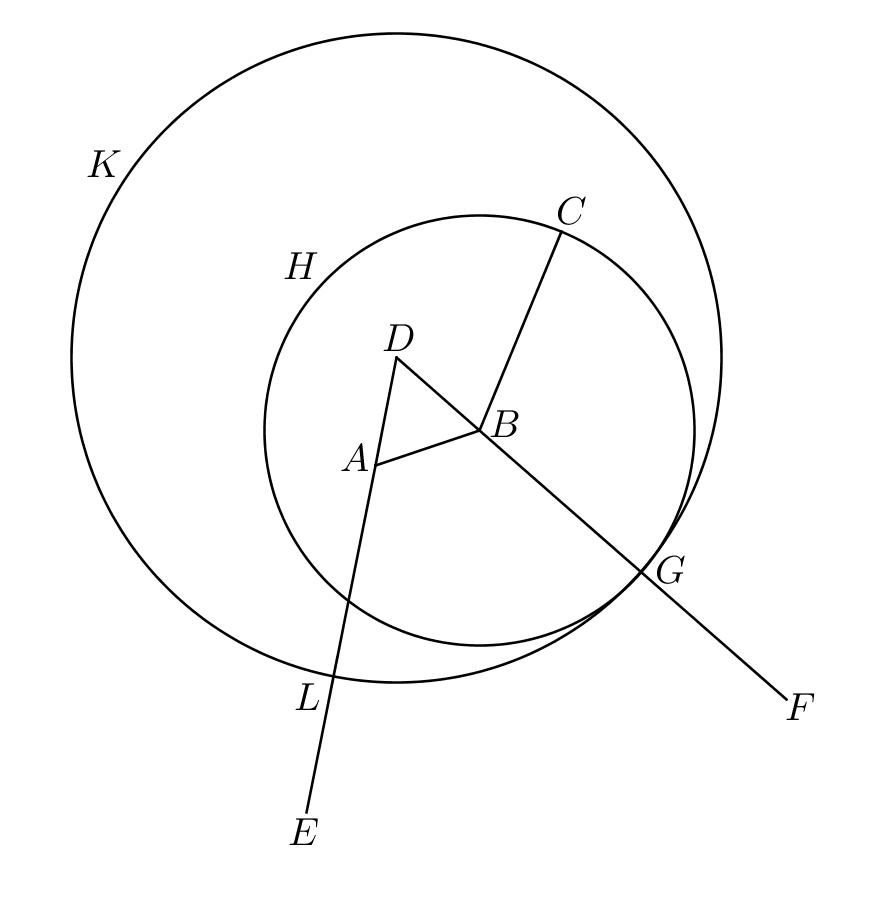}
\end{center}
\caption{I.2: P 15v, F 2r, B 8v, V 10}
\label{I2}
\end{figure}

\begin{proof}
{\em ekthesis}: Let the point $A$ be given and let  the straight line $BC$ be given.

{\em diorismos}: Thus it is required to place at the point $A$ as an extremity 
a straight line equal to the given straight line $BC$.

{\em kataskeu\={e}}: Let the  straight line $AB$ be joined from $A$ to $B$ (Postulate 1).
On this straight let the equilateral triangle $DAB$ be constructed (I.1).
Let the straight line $AE$ be produced in a straight line with $DA$ and let the straight line
$BF$ be produced in a straight line with $DB$ (Postulate 2).
Let the circle $CGH$ be described with center $B$ and distance $BC$ (Postulate 3), intersecting the line
$BF$ at $G$; and
let the circle $GKL$ be described with center $D$ and distance $DG$ (Postulate 3), intersecting the line
$AE$ at $L$.

{\em apodeixis}: Because $B$ is the center of the circle $CGH$, $BC$ is equal to $BG$.
Because $D$ is the center of the circle $GKL$, $DL$ is equal to $DG$.
And $DA$ is equal to $DB$, so when $DA$ is subtracted from $DL$ and 
$DB$ is subtracted from $DG$, the remainders $AL$ and $BG$ are equal (Common Notion 3).
But $BC$ is equal to $BG$, so each of the straight lines $AL,BC$ is equal to $BG$.
Therefore $AL$ is equal to $BC$ (Common Notion 1).

{\em sumperasma}: Therefore the straight line $AB$, equal to the given straight line $BC$,
has been placed at the given point $A$.
\end{proof}

Proclus 222--223 \cite[p.~174]{proclus}:

\begin{quote}
Some problems have no cases, while others have many;
and the same is true of theorems. A proposition is said to have
cases when it has the same force in a variety of diagrams,
that is, can be demonstrated in the same way despite changes
in position, whereas one that succeeds only with a single position
and a single construction is without cases.
\end{quote}

1.3: ``Given two unequal straight lines, to cut off from the 
greater a straight line equal to the less.''

I.4: ``If two triangles have the two sides equal to two sides
respectively, and have the angles contained by the equal straight
lines equal, they will also have the base equal to the base, the
triangle will be equal to the triangle, and the remaining angles
will be equal to the remaining angles respectively, namely those
which the equal sides subtend.''

Proclus 235--236 \cite[pp.~183--184]{proclus}:

\begin{quote}
There are three things proved and two things given about these triangles. One of the given elements
is the equality of two sides (really two given sides, but obviously given in ratio to one another) and the equality
of the angles contained by the equal sides. And the things to be proved are three: the
equality of  base to base, the equality of triangle to triangle,
and the equality of the other angles to the other angles.
Since it would be possible for the triangles to have two sides 
equal to two sides and yet the theorem to be false because the sides
are not equal one to another but one pair to the other pair,
he did not simply say, in his statement of the given, that
the lines are equal, but that they are equal ``respectively.'' For
if it should happen that one of the triangles had one side of
three units and the other of four, while the other triangle had
one side of five units and another of two (the angle included between them being
a right angle), the two sides of the one 
would be equal to the two sides of the other, since their sum
is seven in each case. But this would not show the one triangle
equal to the other; for the area of the former is six, of the latter five.
\end{quote}

Proclus 236 \cite[p.~184]{proclus}:

\begin{quote}
As to the ``base'' of a triangle, when no side has previously been named,
we must suppose it to denote the side towards the observer, but when two sides have already been mentioned,
it must mean the remaining side.
\end{quote}

Proclus 236 \cite[p.~185]{proclus}:

\begin{quote}
Two triangles are said to be equal when their areas are equal. It can happen that
two triangles with equal perimeters have unequal areas because of the inequality of their
angles. ``Area'' I call the space itself which is cut off by the sides of the triangle,
and ``perimeter'' the line composed of the three sides of the triangle.
\end{quote}

Proclus 237--238 \cite[p.~185]{proclus}:

\begin{quote}
Base is said to be equal to base and generally a straight line
to another straight line when the congruence of their extremities
makes the whole of the one line coincide with the whole
of the other. Every straight line coincides with every other,
and in the case of equal lines their extremities also coincide.
A rectilineal angle is said to be equal to a rectilineal angle when,
if one of the sides containing it is placed upon one of 
the sides containing the other, the second side of the first coincides
with the second side of the second.
\end{quote}

Proclus 238 \cite[p.~186]{proclus}:

\begin{quote}
This also must  be understood in advance, that the side that
lies opposite an angle is said to subtend it. Every angle in a 
triangle is contained by two sides of the triangle and subtended by the other.
When the other 
sides fail to coincide, that angle is greater whose side falls
outside, and that angle less whose side falls inside. For in the
one case the one angle includes the other, in the other case it
is included by the other.
\end{quote}

I.7: ``Given two straight lines constructed on a straight line
(from its extremities) and meeting in a point, there cannot be
constructed on the same straight line (from its extremities),
and on the same side of it, two other straight lines meeting in
another point and equal to the former two respectively, namely
each to that which has the same extremity with it.''

I.8: ``If two triangles have the two sides equal to two sides
respectively, and have also the base equal to the base, they will
also have the angles equal which are contained by the equal
straight lines.''

\begin{figure}
\begin{center}
\includegraphics[width=\textwidth]{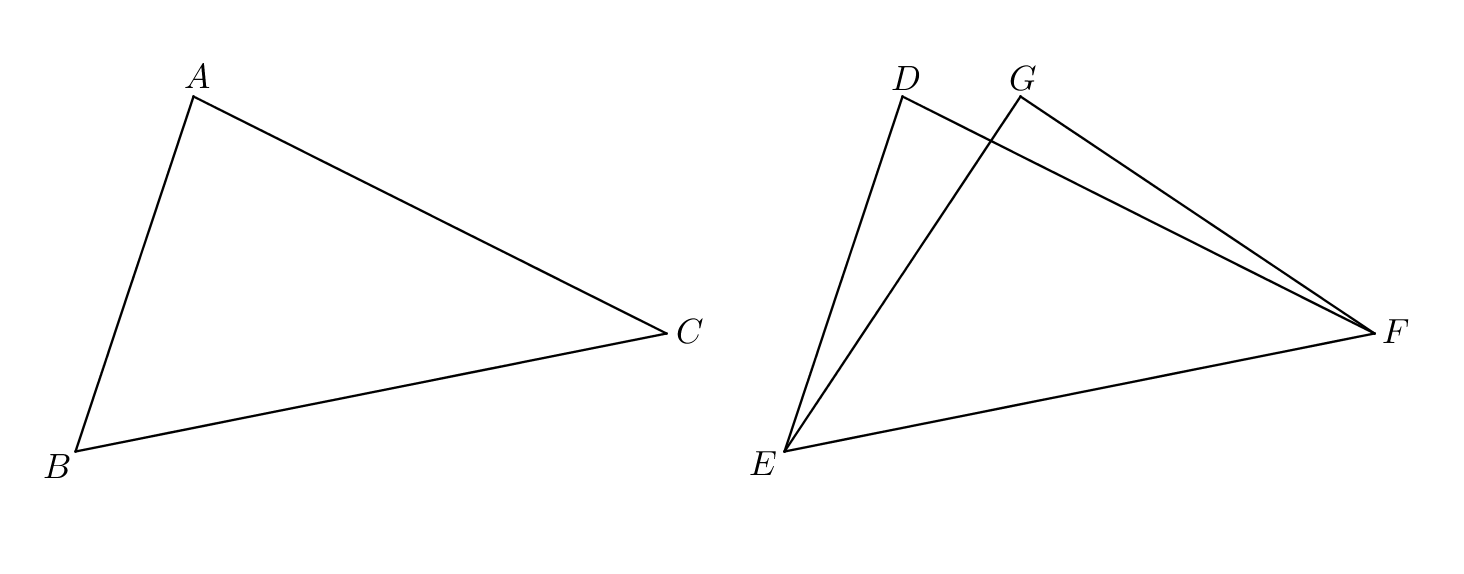}
\end{center}
\caption{I.8: P 21r, F 3v, B 12r, V 17}
\label{I8}
\end{figure}

\begin{proof}
{\em ekthesis}: Let $ABC,DEF$ be two triangles with the two sides $AB,AC$ equal to the two sides $DE,DF$ respectively
and with the base $BC$ equal to the base $EF$. 

{\em diorismos}: I say that the angle $BAC$ is equal to the angle $EDF$.

{\em kataskeu\={e}}: Let the triangle $ABC$ be applied to the triangle $DEF$ with the point $B$ placed on the point $E$
and the straight line $BC$ on $EF$.

{\em apodeixis}: Then the point $C$ will coincide with the point $F$ because
$BC$ is equal to $EF$.
If the base $BC$ coincides with the base $EF$ and the sides $BA,AC$ do not coincide with $ED,DF$ but
 as $EG,EF$ miss them, then there have been constructed on the straight line
$EF$ the straight lines $ED,FD$ meeting at $D$, and on the same side
the straight lines $EG,FG$
meeting at $G$. $AB,AC$ are equal respectively to $DE,DF$ and are also
equal respectively to $GE,GF$, so 
$DE,DF$ are equal respectively to $GE,GF$ (Common Notion 1). But such straight lines cannot be so constructed (I.7). 
Therefore, with the base $BC$ having been applied to the base $ED$, it is impossible for the sides $BA,AC$ not to coincide with the sides
$ED,DF$; therefore the sides $BA,AC$ will coincide with the sides $ED,DF$, and then
the angle $BAC$ will coincide with the angle $EDF$, and therefore will be equal to it (Common Notion 4).

{\em sumperasma}: Therefore, etc.
\end{proof}

I.9: ``To bisect a given rectilineal angle.''

I.10: ``To bisect a given finite straight line.''

I.11: ``To draw a straight line at right angles to a given straight
line from a given point on it.''

I.12: ``To a given infinite straight line, from a given point
which is not on it, to draw a perpendicular straight line.''

\begin{figure}
\begin{center}
\includegraphics[width=\textwidth]{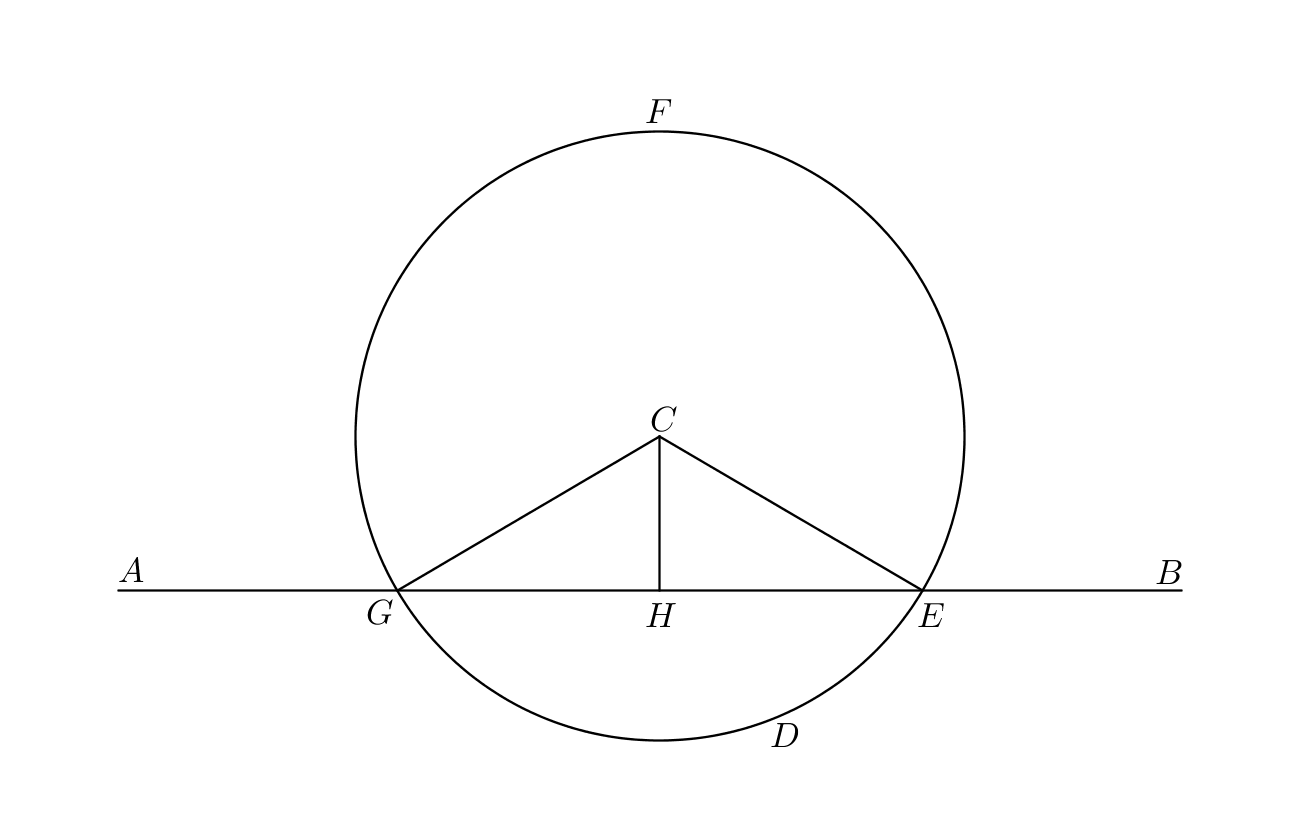}
\end{center}
\caption{I.12: P 23r, F 4v, B 13v, V 19}
\label{I12}
\end{figure}

\begin{proof}
Let $AB$ be the given infinite straight line and let $C$ be the given point not on this line.
Let a point $D$ be taken on the other side of $AB$ and let the circle $EFG$ be described with center $C$ and
distance $CD$ (Postulate 3).
Let the straight line $EG$ be bisected at $H$ (I.10).
Let the straight lines $CG, CH, CE$ be joined (Postulate 1).

$GHC$ and $EHC$ are triangles with a common side $HC$. 
As $GH$ is equal to $EH$, the two sides $GH,HC$ are equal to the two sides
$EH,HC$ respectively;
and because $EFG$ is a circle with center $C$,
$CG$ is equal to $CE$;
therefore the respective bases $CG$ and $CE$ of $GHC$ and $EHC$ are equal and the two
sides $GH,HC$ are equal to the two sides $EH,HC$ respectively;
therefore the angle $GHC$ is equal to the angle $EHC$ (I.8).

The straight line $CH$ set up on the straight line $AB$ has adjacent angles $GHC,EHC$, which are equal; therefore
$GHC,EHC$ are right angles, and $CH$ is perpendicular to $AB$ (Definition 10).
\end{proof}

I.13: ``If a straight line set up on a straight line make angles, it
will make either two right angles or angles equal to two right
angles.''

\begin{figure}
\begin{center}
\includegraphics[width=\textwidth]{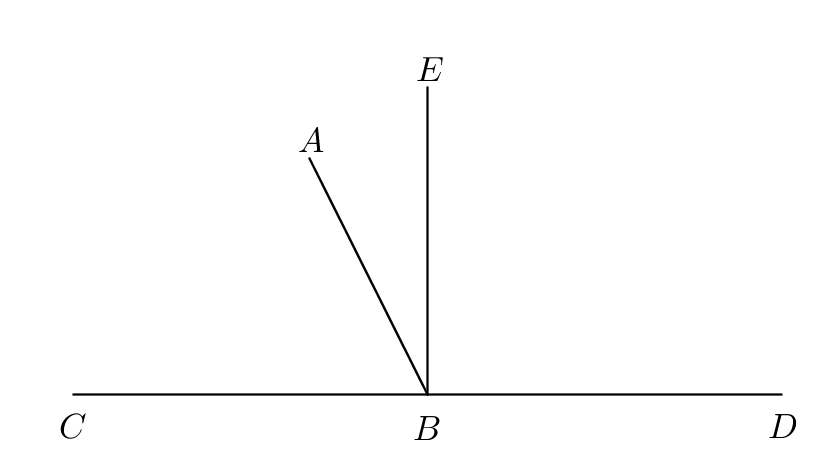}
\end{center}
\caption{I.13: P 23v, F 4v, B 14r, V 20}
\label{I13}
\end{figure}

\begin{proof}
Let a straight line $AB$ be set up on a straight line $CD$ make angles $CBA,ABD$.
If $CBA$ is equal to $ABD$ then $CBA,ABD$ are right angles (Definition 10).
If not, let $BE$ be drawn from $B$ at right angles to $CD$ (I.11); therefore
the angles $CBE, EBD$ are two right angles. 
The angle $CBE$ is equal to the two angles $CBA,ABE$.
Then the angles $CBE,EBD$ are equal to the three angles $CBA,ABE,EBD$ (Common Notion 2).
The angle $DBA$ is equal to the two angles $DBE,EBA$; therefore
the angles $DBA,ABC$ are equal to the three angles $DBE,EBA,ABC$ (Common Notion 2).
But the angles $CBE,EBD$ are also equal to these three angles; therefore the angles
$DBA,ABC$ are equal to the angles $CBE,EBD$ (Common Notion 1).
But the angles $CBE,EBD$ are two right angles; therefore
the angles $DBA,ABC$ are equal to two right angles.
\end{proof}

Proclus 292--293 \cite[pp.~228--229]{proclus}:

\begin{quote}
But what
does he intend when he adds that it makes ``either two right
angles or angles equal to two right angles?'' For when it
makes two right angles, it makes angles equal to two right
angles, since all right angles are equal to one another. Is it
not that the one expression denotes an attribute common to
both equal and unequal angles, the other a property of
equal angles only? Whenever both a general and a special
attribute can be affirmed truly of something, we are
accustomed to indicate its character by the special attribute;
but whenever we cannot hit upon this, we are satisfied with
the general character for the clarification of the things under consideration.
\end{quote}

I.14: ``If with any straight line, and at a point on it, two straight lines
not lying on the same side make the adjacent angles equal to
two right angles, the two straight lines will be in a straight line with one another.''

I.15: ``If two straight lines cut one another, they make the vertical
angles equal to one another.''

\begin{figure}
\begin{center}
\includegraphics[width=\textwidth]{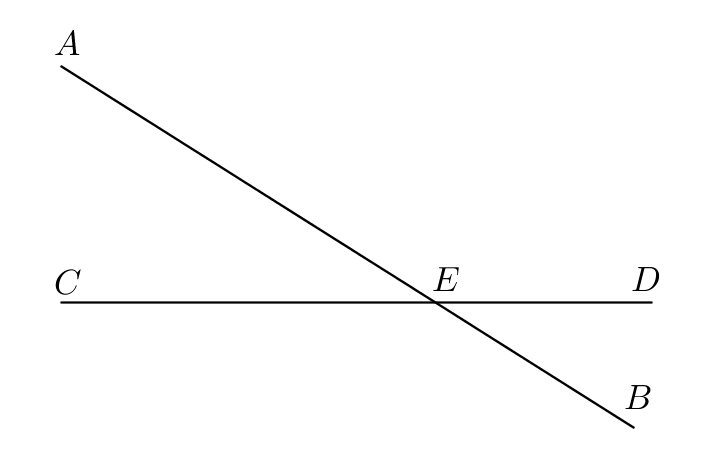}
\end{center}
\caption{I.15: F 5r, B 15r, V 21}
\label{I15}
\end{figure}

\begin{proof}
The straight line $AE$ stands on the straight line $CD$ and makes the angles
$CEA,AED$; thus the angles $CEA,AED$ are equal to two right angles (I.13).
And the straight line $DE$ stands on the straight line $AB$ and makes the angles
$AED,DEB$; thus the angles $AED,DEB$ are equal to two right angles (I.13).
Therefore the angles $CEA,AED$ are equal to the angles $AED,DEB$ (Postulate 4 and Common Notion 1).
Let the angle $AED$ be subtracted from $CEA,AED$ and $AED,DEB$; then the remainders
$CEA$ is equal to $DEB$ (Common Notion 3).
It can be proved in the same way that  $CEB$ is equal to $DEA$.
\end{proof}

Proclus 298 \cite[p.~233]{proclus}:

\begin{quote}
Vertical angles are different from adjacent angles, we say,
in that they arise from the intersection of two straight lines,
whereas adjacent angles are produced when one only of the
two straight lines is divided by the other. That is, if a straight
line, itself undivided, cuts the other with its extremity and
makes two angles, we call these angles adjacent; but if two
straight lines cut each other, they make vertical angles. We
call them so because their vertices come together at the same
point; and their vertices are the points at which the lines
converging make the angles.
\end{quote}

I.16: ``In any triangle, if one of the sides be produced, the exterior
angle is greater than either of the interior and opposite angles.''

I.20: ``In any triangle two sides taken together in any manner
are greater than the remaining one.''

I.22: ``Out of three straight lines, which are equal to three given
straight lines, to construct a triangle: thus it is necessary that
two of the straight lines taken together in any manner should
be greater than the remaining one.''

\begin{figure}
\begin{center}
\includegraphics[width=\textwidth]{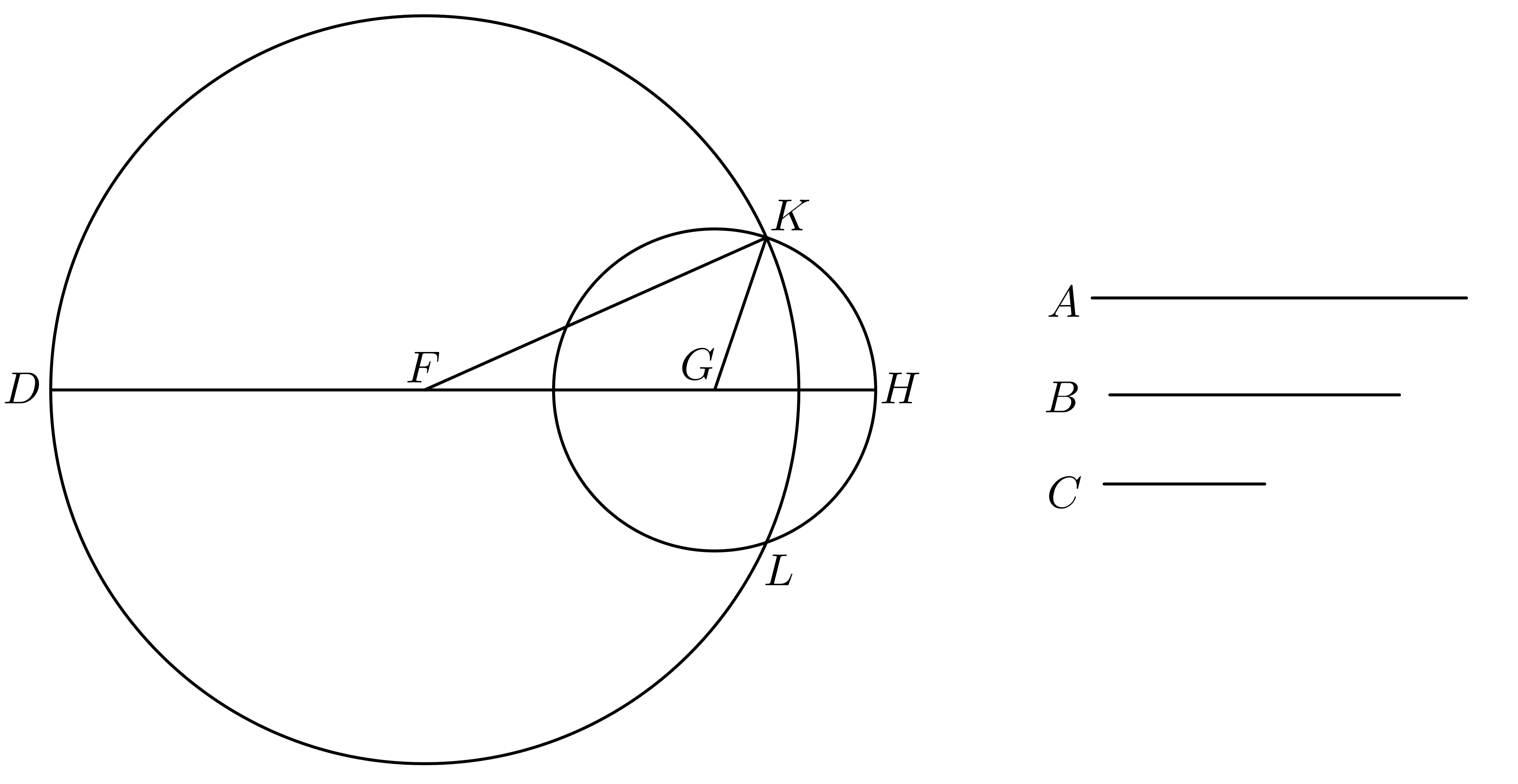}
\end{center}
\caption{I.22: P 25r, F 6v, B 18r, V 26}
\label{I22}
\end{figure}

\begin{proof}
Let the three given straight lines be $A,B,C$, with $A,B$ greater than $C$, and $A,C$ greater than 
$B$, and $B,C$ greater than $A$. Let $DE$ be a straight line that terminates at $D$ and is infinite in the direction
$E$. Let $DF$ be cut off from $DE$ equal to $A$, let $FG$ be cut off from $FE$ equal to $B$, and let
$GH$ be cut off from $GE$ equal to $C$ (I.3).
With center $F$ and distance $FD$ let the circle $DKL$ be described and
with center $G$ and distance $GH$ let the circle $KLH$ be described (Postulate 3), with
$K$ and $L$ the points at which the two circles intersect.
Let $KF,KG$ be joined (Postulate 1).

Because $F$ is the center of the circle $DKL$, $FD$ is equal to $FK$; but $FD$ is equal to $A$, so
$KF$ is equal to $A$ (Common Notion 1).
Because $G$ is the center of the circle $LKH$, 
$GH$ is equal to $GK$; but $GH$ is equal to $C$, so
$KG$ is equal to $C$ (Common Notion 1).
But $FG$ is equal to $B$, so
the three straight lines $KF,FG,GK$ are equal respectively to the three straight lines
$A,B,C$. 
\end{proof}

I.23: ``On a given straight line and at a point on it to construct a 
rectilineal angle equal to a given rectilineal angle.''

\begin{figure}
\begin{center}
\includegraphics[width=\textwidth]{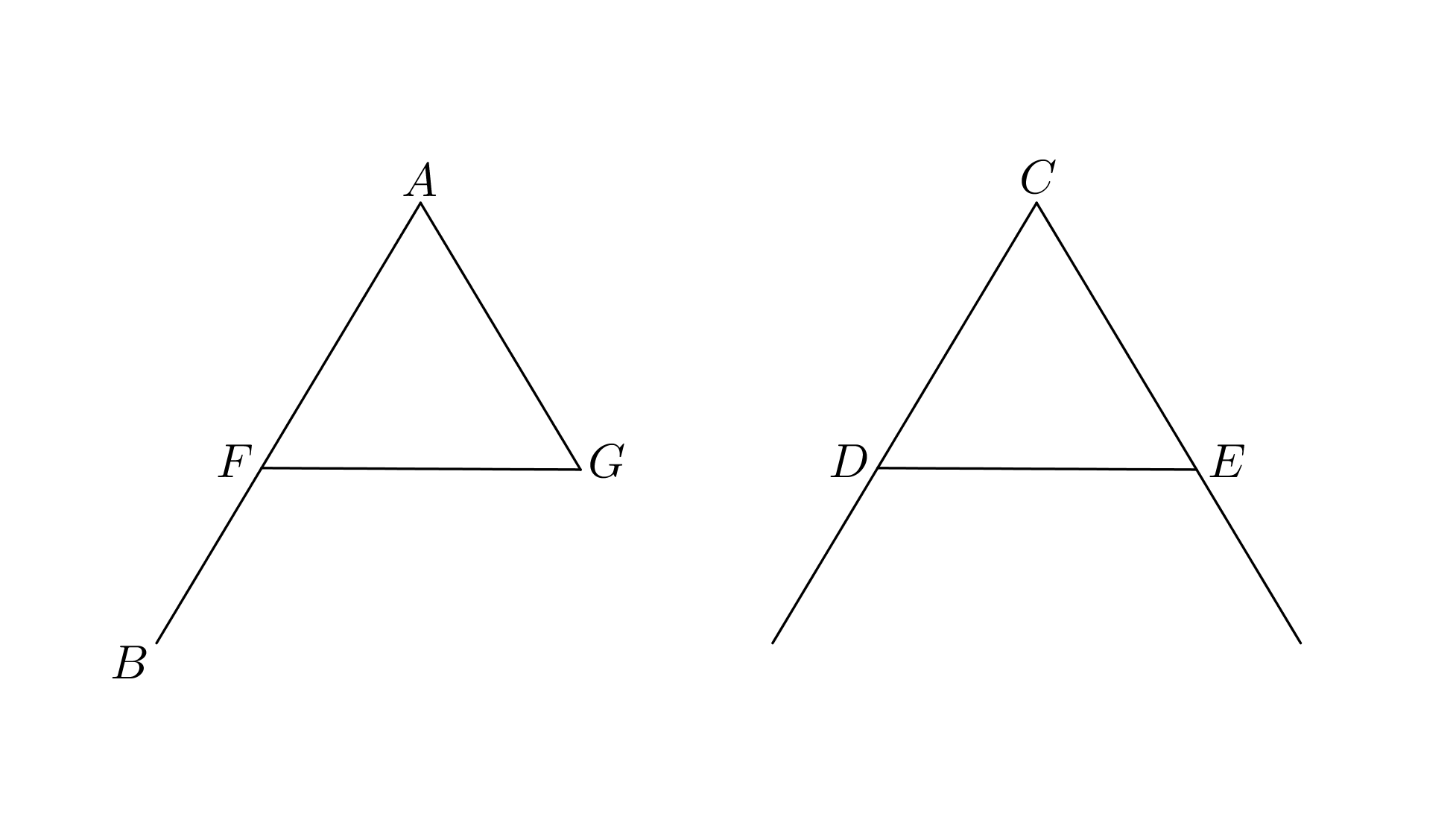}
\end{center}
\caption{I.23: P 25v, F 7r, B 18v, V 27}
\label{I23}
\end{figure}

\begin{proof}
{\em ekthesis}: Let the given straight line be $AB$, let the  point on the line be $A$,  and let the given angle be $DCE$.

{\em diorismos}: Thus it is required to construct on the line $AB$ at the point $A$ an angle equal to the angle $DCE$.

{\em kataskeu\={e}}: On the straight lines $CD,CE$ respectively let points taken by chance be $D,E$. 
Let $DE$ be joined. From the three straight lines $CD,DE,CE$ let the triangle $AFG$ be constructed such that $CD$ is equal to $AF$,
$CE$ to $AG$, and $DE$ to $FG$ (I.22).\footnote{What is used here is more than what I.22 provides.}

{\em apodeixis}: Since the two sides $DC,CE$ are respectively equal to the two sides $FA,AG$ and the base $DE$ is equal to the base
$FG$, the angle contained by the straight lines $DC,CE$ is equal to the angle contained by the straight lines
$FA,AG$ (I.8). That is, the angle $DCE$ is equal to the angle $FAG$.

{\em sumperasma}: On the given straight line $AB$ at the point $A$, the angle $FAG$ has been constructed that is equal to the given
angle $DCE$.
\end{proof}

Scholia for I.23 \cite[pp.~161--162]{euclidisV}.

The figures for I.23 in
Al-Nayrizi \cite[pp.~146--147]{alnayriziI} 
and Adelard of Bath \cite[pp.~49--50]{adelardI} are the same as I.23 in Heiberg.

Proclus 334--335 \cite[pp.~261--262]{proclus} gives the following proof of I.23.

\begin{figure}
\begin{center}
\includegraphics[width=\textwidth]{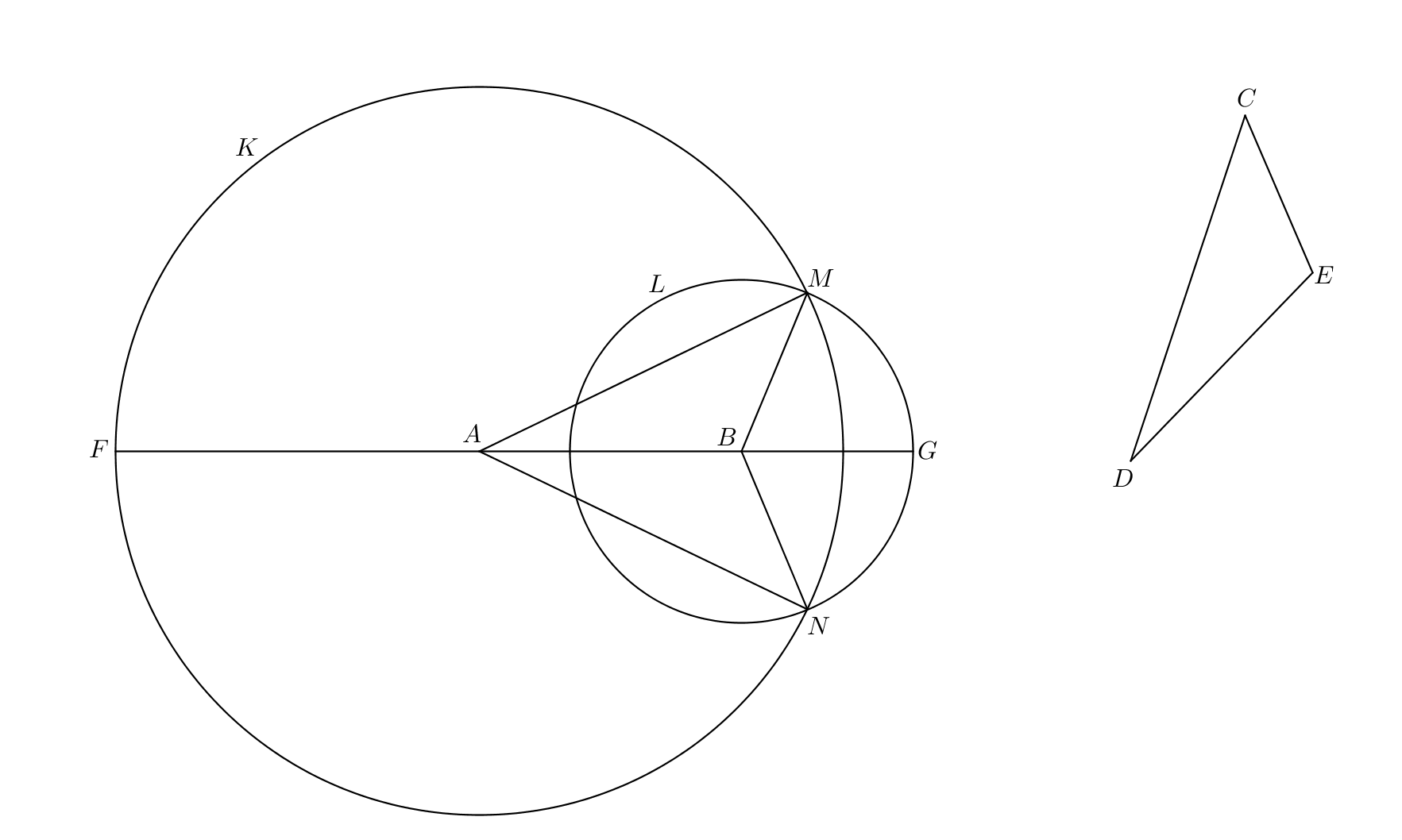}
\end{center}
\caption{Proclus, I.23}
\label{proclusI23}
\end{figure}

\begin{proof}
Let $AB$ be a given straight line, $A$ the given point on it, and $CDE$ the given rectilineal angle, with
$DE$ equal to $AB$.
Let $CE$ be joined and let $AB$ be produced in both directions to $F$ and $G$. Let
$FA$ be cut off equal to $CD$ and let $BG$ be cut off equal to $CE$ (I.3).
With $A$ as center and distance $FA$ let the circle $K$ be described and
with $B$ as center and distanc $BG$ let the circle $L$ be described (Postulate 3).
These circles intersect at the points $M,N$.
Let $MA,MB$ be joined and let $NA,NB$ be joined (Postulate 1).
Because $FA$ is equal to each of $AM,AN$ and $FA$ is equal $CD$, each of $AM,AN$ is equal to $CD$, and likewise
each of $BM,BN$ is equal to $CE$ (Common Notion 1).
But $AB$ is equal to $DE$, so the two lines $AB,AM$ are equal respectively to $DE,DC$ and the base
$BM$ is equal to $CE$, hence the angle
$MAB$ is equal to the angle $CDE$ (I.8); it can likewise be proved that the angle
$NAB$ is equal to the angle $CDE$.
\end{proof}

Johannes de Tinemue \cite[pp.~52--53]{adelardIII} gives a proof of I.23
that does not invoke I.22.

\begin{quote}
Data recta linea supra terminum eius cuilibet angulo proposito
equum angulum designare.
\end{quote}

\begin{figure}
\begin{center}
\includegraphics[width=\textwidth]{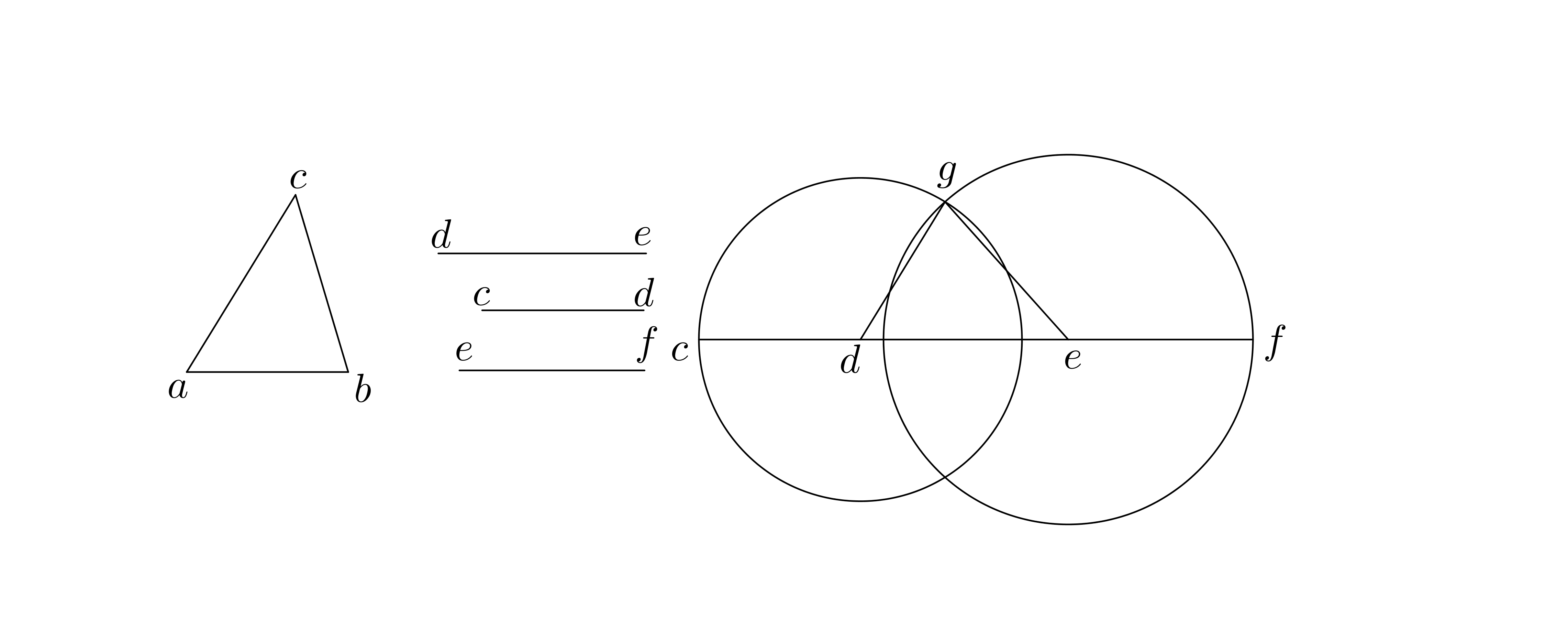}
\end{center}
\caption{Johannes de Tinemue, I.23}
\label{tinemueI23}
\end{figure}

\begin{quote}
Dispositio. Clausis itaque $b$, $c$ terminis interpositione $b c$ linee adequetur
$d e$ ad $a c$. Duabus lineis altrinsecus eidem scilicet $d e$ directe protractis
$c d$ equali $a b$ et $e f$ equali $b c$, deinde secundum premissam ex tribus lineis
ipsis $c d$, $d e$, $e f$ equalibus designetur $d e g$ triangulus.

Ratiocinatio. Age. Si memineris priorum, istorum triangulorum $a b c$,
$d e g$, $a b$, $d g$ mediante $c d$ et $a c$, $d e$ et $b c$, $e g$ mediante $e f$ sunt equalia.
Ergo secundum 8\textsuperscript{am} $d$ et $a$ anguli sunt equales. Sicque super terminum $d$
designatus est angulus equalis $a$. Quod proposuimus.
\end{quote}

Albert the Great \cite[pp.~82--84]{albertus} gives a proof of I.23 that does not invoke I.22.

\begin{figure}
\begin{center}
\includegraphics[width=\textwidth]{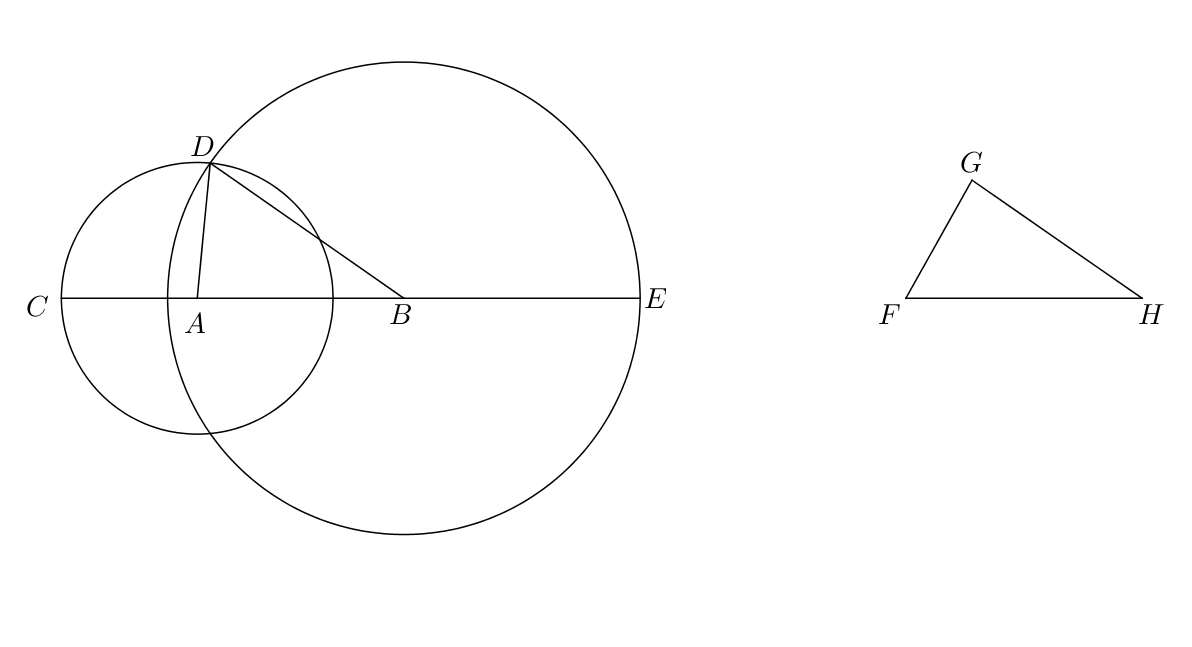}
\end{center}
\caption{Albert the Great, I.23, intersecting circles}
\label{albertI23a}
\end{figure}

\begin{figure}
\begin{center}
\includegraphics[width=\textwidth]{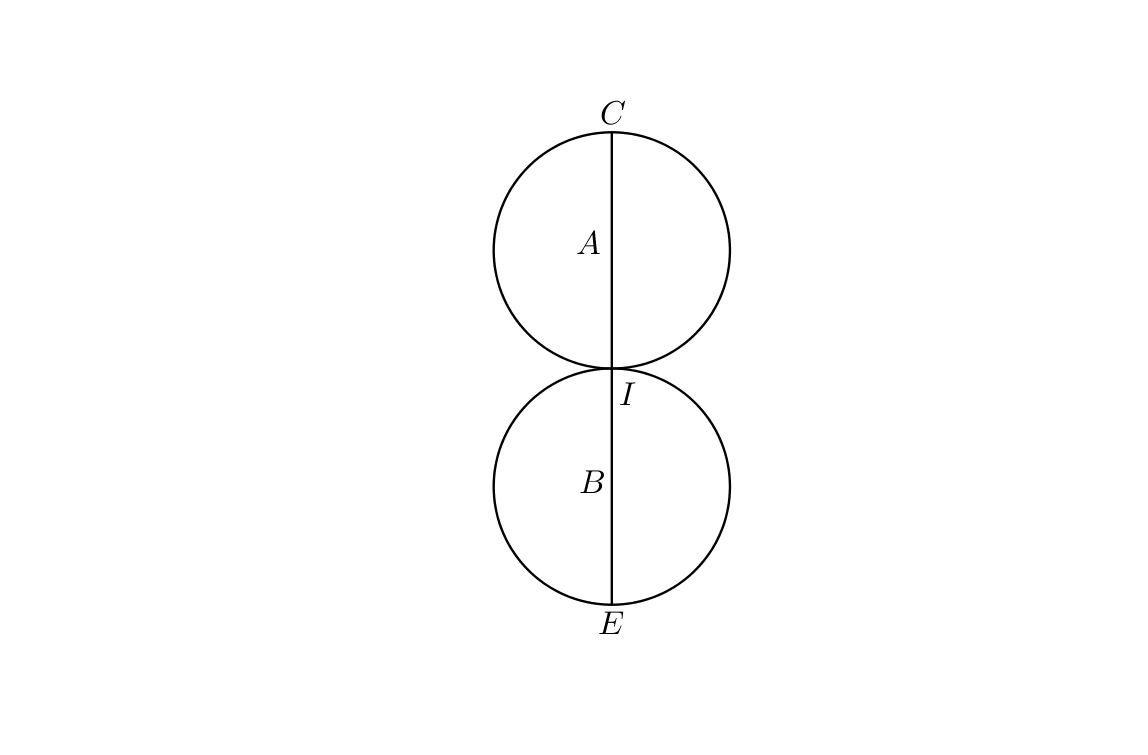}
\end{center}
\caption{Albert the Great, I.23, touching circles}
\label{albertI23b}
\end{figure}

\begin{proof} 
Let the given straight line be $AB$ and let the given angle be $FGH$.
If $AB$ is longer than $GH$ then apply I.3 to cut off $AB$ equal to $GH$.
If $AB$ is shorter than $GH$ then apply I.2 to extend $AB$ to be equal to $GH$.
Extend $AB$ on the side of $A$ to $C$ so that $AC$ is equal to $FG$.
Extend $AB$ on the side of $B$ to $E$ so that $BE$ is equal to $FH$.
Draw the circle with center $A$ and radius $AC$ and the circle with center $B$ and radius $BE$.
These circles either intersect or they do not.

If they intersect let $D$ be one of the points of intersection. 
$AC$ is equal to $FG$ and $AC$ is equal to $AD$ so $AD$ is equal to $FG$. 
Likewise, $BE$ is equal to $FH$ and $BE$ is equal to $BD$ so $BD$ is equal to $FH$.
Finally, $AB$ was made equal to $GH$. Thus the two triangles $DAB,FGH$ have the two
sides  $DA,AB$ are respectively equal to the two sides $FG,GH$ and the 
base $DB$ is equal to the base $FH$, so by I.8 the angle contained by the sides $DA,AB$ is equal to the angle
contained by $FG,GH$, that is, the angle $DAB$ is equal to the given angle $FGH$.

If the circles do not intersect then either they touch each other on line $AB$ or they do not touch each other.
If they touch each other, let $I$ be the point of contact.
Because $I$ lies on the circle $AC$, $AC$ is equal to $AI$. Because $I$ lies on the circle
$BE$, $EB$ is equal to $BI$.
Because the circles touch each other at $I$, $AIB$ is a straight line, so $AB$ is equal to $AC,EB$.
But $AB$ is equal to $GH$, $AC$ is equal to $FG$, and $EB$ is equal to $FH$. So $GH$ is equal to $FG,FH$, which according
to I.20 is absurd. 

If the circles do not touch each other, then by the above reasoning we get that $GH$ is greater than $FG,FH$, which according
to I.20 is absurd.
\end{proof}

Campanus \cite[p.~74]{campanusI}:

\begin{quote}
Data recta linea super terminum eius cuilibet angulo proposito
equum angulum designare.
\end{quote}

\begin{figure}
\begin{center}
\includegraphics[width=\textwidth]{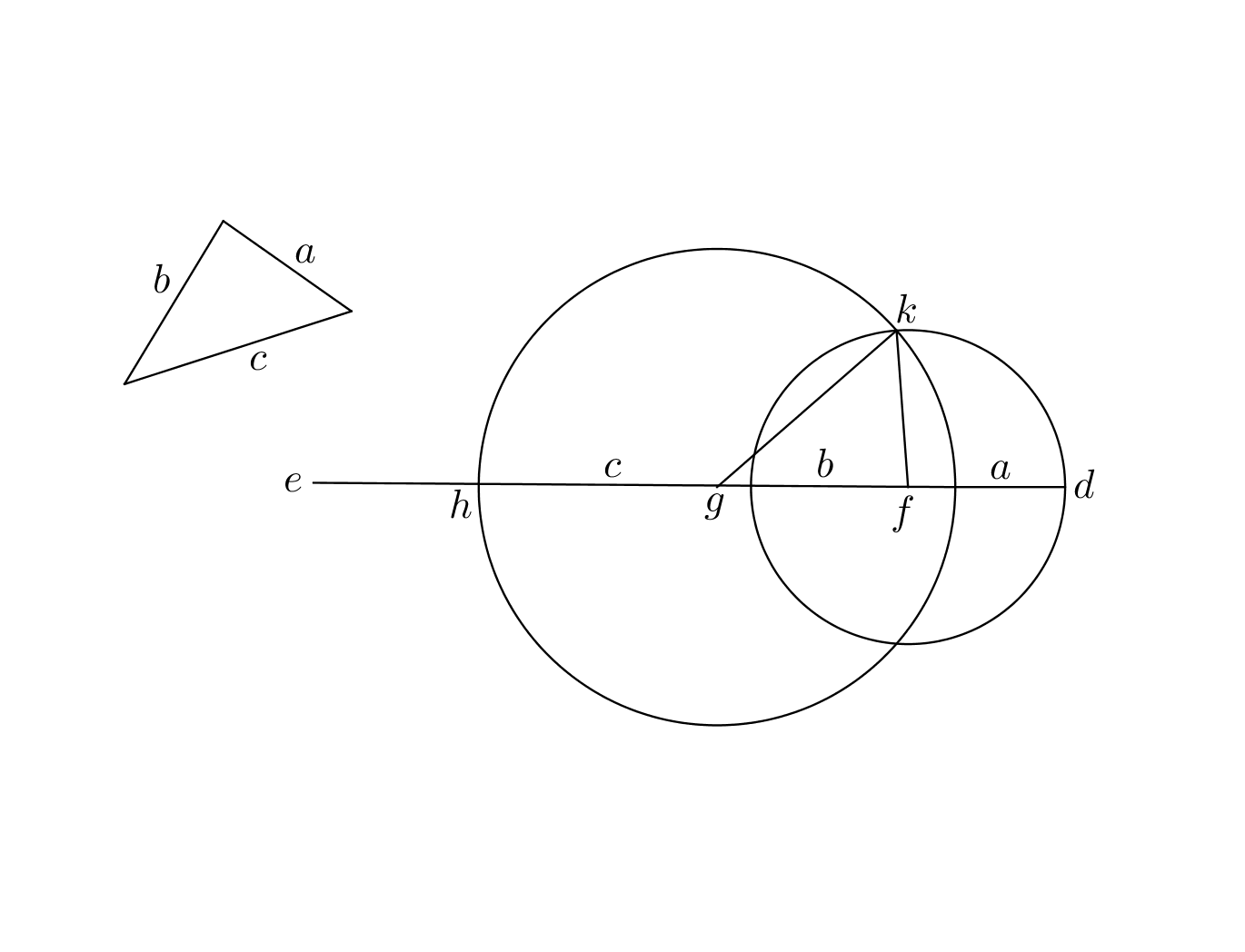}
\end{center}
\caption{Campanus, I.23}
\label{campanusI23}
\end{figure}

\begin{quote}
Sit data linea $fe$ que est in superiori figura et sint linee $b,a$ continentes
angulum datum cui subtendam basim $c$. Supra punctum $f$ linee $ef$
iubemur facere equalem angulum angulo dato ad lineam $ef$. Adiungo $fd$
equalem linee $a$ et ex $fe$ sumo $fg$ equalem $b$ et ex $ge$ sumo $gh$ equalem $c$
et super puncta $f$ et $g$ describo duos circulos $dk$ et $kh$ secundum quantitatem
duarum linearum $fd$ et $gh$ intersecantes se in puncto $k$ sicut docuit
precedens. 
Ductisque lineis $kf$ et $kg$ erunt duo latera $kf$ et $fg$ trianguli
$kfg$ equalia duobus lateribus $a$ et $b$ trianguli $abc$ et basis $gk$ equalis basi
$c$, ergo per 8 angulus $kfg$ equalis erit angulo contento ab $a$ et a $b$. Quod
est propositum.
\end{quote}

\begin{proof}
Let the given line be $fe$ and let the given angle be contained by the sides $b,a$ and subtended by the base
$c$. Extend $fd$ to be equal to $a$; from $fe$ take $fg$ equal to $b$;
 from $ge$ take $gh$ equal to $c$; and at the points $f$ and $g$ describe two circles 
 $dk$ and $kh$. 
It follows from I.22 that 
these circles intersect, and let one of the points of intersection be $k$. 
 Join the lines $kf$ and $kg$. 
 Then the two sides $kf$ and $fg$ of the triangle
 $kfg$ are
 equal to the two sides $a$ and $b$ of the triangle $abc$ respectively,
 and the base $gk$ is equal to the base $c$, so by I.8 the angle $kfg$ is equal to the angle contained
 by $a$ and $b$.  
\end{proof}

Peletarius, {\em Euclidis elementa geometrica demonstrationum libri sex}, 1557, p.~42.

Commandinus, {\em Euclidis Elementorum libri XV}, 1572, f.~17r, in his commentary on I.23, gives the following
construction.

\begin{figure}
\begin{center}
\includegraphics[width=\textwidth]{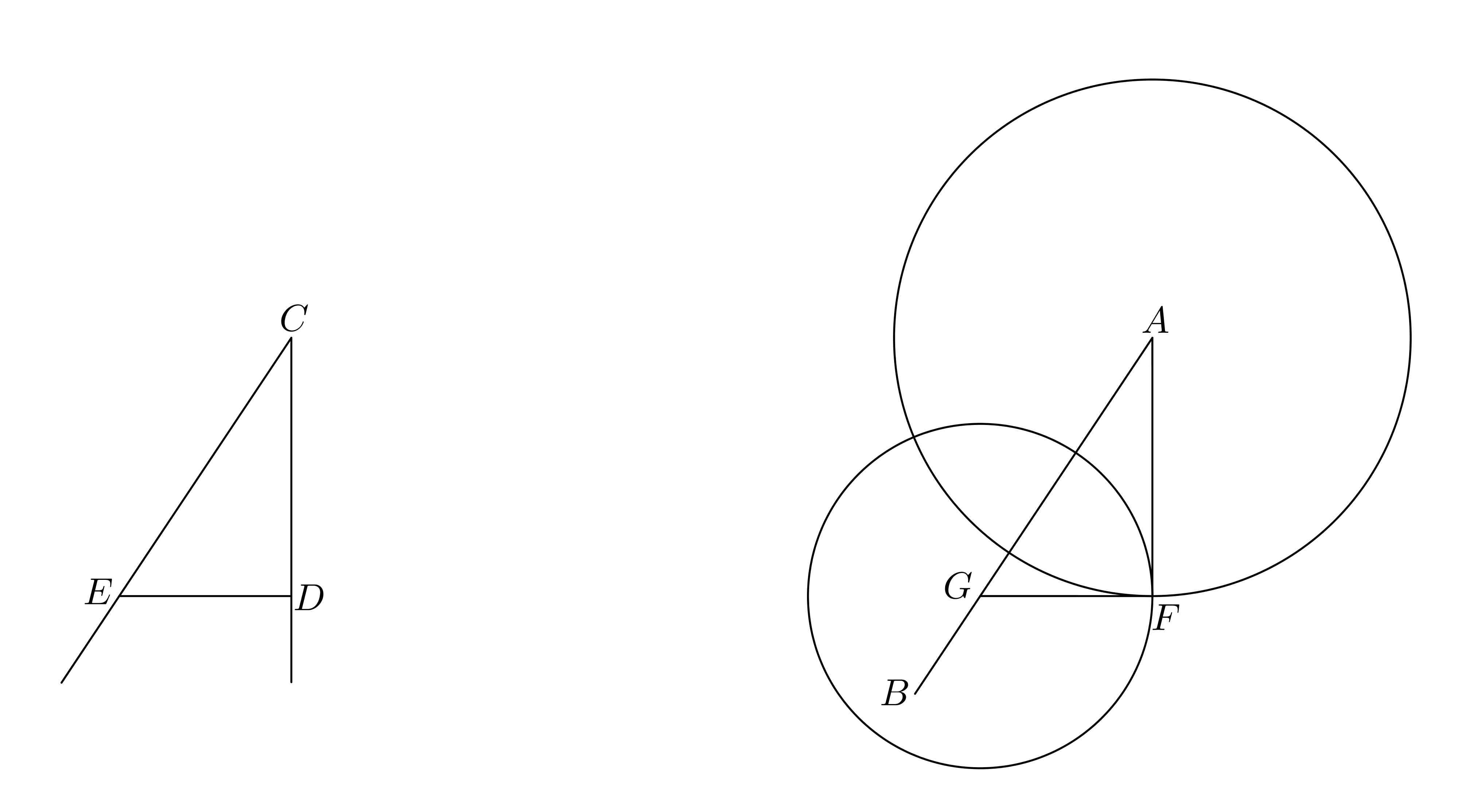}
\end{center}
\caption{Commandinus, I.23}
\label{commandinusI23}
\end{figure}

\begin{proof}
Let $AB$ be the given line and let $A$ be the given point. 
Let $DCE$ be the given angle and let $DE$ be joined.
Cut $AG$ from $AB$ equal to $CE$. 
With center $A$, describe a circle with radius $CD$ and with center $G$ describe a circle
with radius $ED$.
The circles intersect at a point $F$. Join $AF,FG$. Then
$FAG$ is equal to $DCE$. 
\end{proof}

Clavius, {\em Euclidis Elementorum libri XV}, {\em Opera mathematica}, pp.~44--45, 
I.23.

\begin{figure}
\begin{center}
\includegraphics[width=\textwidth]{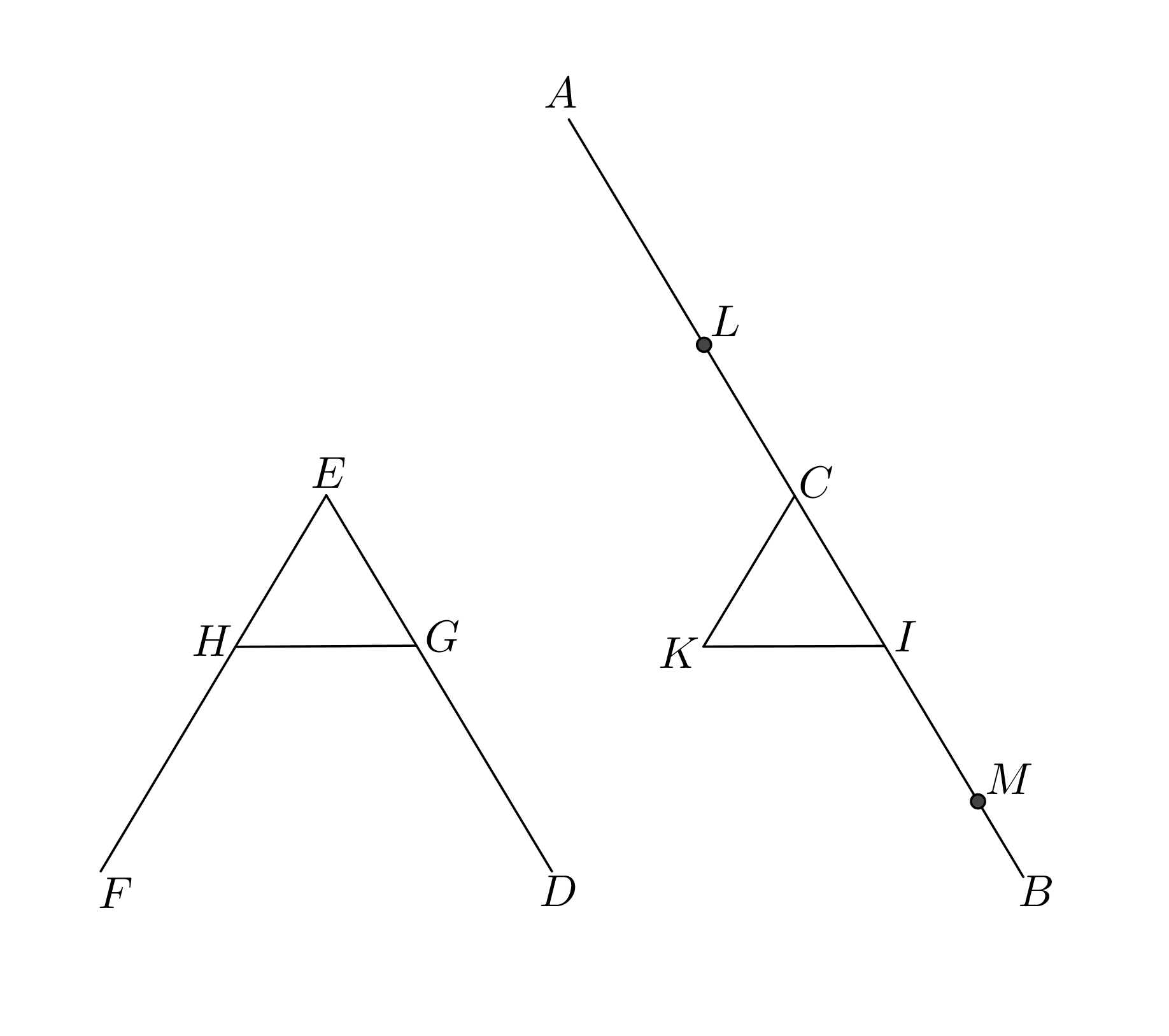}
\end{center}
\caption{Clavius, I.23}
\label{claviusI23}
\end{figure}

\begin{proof}
Let the given line be $AB$, let the given point on the line be $C$, and let the given
angle be $DEF$. 

Take $CI$ equal to $EG$, take $CL$ equal to $EH$, and take $IM$ equal to $GH$.
Describe a circle with center $C$ and radius $CL$ and describe a circle with center $I$ and radius $IM$. These circles intersect at $K$.
\end{proof}

I.26: ``If two triangles have the two angles equal to two angles
respectively, and one side equal to one side, namely, either the
side adjoining the equal angles, or that subtending one of the
equal angles, they will also have the remaining sides equal to
the remaining sides and the remaining angle to the remaining
angle.''

\begin{figure}
\begin{center}
\includegraphics[width=\textwidth]{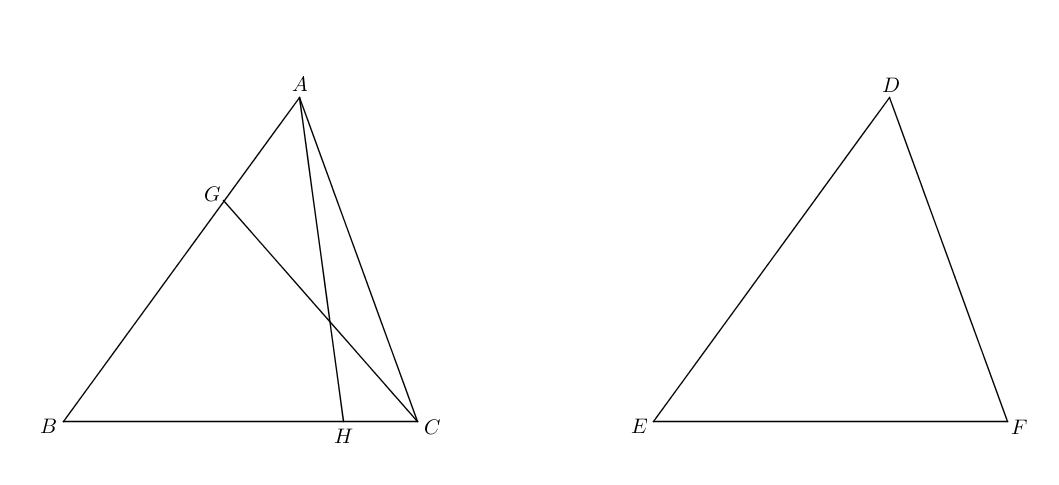}
\end{center}
\caption{I.26: P 27r, F 8r, B 20r, V 30}
\label{I26}
\end{figure}

\begin{proof}
Let $ABC,DEF$ be two triangles having 
the two angles $ABC,BCA$ equal respectively to the two angles
$DEF,EFD$. First, let the sides adjoining the equal angles be equal, namely
let $BC$ be equal to $EF$.

If $AB$ is unequal to $DE$, one is greater and let $AB$ be greater than $DE$.
Let $BG$ be made equal to $DE$ (I.3), and let $GC$ be joined (Postulate 1).
Because $BG$ is equal to $DE$ and $BC$ is equal to $EF$,
and the angle $GBC$ is equal to the angle $DEF$, the
base $GC$ is equal to the base $DF$ and the angles
$GCB,BGC$ are equal respectively to the angles
$DFE,EDF$ (I.4). But by hypothesis the angle $DFE$ is equal to the angle
$BCA$; thus the angle $BCG$ is equal to the angle $BCA$ (Common Notion 1),
namely the less is equal to the greater, which is impossible. Therefore
$AB$ is not unequal to $DE$, so $AB$ is equal to $DE$.

But by hypothesis $BC$ is equal to $EF$, so the sides
$AB,BC$ are equal respectively to the sides $DE,EF$; and
the angle $ABC$ is equal to the angle $DEF$, thus 
the base $AC$ is equal to the base $DF$ and the remaining angle
$BAC$ is equal to the remaining angle $EDF$ (I.4).

Second, let sides subtending equal angles be equal, with $AB$ equal to $DE$.
If $BC$ is unequal to $EF$, one is greater and let $BC$ be greater than $EF$.
Let $BH$ be made equal to $EF$ (I.3), and let $AH$ be joined (Postulate 1).
Then the two sides $AB,BH$ are equal respectively to
$DE,EF$, and by hypothesis the angle $ABH$ is equal to the angle $DEF$;
thus the base $AH$ is equal to the base $DF$, the angle
$BHA$ is equal to the angle $EFD$, and the angle $BAH$ is equal to the angle
$EDF$ (I.4).
But the angle $EFD$ is equal by hypothesis to the angle $BCA$, so the angle
$BHA$ is equal to the angle $BCA$ (Common Notion 1).
Therefore for the triangle $AHC$, the exterior angle $BHA$ is equal
to the interior and opposite angle $BCA$, which is impossible (I.16).
Therefore $BC$ is not unequal to $EF$, so $BC$ is equal to $EF$.

But by hypothesis $AB$ is equal to $DE$. Thus the two sides $AB,BC$ are equal respectively to
the two sides $DE,EF$, and the angle $ABC$ is equal to the angle $DEF$; thus
the base $AC$ is equal to the base $DF$ and the 
remaining angle $BAC$ is equal to the remaining angle $EDF$ (I.4).
\end{proof}

I.27: ``If a straight line falling on two straight lines make the
alternate angles equal to one another, the straight lines will be
parallel to one another.''

\begin{figure}
\begin{center}
\includegraphics{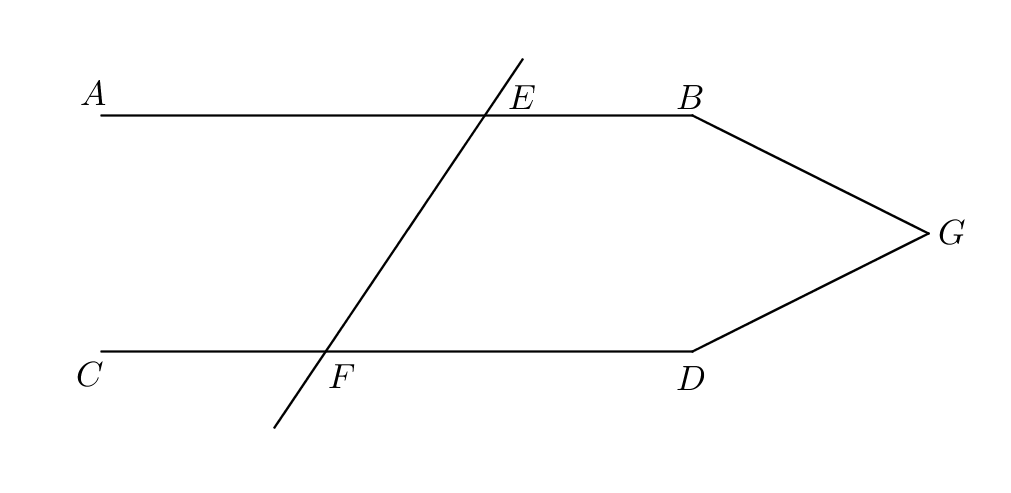}
\end{center}
\caption{I.27: P 27v, F 8v, B 21v, V 31}
\label{I27}
\end{figure}

\begin{proof}
Let the straight line $EF$ falling on the two straight lines $AB,CD$ make the alternate angles
$AEF,EFD$ equal to one another. If $AB$ is not parallel to $CD$, then $AB,CD$ when produced will
meet either in the direction of $B,D$ or in the direction of $A,C$ (Definition 23); let them be produced and meet in the direction
of $B,D$ at $G$.
For the triangle $GEF$, by hypothesis the exterior angle $AEF$ is equal to the interior and opposite angle
$EFG$, which is impossible (I.16).
Therefore $AB,CD$ when produced will not meet in the direction of $B,D$. It can likewise be proved
that $AB,CD$ when produced will not meet in the direction of $A,C$.
Therefore $AB$ is parallel to $CD$.
\end{proof}

Proclus 357 \cite[p.~278]{proclus}:

\begin{quote}
Angles that are produced in different directions and are not adjacent
to one another, but separated by the intersecting line, both of
them within the parallels but differing in that one lies above
and the other below, he calls ``alternate'' angles.
\end{quote}

I.28: ``If a straight line falling on two straight lines make the
exterior angle equal to the interior and opposite angle on the
same side, or the interior angles on the same side equal to two
right angles, the straight lines will be parallel to one another.''

\begin{figure}
\begin{center}
\includegraphics{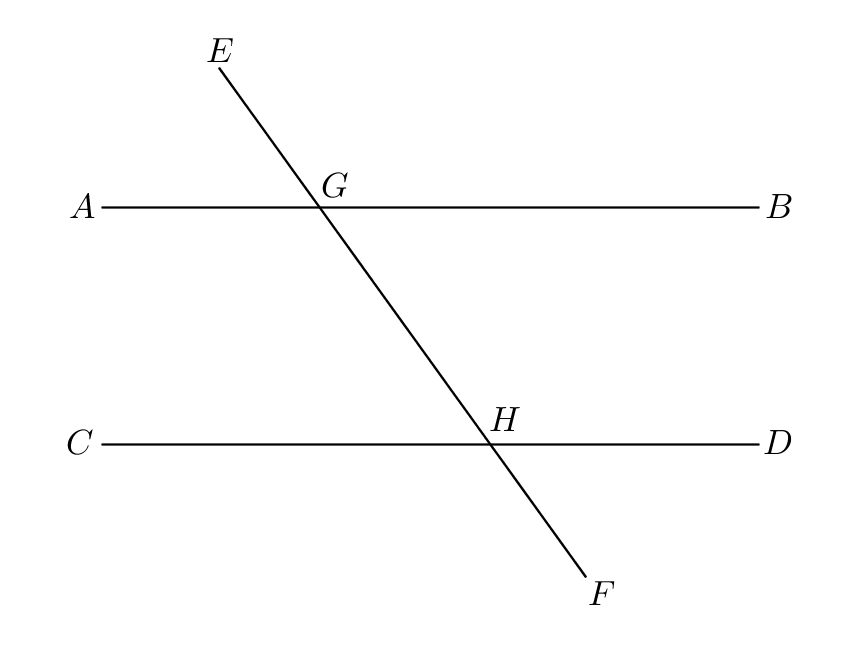}
\end{center}
\caption{I.28: P 28r, F 9r, B 21v, V 32}
\label{I28}
\end{figure}

\begin{proof}
First, let the straight line $EF$ falling on the two straight lines $AB,CD$ make the exterior angle
$EGB$ equal to the interior and opposite angle $GHD$. 
Since the straight lines $EF,AB$ cut one another, the vertical angles $EGB,AGH$ are equal (I.15).
Then since by hypothesis the angle $EGB$ is equal to the angle $GHD$, the angle $AGH$ is 
equal to the angle $GHD$ (Common Notion 1);  the angles $AGH,GHD$ are alternate, so
$AB$ is parallel to $CD$ (I.27).

Second, let the straight line $EF$ falling on the two straight lines $AB,CD$ make the interior angles $BGH,GHD$
equal to two right angles. 
The straight line $GH$ set up on the straight line $AB$ makes the angles $AGH,BGH$, which thus are equal to
two right angles (I.13).
Since the angles $BGH,GHD$ are equal to two right angles and the angles $AGH,BGH$ are equal to two right angles,
the angles $AGH,BGH$ are equal to the angles $BGH,GHD$ (Postulate 4 and Common Notion 1).
Let the angle $BGH$ be subtracted from each;
then the angle $AGH$ is equal to the angle $GHD$ (Common Notion 3).
The alternate angles $AGH,GHD$ are equal, so the straight line $AB$ is parallel to the straight line $CD$ (I.27). 
\end{proof}

I.29: ``A straight line falling on parallel straight lines makes
the alternate angles equal to one another, the exterior angle
equal to the interior and opposite angle, and the interior angles
on the same side equal to two right angles.''

\begin{figure}
\begin{center}
\includegraphics{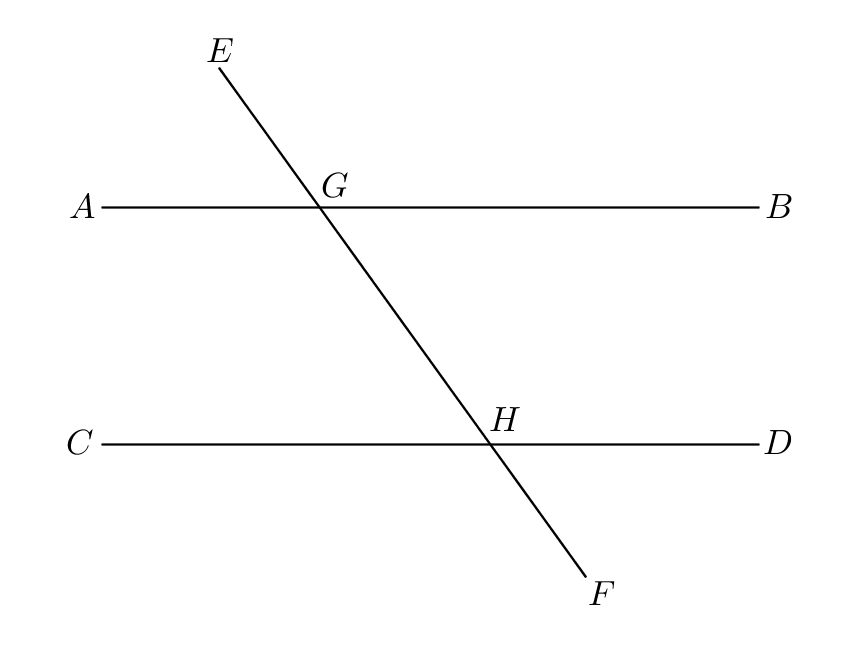}
\end{center}
\caption{I.29: P 28r, F 9r, B 22r, V 33}
\label{I29}
\end{figure}

\begin{proof}
Let the straight line $EF$ fall on the parallel straight lines $AB,CD$. 
If the alternate angles $AGH,GHD$ are unequal then one is greater; let the angle $AGH$ be greater
than the angle $GHD$.
Let the angle $BGH$ be added to each; then the angles $AGH,BGH$ are greater than the
angles $BGH,GHD$. But the angles $AGH,BGH$ are equal to two right angles (I.13); therefore
the angles $BGH,GHD$ are less than two right angles. Thus
the straight line $EF$ falling on the two straight lines $AB,CD$ makes
the interior angles $BGH,GHD$ on the same side less than two right angles, so
the straight lines $AB,CD$ if produced indefinitely meet on 
the same side as are the interior angles $BGH,GHD$ (Postulate 5).
But the straight lines $AB,CD$ are by hypothesis parallel so do not meet;
therefore the angle $AGH$ is not unequal to the angle $GHD$, so the alternate angles
$AGH,GHD$ are equal.

Because the straight line $EF$ cuts the straight line $AB$, 
the vertical angles $AGH,EGB$ are equal (I.15). But the alternate angles
$AGH,GHD$ are equal, so
the exterior angle $EGB$  is equal to the interior and opposite angle $GHD$  (Common Notion 1).

Let the angle $BGH$ be added to each of $EGB,GHD$; then the angles $EGB,BGH$ are 
$BGH,GHD$ (Common Notion 2).
But the angles $EGB,BGH$ are equal to two right angles (I.13).
Therefore the interior angles on the same side $BGH,GHD$ are equal to two right angles (Common Notion 1).
\end{proof}

I.30: ``Straight lines parallel to the same straight line are also
parallel to one another.''

\begin{figure}
\begin{center}
\includegraphics{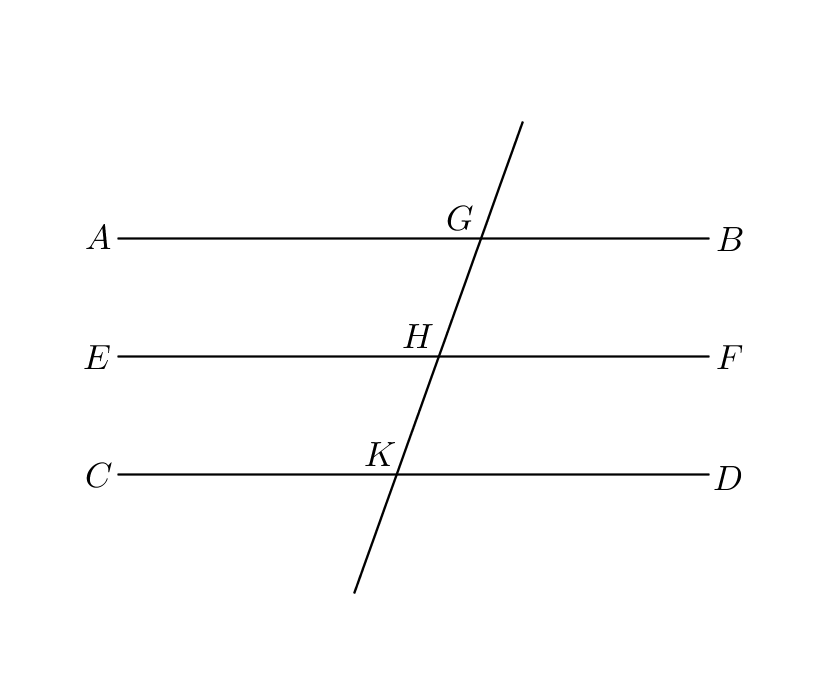}
\end{center}
\caption{I.30: P 28v, F 9v, B 23r, V 34}
\label{I30}
\end{figure}

\begin{proof}
Let each of the straight lines $AB,CD$ be parallel to $EF$. 
Let the straight line $GK$ fall upon these straight lines.
Since the straight line $GK$ has fallen on the parallel straight lines $AB,EF$,
the alternate angles $AGK,GHF$ are equal (I.29).
Likewise, since the straight line $GK$ has fallen on the parallel straight lines $EF,CD$,
the exterior angle $GHF$ is equal to the interior and opposite angle $GKD$ (I.29).
But the angle $AGK$ was proved to be equal to the angle $GHF$, so
the angle $AGK$ is equal to the angle $GKD$ (Common Notion 1).
The line $GK$ falling on the lines $AB,CD$ makes the alternate angles $AGK,GKD$, and because these are equal,
the lines $AB,CD$ are parallel (I.27).
\end{proof}

I.31: ``Through a given point to draw a straight line parallel to a 
given straight line.''

\begin{figure}
\begin{center}
\includegraphics{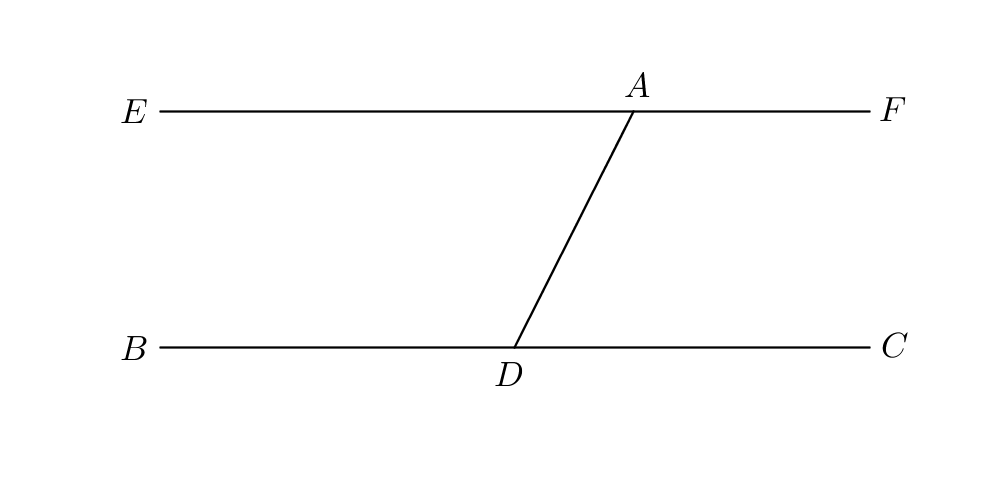}
\end{center}
\caption{I.31: P 28v, F 9v, B 23r, V 34}
\label{I31}
\end{figure}

\begin{proof}
Let $A$ be the given point and let $BC$ be the given straight line. Let $D$ be a point on $BC$ taken by chance
and let $AD$ be joined (Postulate 1).
On the straight line $DA$ at the point $A$ let the angle $DAE$ be constructed equal
to the angle $ADC$ (I.23).
Then let the straight line $AF$ be produced in a straight line with $EA$ (Postulate 2).
The straight line $AD$ falling on the two straight lines $BC,EF$ makes the alternate angles
$EAD,ADC$; but the angles $EAD,ADC$ have been proved to be equal; therefore
the lines $BC,EAF$ are parallel (I.27).
\end{proof}

I.32: ``In any triangle, if one of the sides be produced, the exterior
angle is equal to the two interior and opposite angles, and the
three interior angles of the triangle are equal to two right
angles.''

Proclus 381--382 \cite[pp.~300--301]{proclus}, on {\em Elements} I.32 (Proclus refers to {\em Timaeus} 53c):

\begin{quote}
We can now say that in 
every triangle the three angles are equal to two right angles. 
But we must find a method of discovering for all the other
rectilineal polygonal figures -- for four-angled, five-angled,
and all the succeeding many-sided figures -- how many right
angles their angles are equal to. First of all, we should know
that every rectilineal figure may be divided into triangles, for
the triangle is the source from which all things are constructed,
as Plato teaches us when he says, ``Every rectilineal plane
face is composed of triangles.'' Each rectilineal figure is
divisible into triangles two less in number than the number of
its sides: if it is a four-sided figure, it is divisible into two
triangles; if five-sided, into three; and if six-sided, into four.
For two triangles put together make at once a four-sided
figure, and this difference between the number of the
constituent triangles and the sides of the first figure
composed of triangles is characteristic of all succeeding figures.
Every many-sided figure, therefore, will have two more
sides than the triangles into which it can be resolved. Now
every triangle has been proved to have its angles equal to
two right angles. Therefore the number which is double the
number of the constituent triangles will give the number of
right angles to which the angles of a many-sided figure are
equal. Hence every four-sided figure has angles equal to four
right angles, for it is composed of two triangles; and every
five-sided figure, six right angles; and similarly for the rest.
\end{quote}

I.33: ``The straight lines joining equal and parallel straight
lines (at the extremities which are) in the same directions
(respectively) are themselves also equal and parallel.''

\begin{figure}
\begin{center}
\includegraphics{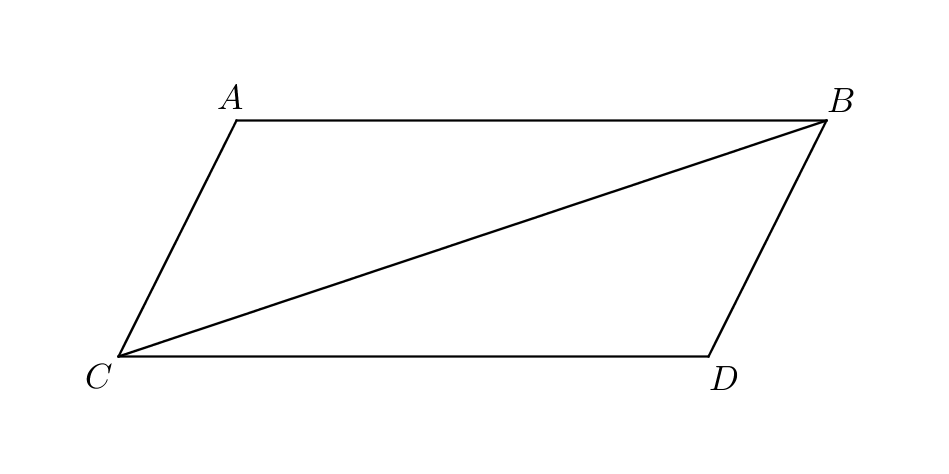}
\end{center}
\caption{I.33: P 30v, F 10r, B 24r, V 36}
\label{I33}
\end{figure}

\begin{proof}
Let $AB,CD$ be equal and parallel and let the straight lines $AC,BD$ joins the points on the same
sides. Let $BC$ be joined (Postulate 1). 
Since $BC$ falls on the parallel lines $AB,CD$, the alternate angles $ABC,BCD$ are equal
(I.29).
By hypothesis, $AB$ is equal to $CD$, and $BC$ is a common side of the triangles
$ABC,DCB$, so 
in the triangles $ABC,DCB$, the two sides $AB,BC$ are equal to the two
sides $CD,BC$, and the angles $ABC,BCD$ contained by the equal sides are equal; thus
the base $AC$ of the triangle $ABC$ is equal to the base $BD$ of the triangle $DCB$,
the angle $ACB$ is equal to the angle $CDB$, and the angle $BAC$ is equal to the angle $BDC$ (I.4).
Then the straight line $BC$ falling on the straight lines $AB,CD$ makes the alternate angles
$ACB,CBD$ equal; therefore $AC$ is parallel to $BD$ (I.27). And $AC$ was proved equal to $BD$, so
the lines $AC,BD$ are equal and parallel.
\end{proof}

I.34: ``In parallelogrammic areas the opposite sides and angles
are equal to one another, and the diameter bisects the areas.''

\begin{figure}
\begin{center}
\includegraphics{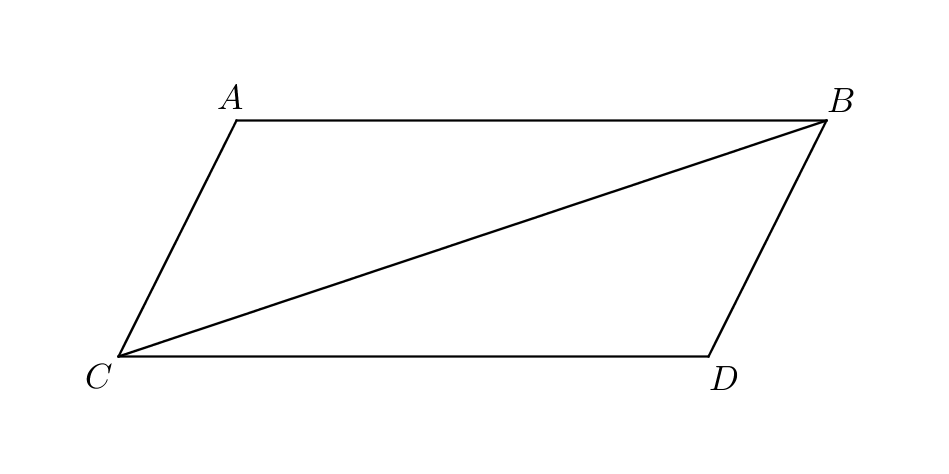}
\end{center}
\caption{I.34: P 32r, F 10r, B 24v, V 37}
\label{I34}
\end{figure}

\begin{proof}
Let $ACDB$ be a parallelogrammic area and let $BC$ be its diameter.
Since the straight line $BC$ falls on the parallel straight lines $AB,CD$,
the alternate angles $ABC,BCD$ are equal (I.29).
And since the straight line $BC$ falls on the parallel straight lines
$AC,BD$, the alternate angles $ACB,CBD$ are equal (I.29).
The two angles $ABC,BCA$ of the triangle $ABC$ are equal respectively to
the two angles $DCB,CBD$ of the triangle $DCB$, and the side $BC$ adjoining the equal angles
is common to both triangles; therefore
the side $AB$ of $ABC$ is equal to the side $CD$ of $DCB$, the
side $AC$ is equal to the side $BD$ of $DCB$, and the angle
$BAC$ of $ABC$ is equal to the angle $CDB$ of $DCB$ (I.26).
The angles $ABC,CBD$ are equal respectively to the angles
$BCD,ACB$;
and the whole angle $ABD$ is equal to the angles $ABC,CBD$, and the whole
angle $ACD$ is equal to the angles $ACB,BCD$; therefore 
the angle $ABD$ is equal to the angle $ACD$ (Common Notion 2).
But the angle $BAC$ was proved equal to the angle $CDB$,
the side $AB$ was proved equal to the side $CD$, and the side $AC$ was proved equal to the side
$BD$; therefore in the parallelogrammic area $ACDB$ the opposite sides are equal to the opposite sides
and the opposite angles are equal to the opposite angles.

For the triangles $ABC,DCB$, the side $AB$ is equal to the side $CD$ and the side $BC$ is common,
the two sides $AB,BC$ are equal respectively to the two sides $DC,CB$;
and the angle $ABC$ is equal to the angle $BCD$;
therefore the base $AC$ of $ABC$ is equal to the base $DB$ of $DCB$ and the triangle
$ABC$ is equal to the triangle $DCB$ (I.4).
Since the triangle $ABC$ is equal to the triangle $DCB$, the diameter $BC$ bisects the parallelogram
$ACDB$.
\end{proof}

Proclus 393--394 \cite[p.~309]{proclus}:

\begin{quote}
It seems also that this very term ``parallelogram'' was coined by the author of the {\em Elements} and that
it was suggested by the preceding theorem. For when he had
shown that the straight lines connecting equal and parallel
lines in the same directions are themselves equal and parallel
he had clearly shown that both pairs of opposite sides, the
connecting and the connected lines, are parallel; and he
rightly called the figure enclosed by parallel lines a ``parallelogram,''
just as he had designated as ``rectilinear'' the
figure enclosed by straight lines.
\end{quote}

I.35: ``Parallelograms which are on the same base and in the 
same parallels are equal to one another.''

\begin{figure}
\begin{center}
\includegraphics{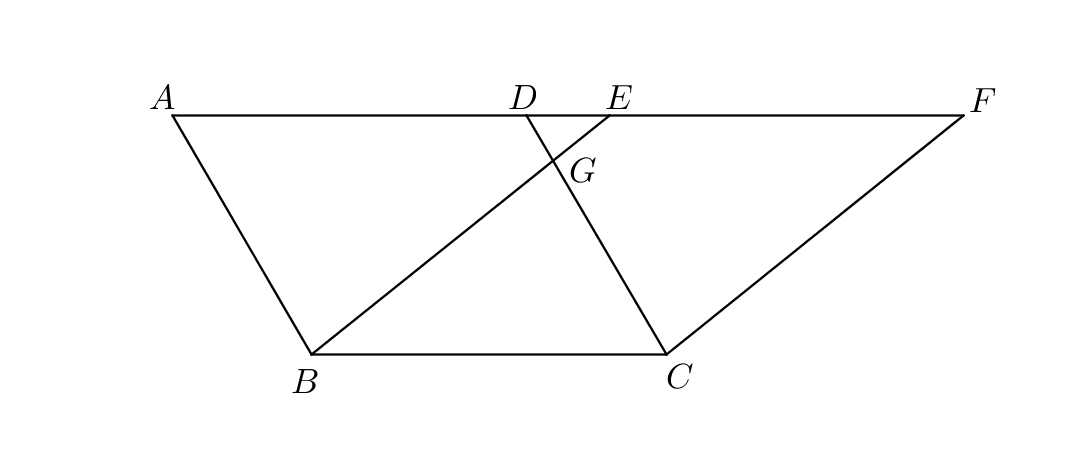}
\end{center}
\caption{I.35: P 32v, F 10v, B 25r, V 37}
\label{I35}
\end{figure}

\begin{proof}
Let $ABCD, EBCF$ be parallelograms on the same base $BC$ and in the same parallels
$AF,BC$. Because $ABCD$ is a parallelogram, the opposite sides $AD,BC$ are equal,
and because
$EBCF$ is a parallelogram, the opposite sides $EF,BC$ are equal (I.34); therefore $AD,EF$ are equal (Common Notion 1).
Let $DE$ be added to each of $AD,EF$; thus $AE$ is equal to $DF$ (Common Notion 2).
And because $ABCD$ is a parallelogram, the opposite sides $AB,DC$ are equal (I.34). 
Therefore the sides $EA,AB$ of the triangle $EAB$ are equal respectively 
to the sides $FD,DC$ of the triangle $FDC$, and the angle; 
and the straight line $FA$ falling on the parallel straight lines 
$DC,AB$ makes the exterior angle $FDC$ equal to the interior and opposite
angle $EAB$ (I.29); therefore the base $EB$ of the triangle $EAB$ is equal to the base
$FC$ of the triangle $FDC$ and the triangle $EAB$ is equal to the triangle $FDC$ (I.4).

Let the triangle $DGE$ be subtracted from each of the triangles $EAB,FDC$; then
the remaining trapezium $ABGD$ is equal to the remaining trapezium $EGCF$ (Common Notion 3).
Let the triangle $GBC$ be added to each of the trapezia $ABGD,EGCF$; then
the  parallelogram $ABCD$ is equal to the parallelogram $EBCF$ (Common Notion 2).
\end{proof}

Heath \cite[p.~327]{euclidI}:

\begin{quote}
It is important to observe that we are in this proposition introduced for
the first time to a new conception of equality between figures. Hitherto we
have had equality in the sense of {\em congruence} only, as applied to straight lines,
angles, and even triangles (cf. I.4). Now, without any explicit reference to
any change in the meaning of the term, figures are inferred to be {\em equal} which
are equal in {\em area} or in {\em content} but need not be of the same {\em form}. No
{\em definition} of equality is anywhere given by Euclid; we are left to infer its
meaning from the few {\em axioms} about ``equal things.''
\end{quote}

Proclus 396--397 \cite[pp.~312--313]{proclus}:

\begin{quote}
It may seem a great puzzle to those inexperienced in this
science that the parallelograms constructed on the same base 
[and between the same parallels] should be equal to one
another. For when the sides of the areas constructed on the
same base can be extended indefinitely -- and we can increase
the length of these sides of the parallelograms as far as we
can extend the parallel lines -- we may well ask how the areas
can remain equal when this happens. For if the breadth is
the same (since the base is identical) while the side becomes
greater, how could the area fail to become greater? This
theorem, then, and the following one about triangles belong
among what are called the ``paradoxical'' theorems in mathematics.
The mathematicians have worked out what they call
the ``locus of paradoxes,'' as the Stoics have done in their
dogmas, and this theorem is included among them. Most
people at least are immediately startled to learn that multiplying
the length of the side does not destroy the equality of
the areas when the base remains the same. The truth is,
nevertheless, that the equality or inequality of the angles
is the factor of greatest weight in determining the increase or
decrease of the area. For the more unequal we make the
angles, the more we decrease the area, if the side and base
remain the same; hence if we are to preserve equality, we
must increase the side.
\end{quote}

Proclus 398 \cite[p.~314]{proclus}:

\begin{quote}
With regard to the theorem before us we must realize that, when it
says the parallelograms are equal, it means the areas, not
the sides, are equal, for the statement is about the included
spaces, the areas; and also that in the demonstration of this
theorem our author for the first time mentions trapezia.
\end{quote}

I.36: ``Parallelograms which are on equal bases and in the same
parallels are equal to one another.''

\begin{figure}
\begin{center}
\includegraphics{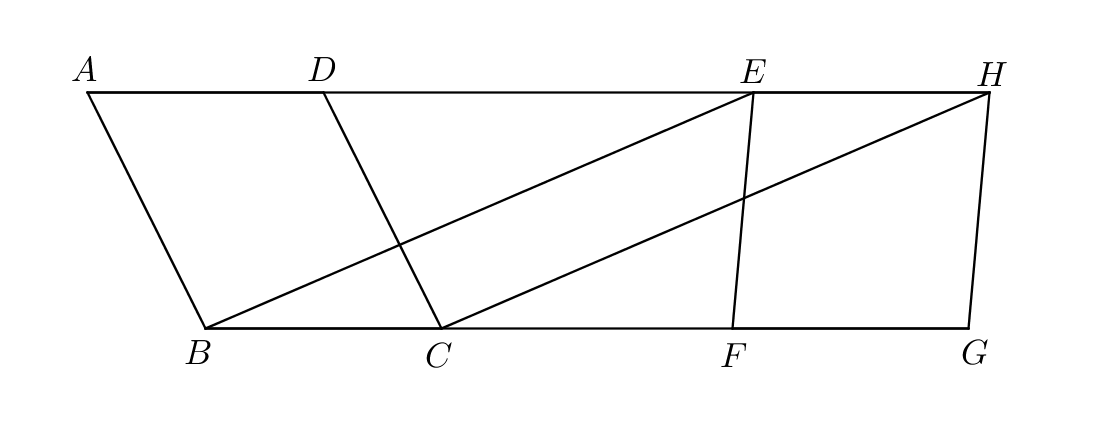}
\end{center}
\caption{I.36: P 32v, F 10v, B 25v, V 38}
\label{I36}
\end{figure}

\begin{proof}
Let $ABCD,EFGH$ be parallelograms which are on equal bases $BC,FG$ and in the same parallels $AH,BG$.
Let $BE,CH$ be joined (Postulate 1). Since by hypothesis $BC$ is equal to $FG$,
and $FG$ is equal to $EH$ (I.34), thus $BC$ is equal to $EH$ (Common Notion 1).
But the straight lines $BC,EH$ are equal and parallel, so the straight lines $EB,HC$ which join the extremities on the same
sides are themselves also equal and parallel (I.33); therefore $EBCH$ is a parallelogram.
Because the parallelograms $EBCH,ABCD$ are on the same base $BC$ and in the same parallels $BC,AH$,
they are equal (I.35).

It can be proved likewise that $EFGH$ is a parallelogram; and because the parallelograms $EFGH,EBCH$ are on the same
base $EH$ and in the same parallels $EH,BG$, they are equal (I.35). 

Therefore the parallelograms $ABCD,EFGH$ are equal (Common Notion 1).
\end{proof}

I.37: ``Triangles which are on the same base and in the same
parallels are equal to one another.''

\begin{figure}
\begin{center}
\includegraphics{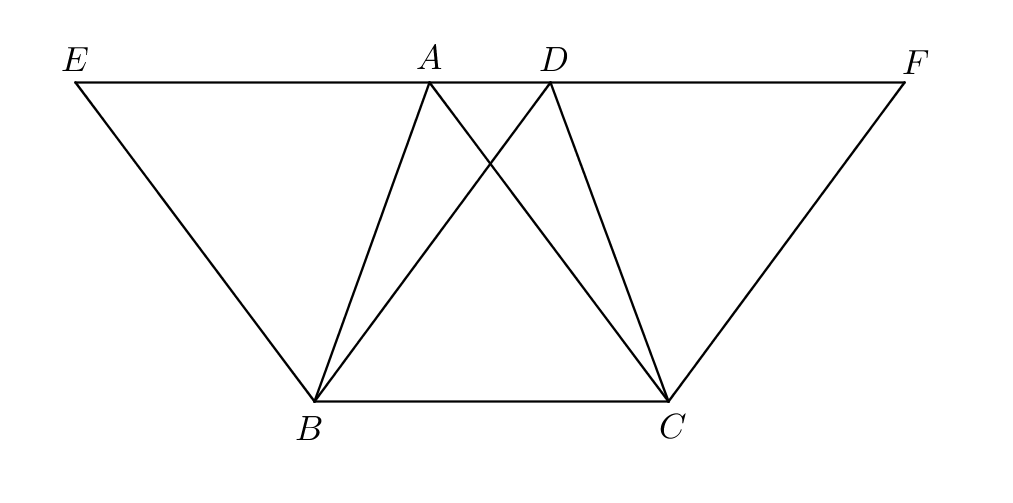}
\end{center}
\caption{I.37: P 33r, F 11r, B 26r, V 39}
\label{I37}
\end{figure}

\begin{proof}
Let $ABC,DBC$ be triangles on the same base $BC$ and in the same parallels $AD,BC$.
Let $AD$ be produced in both directions to $E,F$ (Postulate 2). Then through $B$ let $BE$ be drawn parallel
to $CA$, and through $C$ let $CF$ be drawn parallel to $BD$ (I.31). 
Then each of the figures $EBCA,DBCF$ is a parallelogram; and because they are on the same
base $BC$ and in the same parallels $BC,EF$, the parallelograms $EBCA,DBCF$ are equal (I.35).

In the parallelogram $EBCA$, the diameter $AB$ bisects the area (I.34); thus the triangle
$ABC$ is half of the parallelogram $EBCA$. And in the parallelogram $DBCF$, the diameter
$DC$ bisects the area (I.34); thus the triangle
$DBC$ is half of the parallelogram $DBCF$.
But the parallelogram $EBCA$ is equal to the parallelogram $DBCF$; therefore the triangle
$ABC$ is equal to the triangle $DBC$.
\end{proof}

I.38: ``Triangles which are on equal bases and in the same
parallels are equal to one another.''

\begin{figure}
\begin{center}
\includegraphics{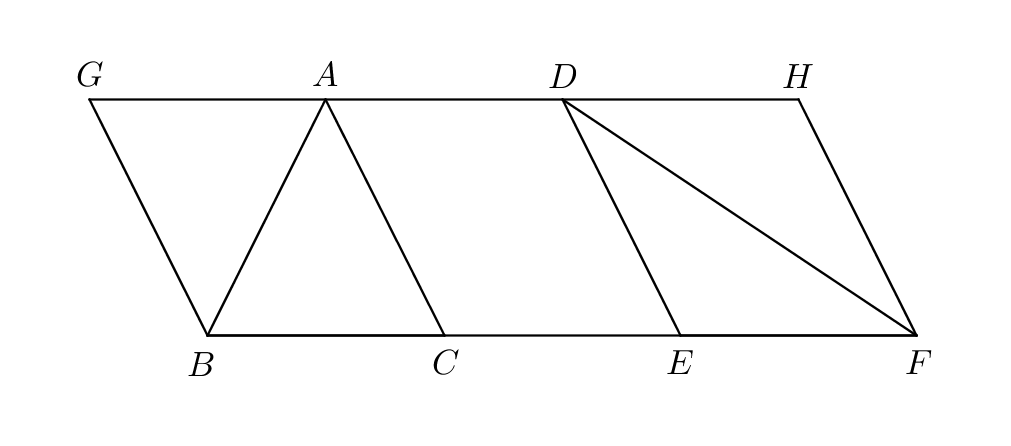}
\end{center}
\caption{I.38:  P 33r, F 11r, B 26v, V 40}
\label{I38}
\end{figure}

\begin{proof}
Let $ABC,DEF$ be triangles on equal bases $BC,EF$ and in the same parallels 
$BF,AD$. 
Let $AD$ be produced in both directions to $G,H$ (Postulate 2). 
Through $B$ let $BG$ be drawn parallel to $CA$ and through $F$ let $FH$ be drawn parallel
to $DE$ (I.31).
Then each of the figures $GBCA,DEFH$ is a parallelogram; 
and because they are on equal bases $BC,EF$ and in the same
parallels $BF,GH$, the parallelograms $GBCA,DEFH$ are equal (I.36).

In the parallelogram $GBCA$, the diameter $AB$ bisects the area (I.34); thus the triangle
$ABC$ is half of the parallelogram $GBCA$.
And in the parallelogram $DEFH$, the diameter $DF$ bisects the area (I.34); thus the triangle
$FED$ is half of the parallelogram $DEFH$. 
But the parallelogram $GBCA$ is equal to the parallelogram $DEFH$; therefore
the triangle $ABC$ is equal to the triangle $FED$.
\end{proof}

I.41: ``If a parallelogram have the same base with a triangle and
be in the same parallels, the parallelogram is double of the
triangle.''

\begin{figure}
\begin{center}
\includegraphics{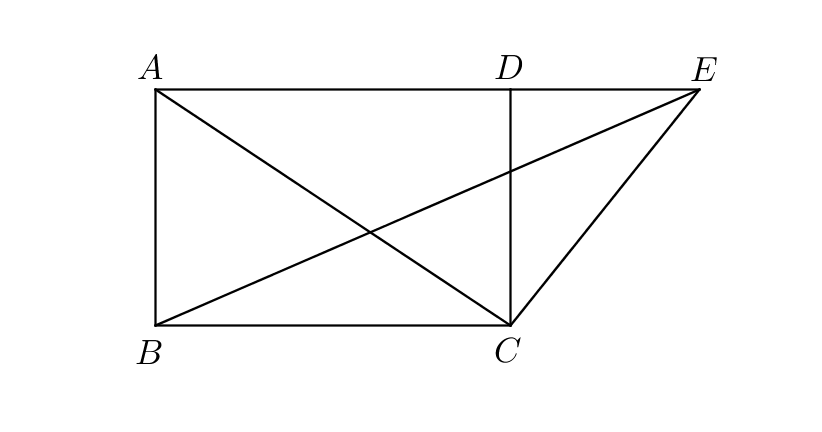}
\end{center}
\caption{I.41:  P 34r, F 11v, B 28r, V 42}
\label{I41}
\end{figure}

\begin{proof}
Let the parallelogram $ABCD$ have the same base $BC$ with the triangle $EBC$ and let it be in
the same parallels $BC,AE$. Let $AC$ be joined (Postulate 1). Then the triangles
$ABC,EBC$ are on the same base $BC$ and in the same parallels $BC,AE$; thus 
the triangle $ABC$ is equal to the triangle $EBC$ (I.37).
But in the parallelogram $ABCD$ the diameter $AC$ bisects the area (I.34); so
the parallelogram $ABCD$ is double of the triangle $ABC$. 
And the triangle $ABC$ has been proved  equal to the triangle $EBC$;
therefore the parallelogram $ABCD$ is double of the triangle $EBC$.
\end{proof}

I.42: ``To construct, in a given rectilineal angle, a parallelogram
equal to a given triangle.''

\begin{figure}
\begin{center}
\includegraphics{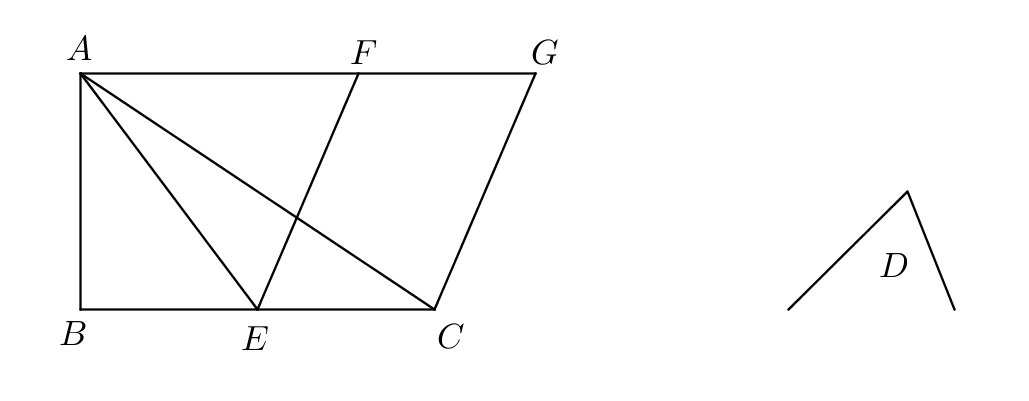}
\end{center}
\caption{I.42:  P 34v, F 12r, B 28r, V 42}
\label{I42}
\end{figure}

\begin{proof}
Let $ABC$ be the given triangle and $D$ the given rectilineal angle.
Let $BC$ be bisected at $E$ (I.10) and let $AE$ be joined (Postulate 1).
On the straight line $EC$ at the point $E$ on it let angle $CEF$ be constructed equal to the angle $D$ (I.23).
Through the point $A$ let the line $AG$ be drawn parallel to $EC$
and through the point $C$ let the line $CG$ be drawn parallel to $EF$ (I.31).
Then the figure $FECG$ is a parallelogram. 

The triangles $ABE,AEC$ are on equal bases $BE,EC$ and are in the same parallels $BC,AG$; therefore 
they are equal (I.38). Thus the triangle $ABC$ is double of the triangle $AEC$.

The parallelogram $FECG$ is on the same base $EC$ with the triangle $AEC$ and is in the same parallels $BC,AG$; therefore
the parallelogram $FECG$ is double of the triangle $AEC$ (I.42). But
the triangle $ABC$ was proved double of the triangle $AEC$; therefore the parallelogram $FECG$ is equal to the 
triangle $ABC$. 
And $FECG$ has the angle $CEF$ equal to the given angle $D$.
\end{proof}

Al-Nayrizi \cite[p.~185]{alnayriziI}, I.42, {\em protasis}:

\begin{quote}
We want to demonstrate how to construct a surface that is a 
parallelogram whose angle is equal to a known angle and [which is]
equal to a known triangle.
\end{quote}

\begin{figure}
\begin{center}
\includegraphics{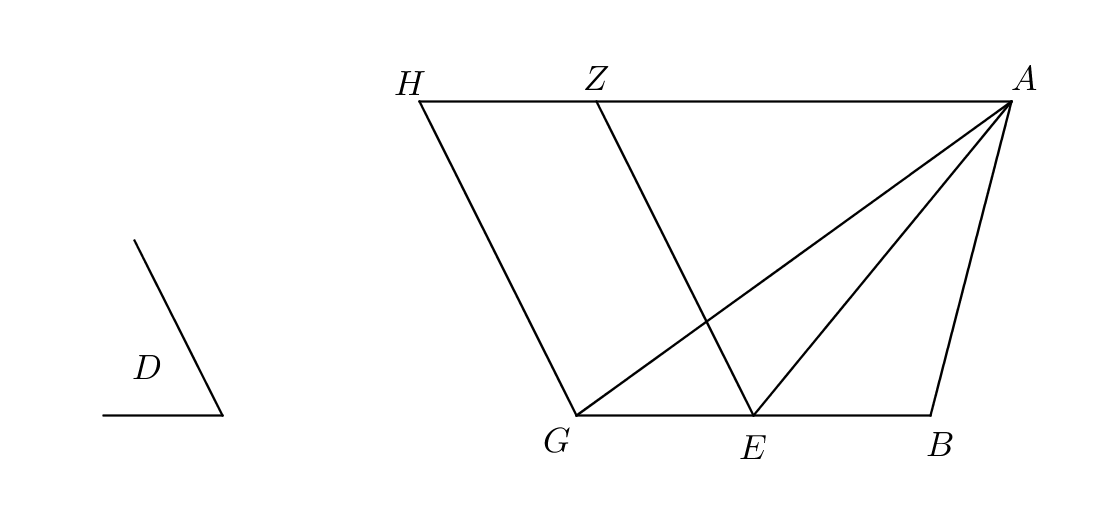}
\end{center}
\caption{Al-Nayrizi, I.42}
\label{nayriziI42}
\end{figure}

\begin{proof}
{\em ekthesis}: So let the known angle be the triangle $D$ and the known triangle  the triangle $ABG$.

{\em diorismos}: We want to construct a parallelogram whose angle is equal to the angle $D$ and which is equal to the triangle $ABG$.

{\em kataskeu\={e}}: So let us choose one of the sides of the triangle and let us cut it into two halves according to I.10. So let us suppose that we cut the side
$BG$ into two halves at the point $E$, and let us draw the line $AE$. Then let us construct at the point $E$ of the line $GE$ an angle
equal to $D$, according to I.23, and let this angle be $GEZ$. 
And let us draw from the point $G$ a line parallel to $EZ$ and from the point $A$ a line parallel to $BG$, according
to I.23; let the line from $G$ parallel to $EZ$ be $GH$ and let the line from $A$ parallel to $BG$ be $AZH$. 

{\em apodeixis}: Then, since the two triangles $ABE,AEG$ are upon equal bases $BE,EG$ and between
the same parallel lines $BG,AH$, the triangle $AEG$ is equal to the triangle $ABE$, according to
I.38. So the triangle $ABG$ is double of the triangle $AEG$. But the surface $GEZH$ is a parallelogram,
and its base $EG$ is the  base of the triangle $AEG$, and the two of them are between the 
two parallel lines $BG,AH$; so
the surface $GEZH$ is double of the triangle $AGE$, according to I.41. 
And we have proved that the triangle $AGB$ is double the triangle $AGE$, and doubles of the same thing are equal,
so the parallelogram $GEZH$ is equal to the triangle $AGB$.

{\em sumperasma}: So we have constructed a surface $GEZH$ that is a parallelogram equal to the known
triangle $ABG$ and whose angle
$GEZ$ is equal to the known angle $D$, which is what we wanted to demonstrate.
\end{proof}

I.43: ``In any parallelogram the complements of the parallelograms
about the diameter are equal to one another.''

\begin{figure}
\begin{center}
\includegraphics{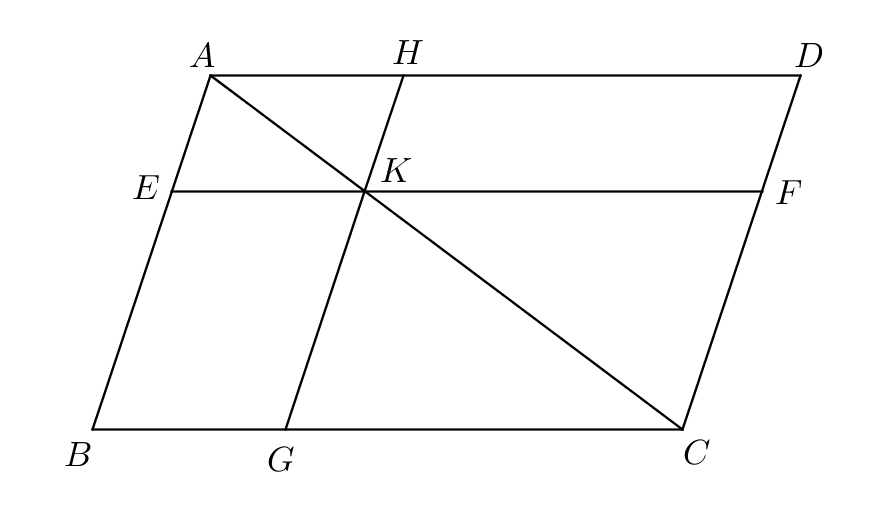}
\end{center}
\caption{I.43: P 34v, F 12r, B 28v, V 43}
\label{I43}
\end{figure}

\begin{proof}
Let $ABCD$ be a parallelogram and $AC$ its diameter; and about
$AC$ let $EH,FG$ be parallelograms and $BK,KD$ the so-called complements.
Since $ABCD$ is a parallelogram, the diameter $AC$ bisects its area (I.34); thus the triangle $ABC$ is equal to the triangle
$ACD$. And since $EH$ is a parallelogram, the diameter $AK$ bisects its area (I.34); thus the triangle $AEK$ is equal to the
triangle $AHK$. And since $FG$ is a parallelogram, the diameter $KC$ bisects its area (I.34); thus the triangle
$KFC$ is equal to the triangle $KGC$
Then, since the triangle $AEK$ is equal to the triangle $AHK$ and the triangle
$KFC$ is equal to the triangle $KGC$, the triangle $AEK$ together with the triangle $KGC$ is equal to the triangle
$AHK$ together with the triangle $KFC$ (Common Notion 2).
The triangle $ABC$ is equal to the triangle $ADC$. Let the triangles $AEK,KGC$ be subtracted from the triangle $ABC$ and let
the triangles $AHK,KFC$ be subtracted from the triangle $ADC$; then the complement $BK$ which remains is equal
to the complement $KD$ which remains (Common Notion 3). 
\end{proof}

Scholia for I.43 \cite[pp.~201--203]{euclidisV}.

Proclus 418--419 \cite[p.~331]{proclus}:

\begin{quote}
The term ``complements'' was derived by the author of the 
{\em Elements} from the thing itself, since complements fill the
whole of the area outside the two parallelograms. This is why
he does not regard it as deserving of special mention in the
Definitions. It would have required a complicated explanation
to make us understand what a parallelogram is and what are
the parallelograms that are constructed about the same diameter
as the whole; for only after these had been explained
would the meaning of ``complement'' have become clear.
Those parallelograms are about the same diameter which have
a segment of the entire diameter as their diameter; otherwise
they are not about the same diameter. 
\end{quote}

Szabo \cite[pp.~344--345]{szabo}:

\begin{quote}
If, for example, the Pythagoreans had
not known that the {\em parapleromata} were  equal, they would not have
been able to develop their method of {\em application of areas}. This method
was mentioned in {\em Example 1} above, where it was pointed out that
finding a fourth proportional to three given {\em numbers} or {\em magnitudes} (i.e.
an $x$ such that $a:b=c:x$) could be construed as the problem of
finding a rectangle with a given side ($a$) which has the same area as a
given rectangle ($bc$). 

If we think of the given rectangle
($bc$) as a {\em parapleroma} and the rectangle which the line $a$ makes with
one of its sides ($c$ or $b$) as a `parallelogram about the diagonal', then the
second `parallelogram about the diagonal' ($bx$ or $cx$) can be constructed 
by extending the side ($c$ or $b$) of the original rectangle until it meets the 
continuation of the diagonal of $ac$ or of $ab$. Now these 
three rectangles together determine uniquely the second {\em parapleroma}
and, in particular, its side $x$ which we set out to find. (This is sometimes
called {\em parabolic} application of areas, because the area $bc$ is {\em applied}
to the line $a$; cf. the word παραβάλλειν.)
\end{quote}

Dalimier \cite{dalimier} on {\em parapleromata}.

I.44: ``To a given straight line to apply, in a given rectilineal
angle, a parallelogram equal to a given triangle.''

\begin{figure}
\begin{center}
\includegraphics{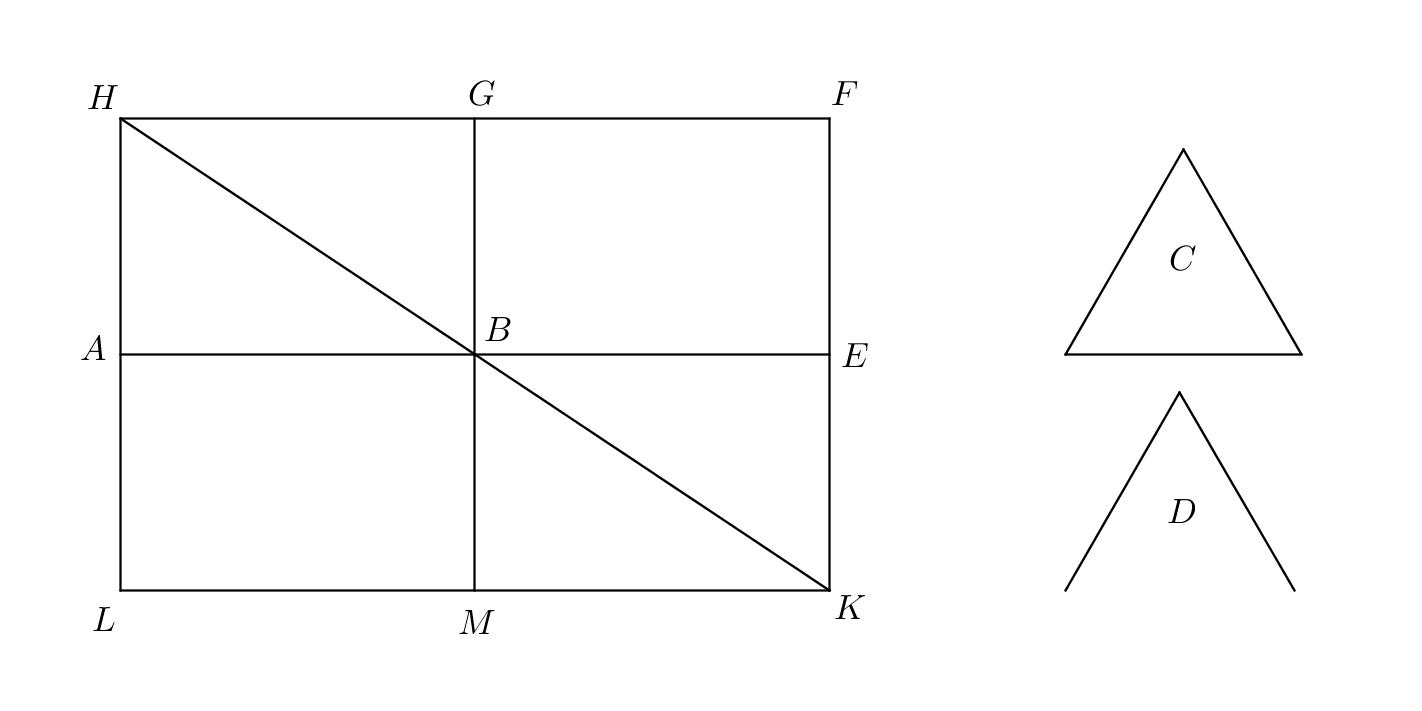}
\end{center}
\caption{I.44:  P 36r, F 12v, B 29r, V 45}
\label{I44}
\end{figure}

\begin{proof}
{\em ekthesis}: Let $AB$ be the given straight line, $C$ the given rectangle, and
$D$ the given rectilineal angle.

{\em diorismos}: Thus, it is required to apply
a parallelogram equal to the given triangle $C$, to the given straight line $AB$, in the
given angle $D$.

{\em kataskeu\={e}}: Let the parallelogram $BEFG$ equal to the triangle $C$
be constructed in the angle $EBG$ equal to $D$ (I.42); let it be placed so that
$BE$ is in a straight line with $AB$.\footnote{What is used here is more than what I.42 provides.}
Through the point $A$ let the straight line $AH$ be drawn parallel either to $BG$ or $EF$ (I.31).

{\em apodeixis}: Then, since the straight line $HF$ falls upon the parallel straight lines,
the interior angles on the same side $AHF,HFE$ are equal to two right angles (I.29).
Therefore the angles $BHG,GFE$ are less than two right angles; thus
the straight line $HF$ falling on the the two straight lines $HB,FE$ 
makes $BHG,GFE$,
the interior angles on the same side, less than two right angles,
so the straight lines $HB,FE$ if produced will meet on the same side of $HF$
as are the angles $BHG,GFE$ (Postulate 5).
Let the lines $HB,FE$ be produced at meet at $K$.
Through the point $K$ let the straight line $KL$ be drawn parallel to either $EA$ or $FH$ (I.31);
then let the straight lines $HA,GB$ be produced to the points $L,M$ (Postulate 2).
Then $HLKF$ be a parallelogram; $HK$ is its diameter, $AG,ME$ are parallelograms about its diameter,
and $LB,BF$ are the so-called complements;
therefore $LB$ is equal to $BF$ (I.43).
But $BF$ is equal to the triangle $C$; therefore $LB$ is also equal to $C$ (Common Notion 1).  
And the lines $AE,GM$ cutting each other make the vertical angles
$GBE,ABM$ equal (I.15); but the angle $GBE$ is equal to the angle $D$,
therefore the angle $ABM$ is also equal to the angle $D$ (Common Notion 1). 

{\em superasma}: Therefore the parallelogram $LB$ equal to the given triangle $C$
has been applied to the given line $AB$ in the angle $ABM$ equal to $D$.
\end{proof}

Scholia for I.44 \cite[pp.~203--207]{euclidisV}.

Vitrac \cite[p.~276]{vitracI} writes:

\begin{quote}
Compte-tenu de l'\'egalit\'e des compl\'ements \'etablie dans la Prop. 43, le
triangle $C$ et l'angle $D$ \'etant donn\'es, il suffit de construire un parall\'elogramme
\'equivalent \`a $C$ admettant un angle \'egal \`a $D$ gr{\^a}ce \`a la Prop. 42,
de la placer de telle fa{\c c}on que l'un de ses c{\^o}t\'es soit un alignement avec la
droite $AB$ donn\'ee -- ce qui suppose d\'ej\`a que l'on accorde une certaine
latitude \`a ce parall\'elogramme quant \`a sa position ou qu l'on autorise le
{$\ll$}d\'eplacement{$\gg$} des figures -- puis de compl\'eter la figure gr{\^a}ce \`a la
th\'eorie des parall\`eles pour obtenir une figure du m{\^e}me type que celle de
la Prop. 43. L'autre compl\'ement sera donc un parall\'elogramme \'equivalent,
construit sur la droite donn\'ee et \'equiangle au premier compl\'ement.
\end{quote}

Al-Nayrizi \cite[p.~188]{alnayriziI}, I.44, {\em protasis}:

\begin{quote}
We want to demonstrate how to construct, upon a known
straight line, a surface that is a parallelogram equal to a known triangle
and whose angle is equal to a known angle.
\end{quote}

\begin{figure}
\begin{center}
\includegraphics{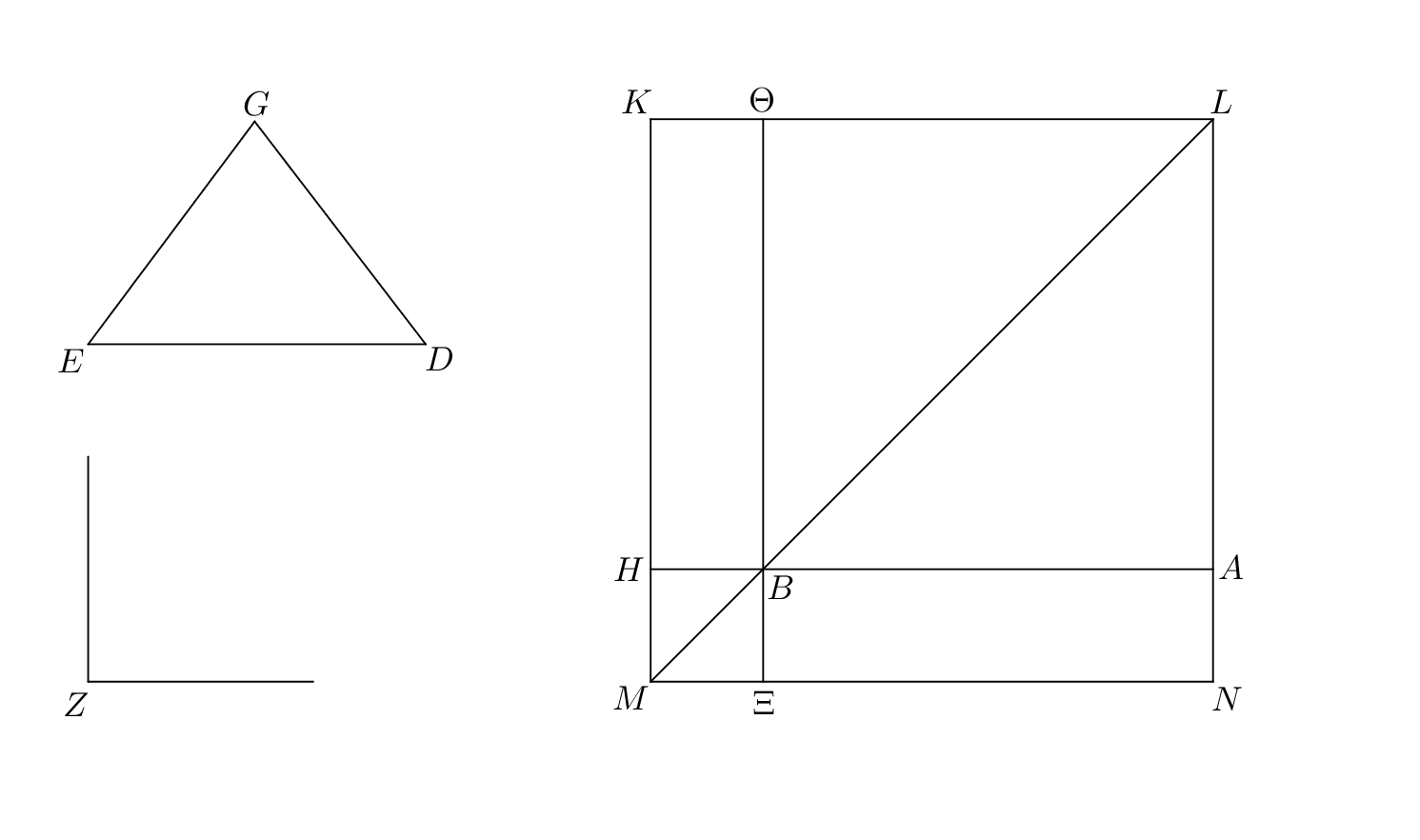}
\end{center}
\caption{Al-Nayrizi, I.44}
\label{nayriziI44}
\end{figure}

\begin{proof}
{\em ekthesis}: So let the known line be the line $AB$, the known triangle the triangle $GDE$, and the known 
angle the angle $Z$.

{\em diorismos}: We want to demonstrate how to construct on the line $AB$ a surface that is a parallelogram
equal to the triangle $GDE$ and whose angle is equal to the angle $Z$.

{\em kataskeu\={e}}: Extend $AB$ (Postulate 2) and 
cut off $BH$ equal to half $DE$ (I.3). 
Use I.42 to construct on $BH$ the parallelogram
$B\Theta KH$ equal to the triangle $GDE$ and the angle $HB\Theta$ equal to the angle
$Z$.\footnote{What is used here is more than what I.42 provides.}
Extend $\Theta K$ to $L$ (Postulate 2). 
Use I.31 to draw through $A$ a line parallel to $B\Theta$. Let $L$ be 
the
intersection of this line and 
$K\Theta L$. Then draw $LB$.

The lines $KH,AL$ are parallel and the line $LK$ falls on them, so according to I.29, 
the interior angles $LKM,KLN$ are equal to two right angles. 
The line $LK$ falling on the lines $LB,KH$ makes the interior angles
$LKM,KLM$, which are less than two right angles. Thus by I.29, the lines $KH,LB$ extended
meet at some point; let $M$ be this point. 
By I.31, draw $MN$ parallel to $KL$. Extend the line $LA$ and let $N$ be the point at which it meets the line $MN$.
Then extend $\Theta B$ and let $\Xi$ be the point at which it meets the line $MN$.

{\em apodeixis}: The surface $LM$ is a parallelogram with diameter $LM$ and the surfaces $A\Theta,\Xi H$ are parallelograms about the diameter. 
By I.43, the complements $NB,BK$ are equal. The parallelogram $B\Theta KH$ was constructed
equal to the triangle $GDE$, so 
the parallelogram  $AB\Xi N$ is equal to the triangle $GDE$.
Now, the lines $AH,\Theta\Xi$ cut each and make vertical angles 
$HB\Theta,AB\Xi$, and by I.15 these angles are equal.
The angle $HB\Theta$ was constructed equal to the angle $Z$, so
the angle $AB\Xi$ is equal to the angle $Z$.

{\em superasma}: We have constructed on the line $AB$ the parallelogram $A\Xi$ equal
to the given triangle $GDE$ and whose angle $AB\Xi$ is equal to the angle $Z$.
\end{proof}

\begin{figure}
\begin{center}
\includegraphics{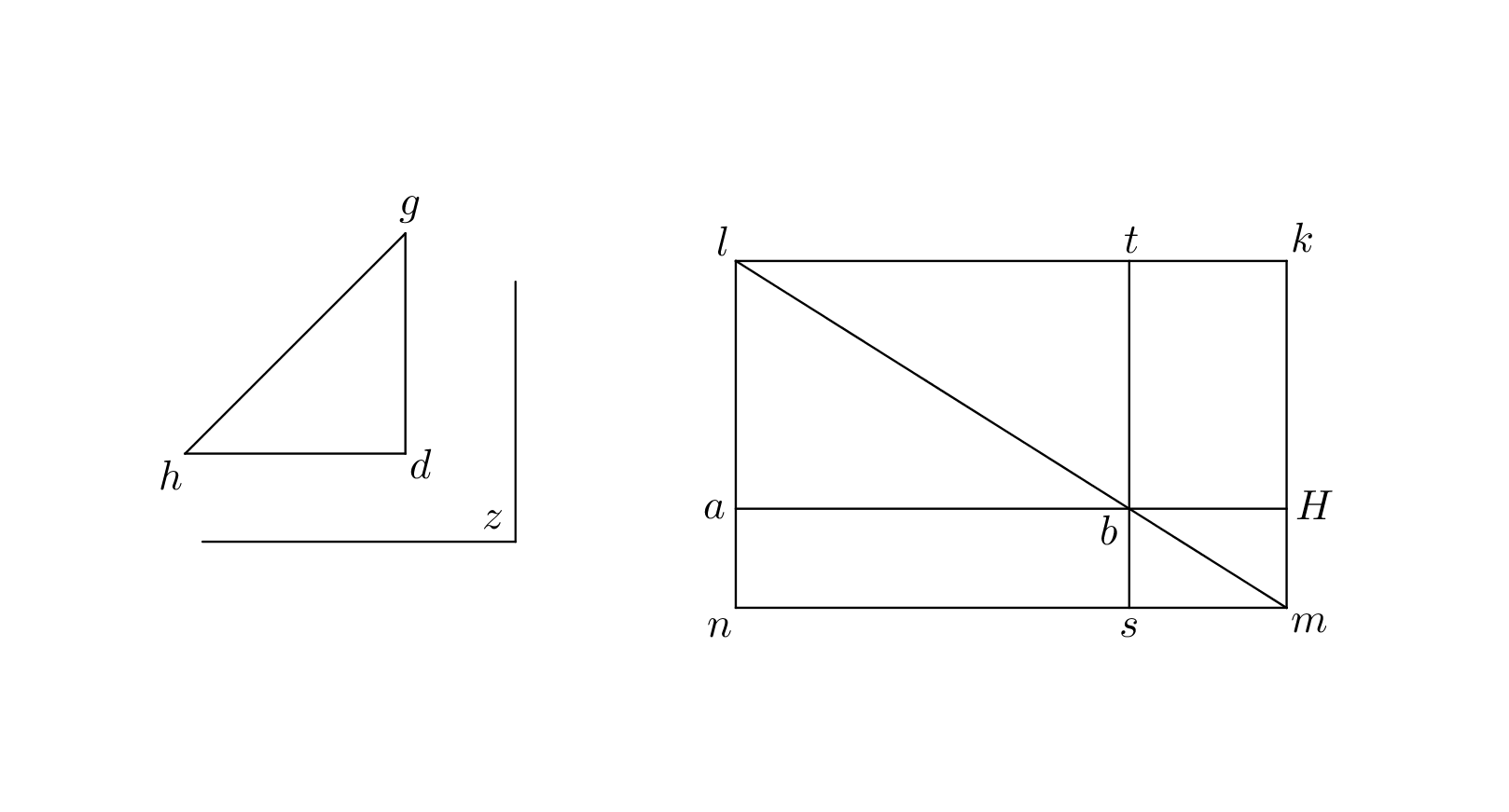}
\end{center}
\caption{Adelard of Bath, I.44}
\label{adelardI44}
\end{figure}

The proof of I.44 in Adelard of Bath \cite[pp.~66--67]{adelardI} is similar to the proof in Al-Nayrizi. The figure for I.44 in Adelard of Bath
is given in Figure \ref{adelardI44}.

Robert of Chester \cite[p.~129]{adelardII}, I.44: 

\begin{quote}
Proposita linea recta super eam superficiem equidistancium laterum, cuius angulus sit angulo assignato equalis, ipsa vero superficies triangulo assignato equalis, designare.
\end{quote}

\begin{figure}
\begin{center}
\includegraphics{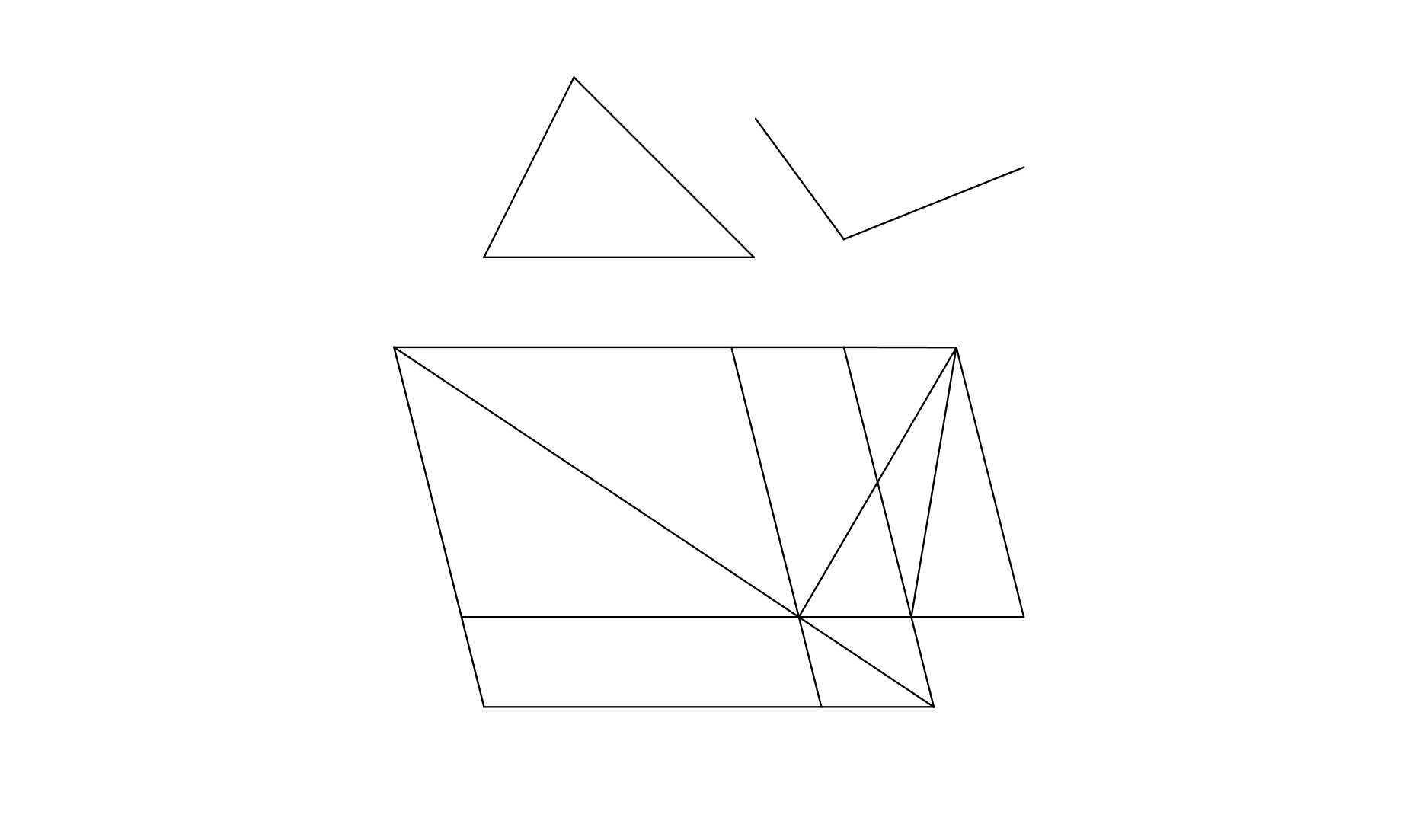}
\end{center}
\caption{Robert of Chester (?), I.44}
\label{chesterI44}
\end{figure}

\begin{quote}
Date linee tamquam dimidium basis dati trianguli adicito directe at- que super adiectam paralellogramum equale dato triangulo, cuius angulus super communem terminum date atque addite linee statutus dato angulo fiat equus, compleaturque figura, in cuius spacio inventum paralellogramum sit unum supplementum et alterum sit super datam lineam. Itaque per XLI\textsuperscript{am} atque per premissam promissum provenire necesse est.

Sic autem super adiectam statues paralellogramum equum dato triangulo: Protrahe eam, donec equetur toti basi dati trigoni. Deinde super duos terminos tocius adiecte linee fac duos angulos equales illis duobus qui sunt super basim dati trianguli, et conclude triangulum, quem ex XXVl\textsuperscript{a} convinces esse equalem dato triangulo. Itaque per XLII\textsuperscript{am} perfice.
\end{quote}

The figure for I.44 in Robert of Chester is given in Figure \ref{chesterI44}.

Johannes de Tinemue \cite[pp.~67--68]{adelardIII}, I.44:

\begin{quote}
Proposita linea recta super eam superficiem equidistantium laterum, 
cuius angulus sit angulo assignato equalis, ipsa vero superficies triangulo
assignato equalis, designare.
\end{quote}

\begin{quote}
Esto exemplum $ab$ linea et $bcd$ triangulus et $e$ angulus.

Dispositio. Protrahatur itaque $ab$ in continuum et directum ad equilitatem
$bc$ basis trianguli et $bc$ quoque dividatur equaliter in $o$ linea dividente
protracta ad $d$. Erit itaque $o$ angulus dexter superior vel maior vel minor
vel equalis $e$ secundum quod tripliciter variatur.
\end{quote}

\begin{figure}
\begin{center}
\includegraphics{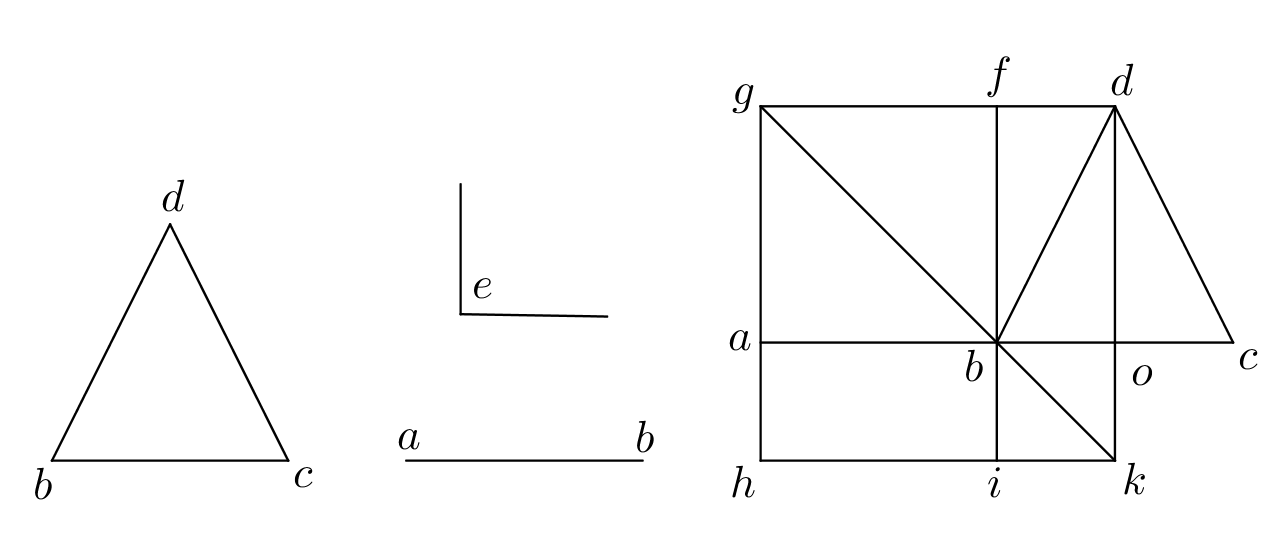}
\end{center}
\caption{Johannes de Tinemue, I.44}
\label{tinemueI44a}
\end{figure}

Figure \ref{tinemueI44a}.

\begin{quote}
Dispositio. Sit itaque primo equalis. Deinde describatur $ad$ parallelogramum
secundum exigentiam $ao, od$. A $b$ vero ducatur $bf$ equidistanter
ad $od$, $gb$ diametro interiecta et protracta in occursum $ga,fb$ vero
protracta in occursum $hk$ et obveniet $hd$ parallelogramum quasi distinctum
per $af$ et $io$ parallelograma circa diametrum et $hb,bd$ supplementa.
Ergo secundum premissam $hb,bd$ supplementa sunt equalia. Sed $bd$ parallelogramum
et $bcd$ triangulus sunt equalia secundum tenorem antepremisse
ratiocinationis, ergo $hb$ et $bcd$ sunt equalia. Preterea $o$ angulus dexter
superior est equalis $e$ angulo dato. Sed $o$ sinister inferior est equalis $o$ dextra
superiori secundum \textsc{xv}\textsuperscript{am} et $b$ extrinsecus reest equalis $o$
inferiori secundum
29\textsuperscript{am}. Ergo $b$ angulus est equalis $e$ sicque super $ab$ lineam describitur $ai$
parallelogramum equale $bcd$ cuius $b$ angulus est equalis $e$ angulo proposito.
Quod proposuiumus.
\end{quote}

\begin{figure}
\begin{center}
\includegraphics{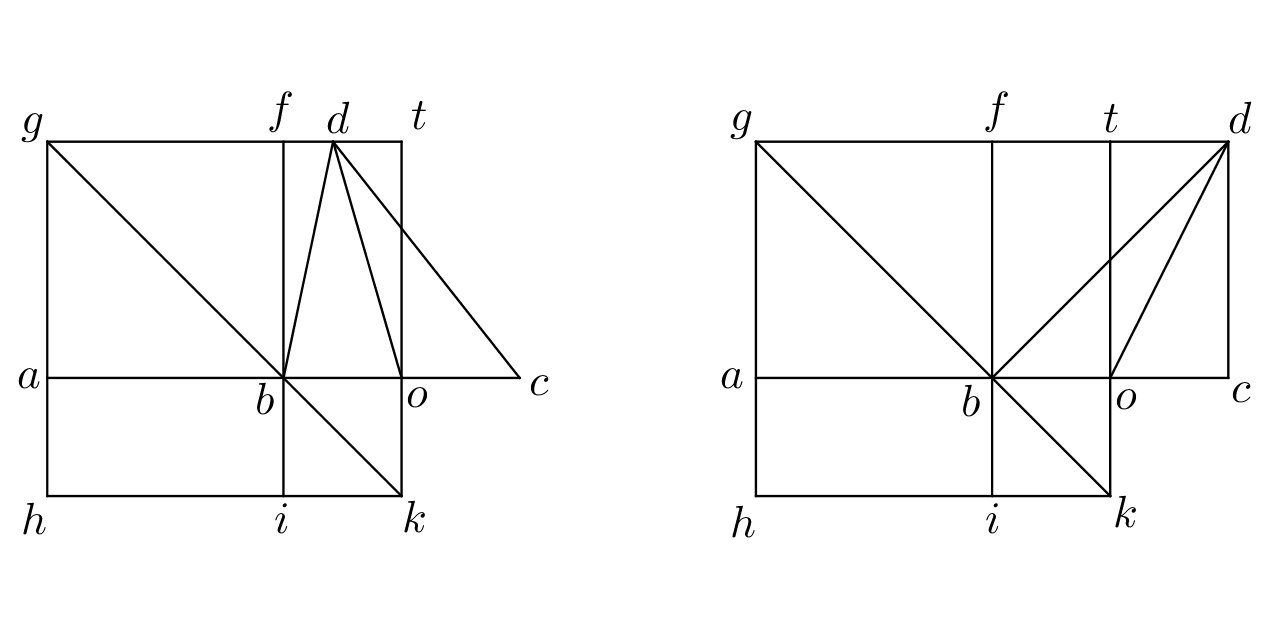}
\end{center}
\caption{Johannes de Tinemue, I.44}
\label{tinemueI44b}
\end{figure}

Figure \ref{tinemueI44b}.

\begin{quote}
Sit deinde $o$ dexter superior maior quam $e$ et ad eius equalitatem 
resecetur $ot$ linea protracta in occursum $dt$ linee equidistantis ad $bc$. A 
$b$ angulo vero educatur $bf$ equidistans ad $ot$. Deinde tota dispositio inclinetur
secundum $ot$ dextrorsum et obveniet propositum secundum superiorem
ratiocinationem ad consequentiam protractam.

Sit denique $o$ angulus superior dexter minor quam $e$ et augeatur ad equalitatem
eiusdem $e$. Deinde secundum priorem ratiocinationem disponendi
tota machina inclinetur sinistrorsum secundum $ot$. Et exibit propositum.
\end{quote}

\begin{proof}
Let $ab$ be the line, let $bcd$ be the triangle, and let $e$ be the angle.

Extend $ab$ in a straight line so that $bc$ is equal to the base of the triangle, let
$bc$ be divided equally at $o$, and let a line be drawn from $o$ to $d$. Then 
the top right angle at  $o$ is either greater than or less than or equal to $e$.

First take the case where the angle is equal to $e$.  
Let the parallelogram $ad$ be described with the sides $ao,od$. 
Let $bf$ be drawn from $b$ parallel to $od$. 
Extend $gb$ to intersect $od$ at some point $k$. 
Draw $kh$ parallel to $gd$ and suppose it meets the line $ga$ at some point $h$.
Extend $fb$ to meet $hk$ at some point $i$. Then 
the parallelogram $gk$ is divided into the parallelograms $af$ and $io$ about the diameter
and the complements $hb,bd$. So according
to I.43, the complements $hb,bd$ are equal. But
the parallelogram $bd$ and the triangle $bcd$ are equal. Hence
the parallelogram $hb$ and the triangle $bcd$ are equal.
Moreover the top right angle at $o$  is equal to the given angle $e$. But the bottom left angle at $o$ is
equal to the top right angle at $o$ by I.15; so $bok$ is equal to the given angle $e$.
And the angle $abi$ is equal to the angle $bok$ by I.29; so $abi$ is equal to the given angle $e$, i.e.,
the angle $b$ is equal to the given angle $e$.
Therefore, on the line $ab$ the parallelogram $ai$ has been constructed equal to the given triangle
$bcd$ and with the angle at $b$ equal to the given angle $e$.
\end{proof}

Hermann of Carinthia \cite{hermann}

Gerard of Cremona \cite{gerard}

Campanus \cite[pp.~91--92]{campanusI}:

\begin{quote}
Proposita recta linea super eam superficiem equidistantium laterum,
cuius angulus sit angulo assignato equalis, ipsa vero superficies triangulo
assignato equalis, designare.
\end{quote}

\begin{figure}
\begin{center}
\includegraphics{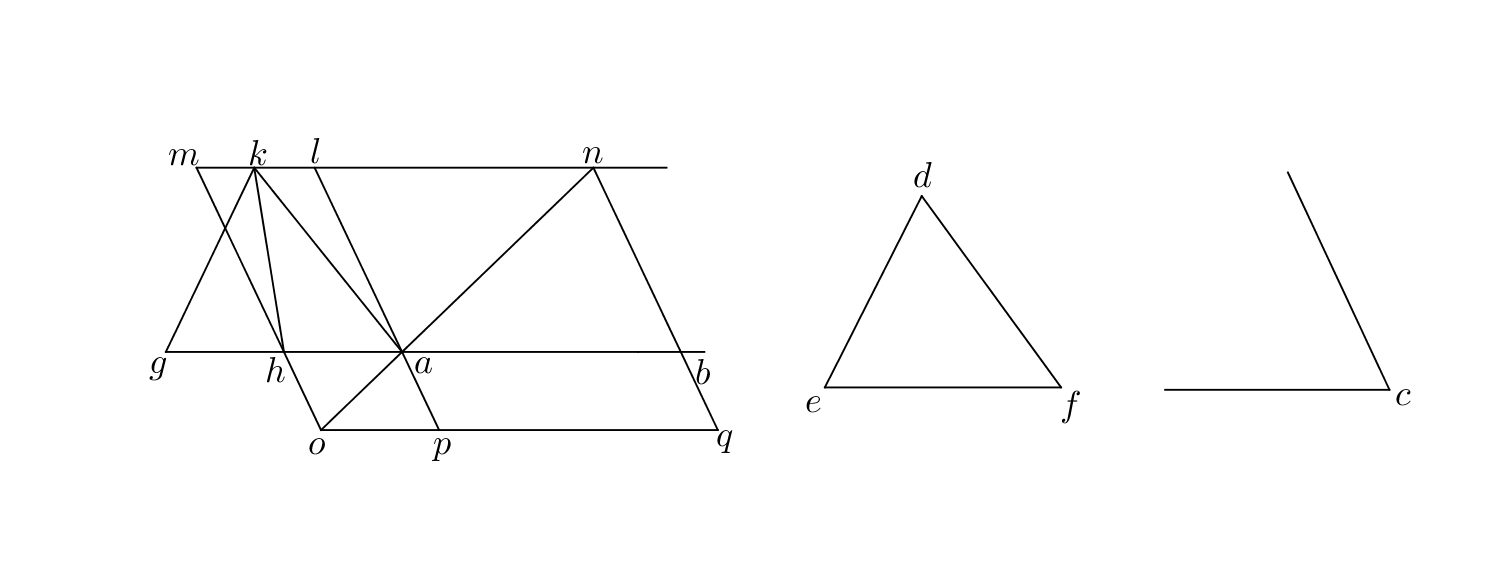}
\end{center}
\caption{Campanus, I.44}
\label{campanusI44}
\end{figure}

\begin{quote}
Designare superficiem equidistantium laterum super lineam aliquam est
lineam ipsam facere latus unum ipsius superficiei. Sit ergo data linea $ab$ et
datus angulus $c$ et datus triangulus $def$. Super lineam $ab$ volo designare
superficiem unam equidistantium laterum ita quod linea $ab$ sit unum ex lateribus
eius cuius uterque duorum angulorum contra se positorum sit
equalis angulo $c$ et ipsa totalis superficies sit equalis triangulo $def$. Differt
autem hec a 42 quia hic datum latus unius superficiei describende scilicet
linea $ab$, ibi autem nullum. Cum ergo hoc volo facere ad lineam $ab$,
adiungo secundum rectitudinem lineam $ag$ quam pono equalem linee $ef$
basi trianguli dati super quam constituo triangulum unum ei equalem et
equilaterum. Quod hoc modo facio.
Constituo angulum $agk$ equalem angulo $e$
et angulum $gak$ equalem angulo $f$ per 23.
Et quia $ga$ posita fuerat equalis $ef$, erit per 26 triangulus $gak$ equalis et equilaterus triangulo $efd$.
Dividam ergo $ga$ per equalia in puncto $h$ et protraham $kh$ et producam
a puncto $k$ lineam $mkn$ equidistantem linee $gh$ eritque per 38 triangulus
$ahk$ equalis triangulo $ghk$.
Tunc super punctum $a$ linee $ga$ faciam
angulum $gal$ equalem angulo $c$ dato et complebo super basim $ah$ et inter
lineas $gh$ et $mn$ equidistantes superficiem equidistantium laterum $mlha$
que per 41 dupla erit ad triangulum $kha$, quare equalis totali triangulo 
$kga$ quare et triangulo $def$ proposito. Protraham ergo $bn$ equidistantem
$al$ et producam diametrum $na$ quam protraham quousque concurrat cum
$mh$ in puncto $o$ et complebo superficiem equidistantium laterum $monq$ et 
protraham $la$ usque ad $p$. Eritque per precedentem supplementum $abpq$
equale supplemento $mlha$, quare triangulo $def$ et quia per 15 angulus
$lah$ est equalis angulo $bap$ et ideo angulus $bap$ equalis angulo $c$, patet
super datam lineam $ab$ descriptam superficiem esse equidistantium laterum
$abpq$ equalem dato triangulo $def$ cuius uterque duorum angulorum
contra se positorum qui sunt $a$ et $q$ est equalis dato angulo $c$. Quod fuit
propositum.
\end{quote}

\begin{proof}
Let $ab$ be the given line, $c$ the given angle, and $def$ the given triangle.
Extend the line $ab$ by the line $ag$ equal to $ef$. 
Apply I.23 to make angle $agk$ equal to angle $e$;
apply I.23 to make angle $gak$ equal to angle $f$; 
let $k$ be the point of intersection of the lines $gk$ and $ak$. 
Because $ga$ was put equal to $ef$, by I.27, the triangle $gak$ is equal to the triangle $efd$.
Divide $ga$ into equal halves at the point $h$ and join $kh$; draw through $k$ a line $mkn$ parallel to the line $gh$;
then by I.38, the triangle $ahk$ is equal to the triangle $ghk$. Therefore the triangle $kga$ is double the triangle $kha$.
At the point $a$ on the line $ga$ make an angle $gal$ equal to the given angle $c$. Then on the base
$ah$ and between the parallel lines $gh$ and $mn$ complete the parallelogram $mlha$; by I.41 this is double the triangle $kha$ and therefore
the parallelogram $mlha$ is equal to the whole triangle $kga$ which itself is equal to the given triangle $def$.
Draw the line  $bn$ parallel to the line $al$ and produce $na$ until it intersects $mh$, and let $o$ be the point of intersection.
Then complete the parallelogram $monq$ and produce $la$ to $p$. Then by I.43,
the complement $abpq$ is equal to the complement $mlha$, hence
the parallelogram $abpq$ is equal to the given triangle $def$. 
By I.15, the lines $bh,lp$ cutting each other make equal vertical angles 
$lah,bap$.  But $gal$ is equal to the given angle $c$, so $bap$ is equal to the given angle $c$. 
Therefore, on the given line $ab$ the parallelogram $abpq$ has been described that is equal to the given triangle $def$
and both of the opposite angles $a$ and $q$ are equal to the given angle $c$. 
\end{proof}

Clavius, {\em Euclidis Elementorum libri XV}, {\em Opera mathematica}, pp.~72--73, scholium to I.43, gives two different
constructions for I.44.

\begin{proof}
{\em Ad datam rectam lineam, dato triangulo aequale
parallelogrammum applicare in dato angulo rectilineo.}

Let the given line be $A$, let the given triangle be $B$, and let
the given angle be $C$.
{\em Quod si quis optet, lineam ipsam $A$, datam, esse unum
latus parallelogrammi, non difficile erit
transferre parallelogrammum $FMLH$, ad rectam $A$, ex
iis, quae in scholio propos. 31. huius lib. docuimus.}
\end{proof}

Proclus 419--420 \cite[pp.~332--333]{proclus}:

\begin{quote}
Eudemus and his school tell us that these things -- the
application ({\em parabole}) of areas, their exceeding ({\em huperbole}),
and their falling short ({\em elleipsis}) -- are ancient discoveries of
the Pythagorean muse. It is from these procedures that later
geometers took these terms and applied them to the so-called
conic lines, calling one of them ``parabola,'' another ``hyperbola,''
and the third ``ellipse,'' although these godlike men
of old saw the significance of these terms in the describing of
plane areas along a finite straight line. For when, given a
straight line, you make the given area extend along the whole
of the line, they say you ``apply'' the area; when you make
the length of the area greater than the straight line itself,
then it ``exceeds''; and when less, so that there is a part of the
line extending beyond the area described, then it ``falls
short.'' Euclid too in his sixth book speaks in this sense of
``exceeding'' and ``falling short''; but here he needed ``application,''
since he wished to apply to a given straight line
an area equal to a given triangle, in order that we might be
able not only to construct a parallelogram equal to a given
triangle, but also to apply it to a given finite straight line. For
example, when a triangle is given having an area of twelve
feet and we posit a straight line whose length is four feet, we
apply to the straight line an area equal to the triangle when
we take its length as the whole four feet and find how many
feet in breadth it must be in order that the parallelogram may
be equal to the triangle. Then when we have found, let us say,
a breadth of three feet and multiplied the length by the
breadth, we shall have the area, that is, if the angle assumed 
is a right angle. Something like this is the method of  ``application''
which has come down to us from the Pythagoreans.
\end{quote}

Proclus 421 \cite[pp.~333--334]{proclus}:

\begin{quote}
Application and construction are not the same
thing, as we have said. Construction brings the whole figure
into being, both its area and all its sides, whereas application
starts with one side given and constructs the area along it,
neither falling short of the length of the line nor exceeding it,
but using it as one of the sides enclosing the area.
\end{quote}

I.45: ``To construct, in a given rectilineal angle, a parallelogram
equal to a given rectilineal figure.''

\begin{figure}
\begin{center}
\includegraphics{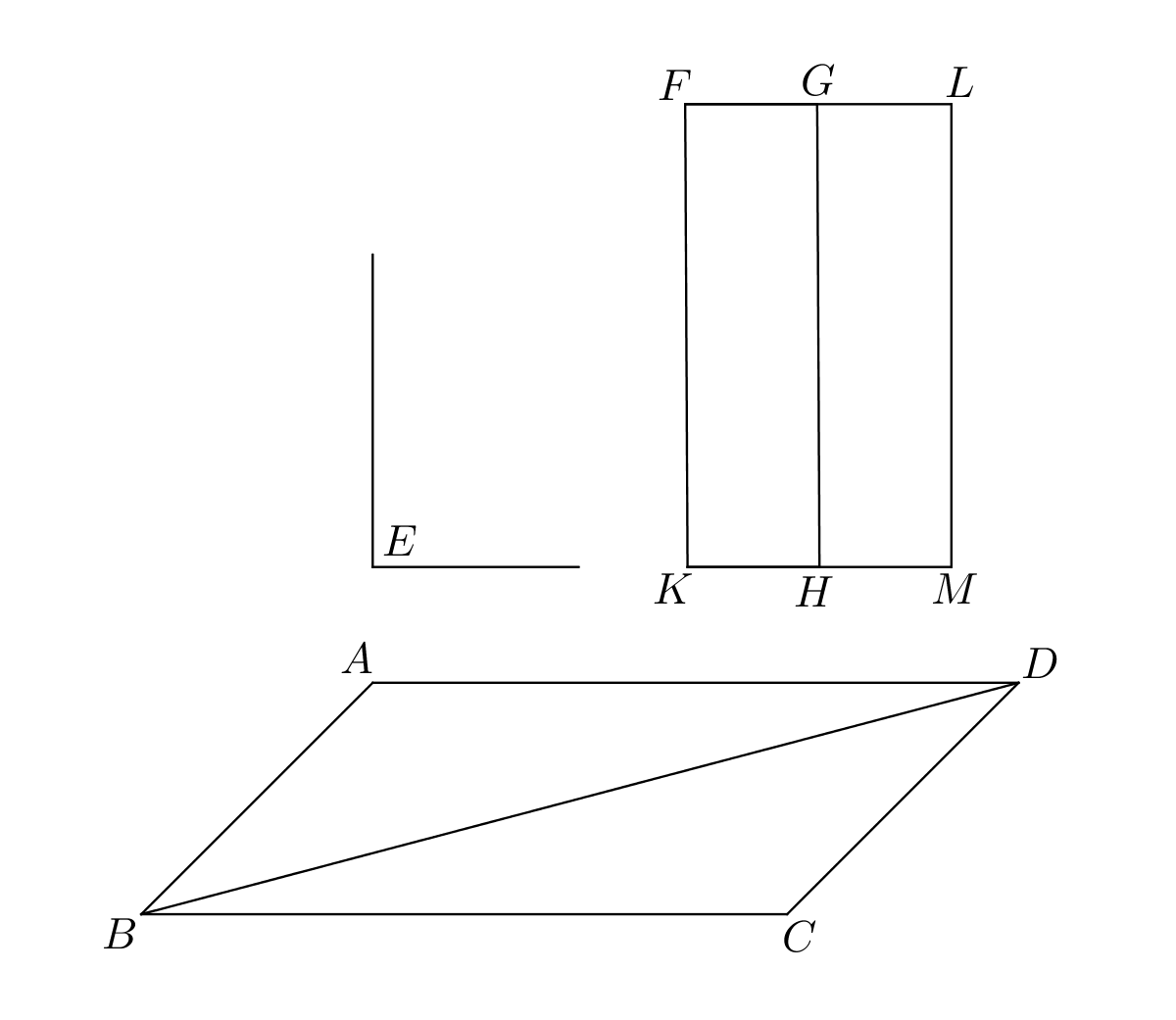}
\end{center}
\caption{I.45: P 38r, F 13r, B 30r, V 46}
\label{I45}
\end{figure}

Scholia for I.45 \cite[pp.~207--209]{euclidisV}.

Proclus 422--423 \cite[pp.~334--335]{proclus}, on {\em Elements} I.45:

\begin{quote}
For any 
rectilineal figure, as we said earlier, is as such divisible into
triangles, and we have given the method by which the number
of its triangles can be found. Therefore by dividing the
given rectilineal figure into triangles and constructing a
parallelogram equal to one of them, then applying parallelograms
equal to the others along the given straight line -- that
line to which we made the first application -- we shall have the
parallelogram composed of them equal to the rectilineal figure
composed of the triangles, and the assigned task will have been
accomplished. That is, if the rectilineal figure has ten sides,
we shall divide it into eight triangles, construct a parallelogram
equal to one of them, and then by applying in seven
steps parallelograms equal to each of the others, we shall
have what we wanted.

It is my opinion that this problem is what led the ancients
to attempt the squaring of the circle. For if a parallelogram
can be found equal to any rectilineal figure, it is worth
inquiring whether it is not possible to prove that a rectilineal 
figure is equal to a circular area. Indeed Archimedes proved
that a circle is equal to a right-angled triangle when its radius
is equal to one of the sides about the right angle and its
perimeter is equal to the base.
\end{quote}

Vitrac \cite[p.~278]{vitracI}:

\begin{quote}
La Prop. 45 montre donc comment construire, dans un angle rectiligne
donn\'e, un parall\'elogramme \'equivalent \`a une figure rectiligne donn\'ee : celle-ci
est divis\'ee en triangles et l'on construit pour chacun d'eux un parall\'elogramme
\'equivalent; il faut les agencer de mani\`ere \`a ce qu'ils forment,
ensemble, un parall\'elogramme, ce qu'Euclide v\'erifie soigneusement gr{\^a}ce
aux angles (Prop. 14) et \`a la theorie des parall\`eles (Prop. 29, 30, 33, 34).
En fait la question de l'agencement est r\'esolue seulement pour deux parall\'elogrammes,
mais le lecteur doit comprendre que la m\'ethode peut donner
lieu \`a it\'eration un nombre $n$ de fois, bien qu'\`a chaque \'etape on ne traite
que deux parall\'elogrammes \`a la fois. Dans les Livres arithm\'etiques, nous
retrouverons cette fa{\c c}on de proc\'eder qui n'est certainement pas un raisonnement
par induction, mais qui, en indiquant un sch\'ema g\'en\'eral par le
biais de la premi\`ere \'etape, justifie n\'eanmoins une \'enonciation qui, quant \`a
elle, est bien «universelle».
\end{quote}

I.46: ``On a given straight line to describe a square.''

\begin{figure}
\begin{center}
\includegraphics{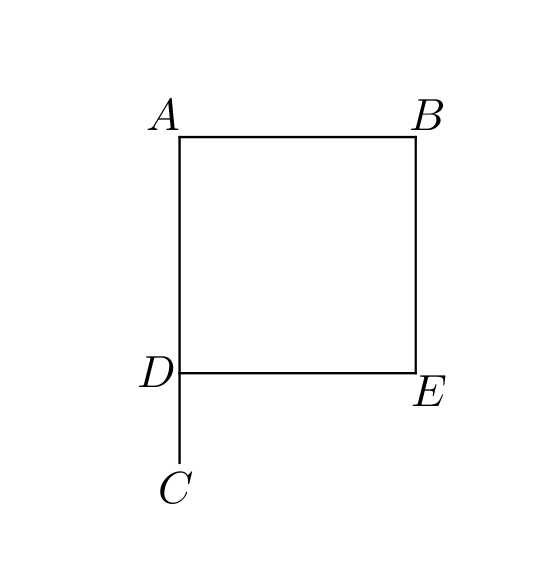}
\end{center}
\caption{I.46: P 38r, F 13r, B 30v, V 46}
\label{I46}
\end{figure}

\begin{proof}
{\em ekthesis}: Let $AB$ be the given straight line.

{\em diorismos}: It is required to describe a square on the straight line $AB$.

{\em kataskeu\={e}}: By I.11, draw a stright line $AC$ at right angles to the straight line $AB$.
By I.3, cut off $AD$ from $AC$ equal to $AB$.
By I.31, draw $DE$ parallel to $AB$ and by I.3 make this equal to $AB$.

{\em apodeixis}: 

{\em sumperasma}: 
\end{proof}

Scholia for I.46 \cite[pp.~209--212]{euclidisV}.

Campanus \cite[p.~92]{campanusI}, I.46 (I.45 in Campanus; Campanus does not include our I.45):

\begin{quote}
Ex data linea quadratam describere.
\end{quote}

\begin{figure}
\begin{center}
\includegraphics{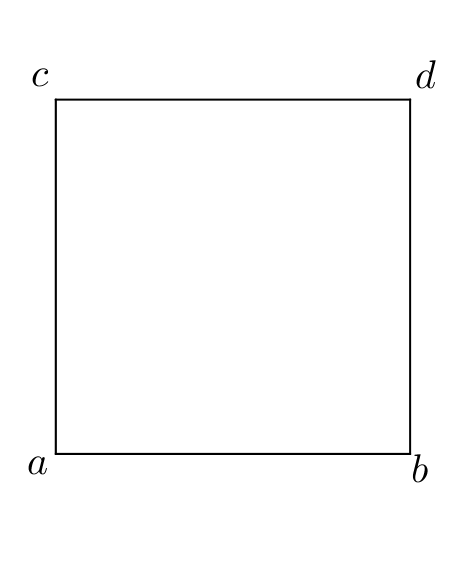}
\end{center}
\caption{Campanus, I.46 (I.45 in Campanus)}
\label{campanusI45}
\end{figure}

\begin{proof}
Let the given line be $ab$.
It is required to construct a square on the line $ab$.
By I.11, at the point $a$ let the line $ac$ be drawn at right angles to the line $ab$; likewise
by I.11, at the point $b$ let $bd$ be drawn at right angles to the line $ab$.
The line $ab$ falling on the lines $ac, bd$ makes interior angles $cab$ and $dba$ on one side; each is a right angle, so the interior angles
are equal to two right angles and then by I.28, the lines $ca,bd$ are parallel.
By I.3, cut down each of these to equal $ab$. 
Join $cd$. By I.33, because the lines $ca,bd$ are equal and parallel, the lines $cd,ab$ are equal and parallel.
By I.29, because both of the two angles $a$ and $b$ are right, both of the two angles $c$ and $d$ are right. Therefore by definition $abcd$
is a square.  

Another proof. By I.11, let $ac$ be perpendicular to $ab$ and let it be equal to $ab$. By I.31, at the point $c$, draw $cd$ parallel
with $ab$ and let it be equal to $ab$. Then draw $db$, which by I.33 it equal and parallel to $ac$.
By I.29 all the angles are right angles, and so by definition $abcd$ is a square.
\end{proof}

\section{Other writings}
Talking about Postulates 1--3 of the {\em Elements}, Mueller \cite[p.~16]{mueller} says:

\begin{quote}
Much more insight is obtained by examining the central proposition or propositions of book I and showing how the
book builds to its or their proof. In fact, almost the entire content of book I can be explained by reference to the construction of a 
parallelogram in a given angle and equal (in area) to a given rectilineal figure in proposition 45. This proposition makes it possible to represent
any rectilineal area as a rectangle. Euclid could have proved a stronger result, namely that any rectilineal area can be represented as a rectangle with a given base.
(Compare I,44.) From our point of view this result would be more interesting, since the areas of rectangles on equal bases are proportional
to the lengths of their sides. For the Greeks, however, the important representation of an area seems to be as a square, and I,45 is sufficient
for Euclid to be able to show in II,14 how to construct a square equal to any given rectilineal figure. This proposition represents
the true culmination of the geometry of the area of rectilineal figures. Euclid postpones it to book II because its proof involves methods which he introduces there and which he wishes, presumably for purposes of exposition, to separate from the methods of book I.
\end{quote}

\begin{description}
\item[I.2] Construct at a point a line equal to given line.
\item[I.11--12] Construct at a point a line perpendicular to given line.
\item[I.23] Construct at a point and on a line an angle equal to a given angle.
\item[I.31] Construct at a point a line parallel to a given line.
\item[I.35] Parallelograms on the same base in the same parallels are equal.
\item[I.41] A parallelogram on the same base and in the same parallels as a given triangle is double the given triangle.
\item[I.42] Construct in an angle a parallelogram equal to a given triangle.
\item[I.43] In a parallelogram, given two parallelograms about the diameter, the complements are equal.
\item[I.44] Construct in an angle and applied to a line a parallelogram equal to a given triangle.
\item[I.45] Triangulation and induction, invokes I.44.
\end{description}

Mugler \cite[p.~324]{mugler}: {\em paraballo}, {\em paraballein}

Plato, {\em Meno} 86e--87a \cite[p.~140]{guthrie}:

\begin{quote}
Just grant me one small relaxation of your sway, and allow me, in considering whether or not
it can be taught, to make use of a hypothesis -- the sort of thing, I mean, that geometers often
use in their inquiries. When they are asked, for example, about a
given area, whether it is possible
for this area to be inscribed as a triangle in a given circle, they will probably reply: `I don't know yet
whether it fulfils the conditions, but I think I have a hypothesis which will help us in the matter.
It is this. If the area is such that, when one has applied it [sc. as a rectangle] to the given line [i.e. the
diameter] of the circle, it is deficient by another rectangle similar to the one which is applied, then, I should
say, one result follows; if not, the result is different. If you ask me, then, about the inscription of the figure
in the circle -- whether it is possible or not -- I am ready to answer you in this hypothetical way.'
\end{quote}

Bluck \cite[p.~442]{bluck} writes in his commentary on the {\em Meno}: ``If a rectangle $ABCD$ is applied to a line $BH$ which is 
greater than the base of the rectangle, it was said to `fall short' by the area enclosed when $DCH$ is completed
as a rectangle. The same is true, {\em mutatis mutandis}, in the case of any parallelogram. This appears
not only from Euclid, but from a passage of Proclus ({\em Comm. in Eucl.} I, 44), in which this use of
{\em ἐλλείπειν} is attributed, on the authority of {\em οί περὶ τὸν Εὒδημον}, to early Pythagoreans.''

\begin{figure}
\begin{center}
\includegraphics{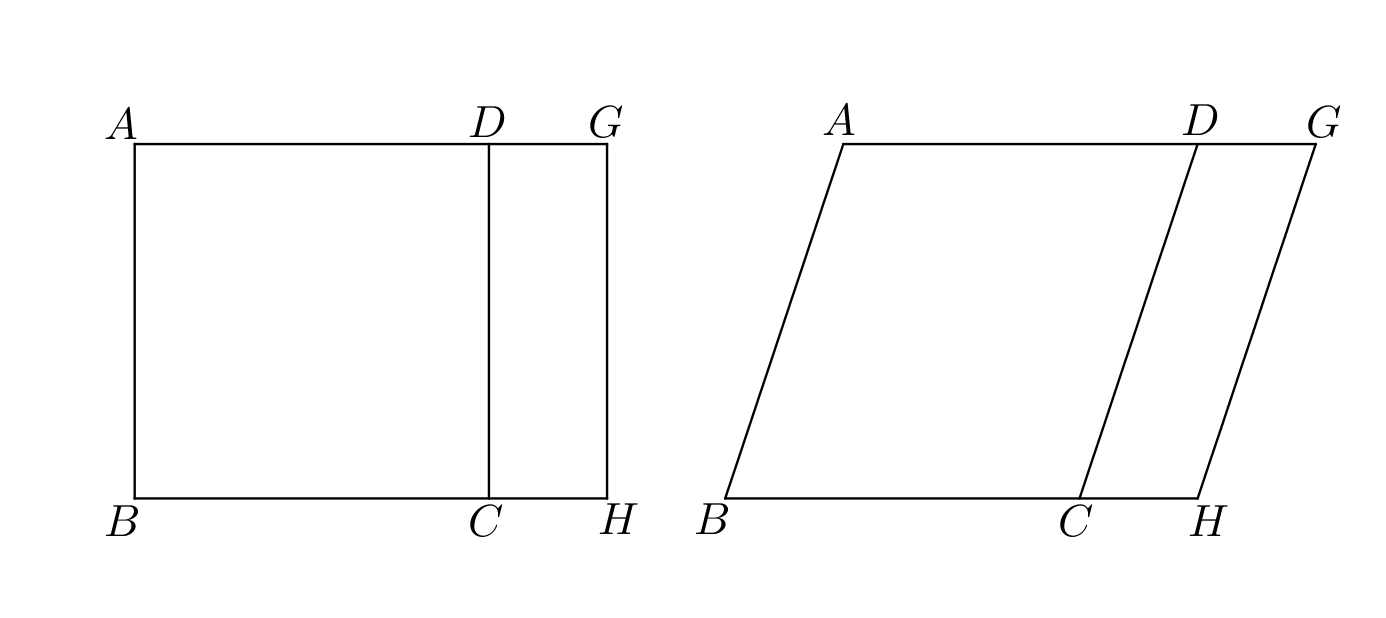}
\end{center}
\caption{Bluck, on {\em Meno} 86e--87a}
\label{bluck}
\end{figure}

Klein \cite[p.~206]{klein} gives the following translation of this passage:

\begin{quote}
Geometricians do often adopt the following kind of procedure. If, for example, one of them has to answer the question whether a certain amount of space (whatever its -- rectilinear -- boundaries) is capable of being fitted as a triangle into a given circular area (so that the three vertices will touch the circumference of that circular area), he may say: while I do not know whether this particular amount of space has that capability I believe I have something of a supposition ({\em hosper \dots tina hypothesin}) at hand which might be useful for the purpose. It is this: if that amount of space (which can always be transformed into a triangular or rectangular area) were to be such that he who ``stretches it along'' ({\em parateinanta}) its
({\em autou}) given line ``runs short'' ({\em elleipein}) of a space like the very one which had been ``stretched along'' (the given line), then, it seems to me, one thing would be the result, and another again, if it were impossible for him to go through this experience. And so I am disposed to tell you what will happen with regard to the inscription of your amount of space ({\em autou}) into the circle, whether it is impossible or not impossible, by way of ``hypothesizing'' ({\em hypothemenos}).
\end{quote}

Plato, {\em Republic} VII, 526e--527b \cite[p.~244]{cornfordrepublic}:

\begin{quote}
Socrates: So geometry 
will be suitable or not, according as it makes us contemplate reality
or the world of change.

Glaucon: That is our view.

Socrates: In this respect, then, no one who has even a slight acquaintance
with geometry will deny that the nature of this science is in flat
contradiction with the absurd language used by mathematicians,
for want of better terms. They constantly talk of `operations' like
`squaring,'  `applying,'  `adding,'
and so on, as if the whole subject were to
{\em do} something, whereas the true purpose of the whole subject is
knowledge -- knowledge, moreover, of what eternally exists, not of
anything that comes to be this or that at some time and ceases to be.
\end{quote}

Plutarch, {\em Quaestiones Convivales} VIII.2.4, 720A \cite[p.~177]{LCL335}:

\begin{quote}
Among the most geometrical theorems, or rather
problems, is this -- given two figures, to apply a third
equal to the one and similar to the other; it was
in virtue of this discovery they say Pythagoras
sacrificed. This is unquestionably more subtle and
elegant than the theorem which he proved that the
square on the hypotenuse is equal to the squares on
the sides about the right angle.
\end{quote}

Plutarch, {\em Non posse suaviter vivi secundum Epicurum} 1094B \cite[pp.~65--67]{LCL428}:

\begin{quote}
Our love of pleasure, to be sure, takes many forms and is enterprising enough;
but no one has so far upon having his way with the
woman he loves been so overjoyed that he sacrificed
an ox, nor has anyone prayed to die on the spot if he
could only eat his fill of royal meat or cakes; whereas
Eudoxus prayed to be consumed in flames like
Pha\"ethon if he could but stand next to the sun and
ascertain the shape, size, and composition of the
planets, and when Pythagoras discovered his theorem
he sacrificed an ox in honour of the occasion, as 
Apollodorus says:
\begin{quote}
When for the famous proof Pythagoras\\
Offered an ox in splendid sacrifice--
\end{quote}
whether it was the theorem that the square on the
hypotenuse is equal to the sum of the squares on 
the sides of the right  angle  or a problem about the
application of a given area.
\end{quote}

Aristotle, {\em De anima} II.2, 413a13--20 \cite[p.~191]{aristotle}:

\begin{quote}
It is not enough that the defining statement should make clear
the bare fact as most definitions do; it should also include and
exhibit the cause. As things are, what is stated in definitions is
usually of the nature of a conclusion. For instance, what is 
`squaring'? The construction of an equilateral rectangle equal (in area)
to a given oblong (rectangle). Such a definition is a statement of
the conclusion, whereas, if you say that squaring is the finding of
a mean proportional, you state the cause of the thing defined.
\end{quote}

Philoponus \cite[p.~34]{philoponus1} writes the following in his commentary on this passage:

\begin{quote}
It is clear that he is speaking of squaring an oblong. A square is both equilateral and a rectangle, that is, an area that has both the four sides equal and the angles right angles; an oblong is rectangular, indeed, but not equilateral. Those, then, who wish to square the oblong seek a mean proportional. What sort of thing do I mean? Let there be an oblong area having one side of eight cubits and the other of two. Clearly the whole is of 16 [square] cubits. For every quadrilateral is measured by multiplying side by side. If, therefore, we wish to
make a square equal to this oblong area, so as to be 16 cubits, the size the oblong was, we must find the mean proportional of the two sides of the oblong, so that it may have that ratio to the greater side, which was of 8 cubits, which the side of the oblong which was of 2 cubits has to it, the mean. Such [a mean] would be of 4 cubits. For that same ratio which 4 has to 8, 2 has to 4: each is half the greater. This is the mean proportional. On this, therefore, will be inscribed a square area of 16 cubits equal to the oblong. And thus should we do with every oblong when we want to inscribe a square equal to it. For again, if there should be an oblong having one side of 16 cubits and the other of 4, it inscribes an area, clearly, of 64. If you should want to make a square equal to this, seek a mean proportional. That is of 8 cubits. For 8 times 8 is 64. For just as the 16 cubit side of the oblong is the double of the 8 cubit [side] that has been found, so too this is the double of the remaining side of the oblong, which was of 4 cubits.

`Finding the mean?, Alexander says, `is shown in the second book of Euclid?. But it is not. Nothing of this sort is shown there, but in the sixth. There it is shown: `Given two straight lines, to find the mean proportional' [{\em Elements} 6.13], and ?If three straight lines are proportional, the [rectangle] contained by the extremes is equal to the
[square] on the middle? [{\em Elements} 6.17].
\end{quote}

Aristotle, {\em Metaphysics} III.2, 996b18--21 \cite[pp.~191--192]{aristotle}:

\begin{quote}
Again, in the case of other things, namely those which are the
subject of demonstration, we consider that we possess knowledge
of a particular thing when we know what it is [i.e. its definition],
e.g. we know what squaring is because (we know that) it is the
finding of the mean proportional.
\end{quote}

Iamblichus, {\em De communi mathematica scientia}, l.~21, Chapter XXIV \cite[p.~75]{festa}. Dillon and Urmson translation \cite[pp. 81-82]{DCMS}:

\begin{quote}
Next, we should discuss the group practice ({\em sunêtheia}) in the Pythagoreans'
occupation with mathematical pursuits. So then, they separated mathematical
arguments from the realm of the sensible, and converted the reasoning
faculty to an assurance of incorporeal being through mathematics. They
used it as a means of transport to the intelligible, and especially investigated
what mathematics contains akin to the pure forms and unitary reason principles.
Therefore they made use of this procedure in their inquiries, once
and for all separating the knowledge of mathematics from the common and
widespread level of understanding, and accordingly transmitted it as a secret.
They shared the knowledge of it with very few, and if some leak to the multitude
occurred, they abominated it as an impiety. Consequently they repelled also
those outside the group practice as being unworthy to share in mathematics.
For Pythagoras considered that they should not share their mathematical
knowledge with all, but only with those with whom one would share the whole
of one’s life. He did not accept people into this fellowship at random nor
casually, but he tested them over a long period, rejecting the unworthy.

Also, he did not communicate outside the group practice the advances made
by himself, but kept the discussions of mathematics secret from others; but he
brought about great advances among those called Pythagoreans by their
association with him, both in the understanding of mathematics and in the study
of geometry; also one would find that the basis of nearly all later advances came
to us from him. He esteemed in mathematics not, as do some later persons, the
power by which they will be able to solve a problem, but the insights themselves;
and, of these, not those most difficult to discover, like most of those later, but
those of them in which one could recognize order or some natural property.

This was their attitude because they thought that the principles of the whole
of nature were present in these, and that these were especially easy to grasp
-- what they were and how many -- because they were concerned with a nature
that was invariant and free from change, and was also simple. Therefore they did
not deal with issues that involved raising problems except those concerning
basic elements, such as geometrical application ({\em parabolé}) and the squaring
<of the circle>. Nor did they make it their business in their investigations to
exhaust every issue and leave no possibility unexplored, but sought only to see
the principles themselves in each case. They made a training in these sciences
and a rational working out that was exact and theoretical into a proper science.
They assigned a suitable ordering within the sciences, and they accepted few
things as principles and developed them. Also they especially perfected the
most valuable and most profound of the theorems; they applied theorems
differently studied to other matters, and set them in such an order as to propound
the more simple <first>, those requiring composition next. Also they arranged
them in accordance with the nature of things, suitably to our power and to the
worth of those who received them, in a way leading to excellence, in accord
with the whole of education and as proper for the purification of the soul.

Such, then, is the information one might give about this Pythagorean
method of procedure, about which we shall say more in the accounts that will
be given of each of the branches of mathematics individually. 
\end{quote}

Becker \cite[p.~59f.]{becker}

Knorr \cite{knorr}

Burkert \cite[p.~452]{burkert}:

\begin{quote}
The application of areas was known to Plato, but Hippocrates of Chios, for a problem soluble
by this method, used the method of ``inclination'' or ``verging'' (νεῦσις); it looks as though
the application of areas was at least not fully developed in Hippocrates.
\end{quote}

Archibald's reconstruction of Euclid's {\em Divisions of Figures} \cite{divisions}

Pappus, {\em Collection} 7.30--31 \cite[pp.~114--116]{pappus7a}:

\begin{quote}
[30] Apollonius, filled out Euclid's four books of {\em Conics} and added on another four, handing down eight volumes of {\em Conics}. Aristaeus, who wrote the five volumes of
{\em Solid Loci}, which have been transmitted until the present immediately following the {\em Conics}, and Apollonius's (other) predecessors, named the first of the three conic curves `section of an acute-angled cone', the second 1of a right-angled', the third of an obtuse-angled'. But since the three curves occur in each of these three cones, when cut
 variously, Apollonius was apparently at a loss to know why on earth his predecessors selectively named the one `section of an acute-angled cone' when it can also
 be (a section) of a right-angled and obtuse-angled one, the second (cone), and the third `of an obtuse-angled' when it can be of an
 acute-angled and a right-angled (cone), so, replacing the names, he called the (section) of an acute-angled (cone) `ellipse', that of a right-angled `parabola', and that of an
 obtuse-angled `hyperbola', each from a certain property of its own. For a certain area applied to a certain line, in the section of an acute-angled cone, falls short by a square, in that of an obtuse-angled (cone) exceeds by a square, but in that of a right-angled (cone) neither falls short nor exceeds.
 
 [31] This was his notion because he did not perceive that by a certain single way of having the plane cut the cone in generating the curves, a different one of the curves is produced in each of the cones, and they named it from the property of the cone. For if the cutting plane is drawn parallel to one side of the cone, one only of the three curves is formed, always the same one, which Aristaeus named a section of the (kind of) cone that was cut.
\end{quote}

Taisbak on application of areas \cite{taisbak2003}

Euclid's {\em Data} \cite{data}, Data 57

Friberg \cite{amazing}: pp.~114

\bibliographystyle{plain}
\bibliography{euclidI}

\begin{thebibliography}{10}

\bibitem{divisions}
Raymond~Clare Archibald.
\newblock {\em Euclid's Book on Division of Figures, with a Restoration Based
  on {W}oepcke's Text and on the {{\em Practica Geometriae} of Leonardo
  Pisano}}.
\newblock Cambridge University Press, 1915.

\bibitem{becker}
Oskar Becker.
\newblock {\em {Das mathematische Denken der Antike}}, volume~3 of {\em
  Studienhefte zur Altertumswissenschaft}.
\newblock Vandenhoeck \& Ruprecht, G\"ottingen, second edition, 1966.

\bibitem{bluck}
R.~S. Bluck.
\newblock {\em {Plato's Meno}}.
\newblock Cambridge University Press, 1961.

\bibitem{burkert}
Walter Burkert.
\newblock {\em Lore and Science in Ancient {P}ythagoreanism}.
\newblock Harvard University Press, 1972.
\newblock Translated from the German by Edwin L. Minar, Jr.

\bibitem{hermann}
H.~L.~L. Busard, editor.
\newblock {\em {The translation of the Elements of Euclid from the Arabic into
  Latin by Hermann of Carinthia (?)}}.
\newblock Brill, Leiden, 1968.

\bibitem{adelardI}
H.~L.~L. Busard, editor.
\newblock {\em {The first Latin translation of Euclid's {\em Elements} commonly
  ascribed to Adelard of Bath}}, volume~64 of {\em Studies and Texts}.
\newblock Pontifical Institute of Mediaeval Studies, Toronto, 1983.

\bibitem{adelardIII}
H.~L.~L. Busard, editor.
\newblock {\em {Johannes de Tinemue's redaction of Euclid's Elements, the
  so-called Adelard III version}}, volume 45,1 of {\em Boethius: Texte und
  Abhandlungen zur Geschichte der Mathematik und der Naturwissenschaften}.
\newblock Franz Steiner Verlag, Stuttgart, 2001.

\bibitem{campanusI}
H.~L.~L. Busard, editor.
\newblock {\em {Campanus of Novara and Euclid's {\em Elements}, volume I}}.
\newblock Franz Steiner Verlag, Stuttgart, 2005.

\bibitem{adelardII}
Hubert L.~L. Busard and Menso Folkerts, editors.
\newblock {\em Robert of {C}hester's (?) redaction of {E}uclid's {{\em
  Elements}}, the so-called {Adelard II} version. {Volume I}}, volume~8 of {\em
  Science Networks. Historical Studies}.
\newblock Springer Basel AG, 1992.

\bibitem{philoponus1}
William Charlton.
\newblock {\em {Philoponus: On Aristotle On the Soul 2.1--6}}.
\newblock Ancient Commentators on Aristotle. Bloomsbury, 2005.

\bibitem{cornfordrepublic}
Francis~Macdonald Cornford.
\newblock {\em The {{\em Republic}} of {P}lato, translated with introduction
  and notes}.
\newblock Oxford University Press, 1941.

\bibitem{dalimier}
Catherine Dalimier.
\newblock Apollonios {D}yscole sur la fonction des conjonctions expl\'etives.
\newblock {\em Revue des \'Etudes Grecques}, 112:719--730, 1999.

\bibitem{DCMS}
J.~O. Dillon, John Dillon~Urmson.
\newblock {\em {Iamblichus: On the General Science of Mathematics}}.
\newblock Ancient Commentators on Aristotle. Bloomsbury Academic, 2020.

\bibitem{LCL428}
Benedict Einarson and Phillip~H. De~Lacy.
\newblock {\em {Plutarch, Moralia, volume XIV}}.
\newblock Number 428 in Loeb Classical Library. Harvard University Press, 1967.

\bibitem{festa}
Nicola Festa, editor.
\newblock {\em {Iamblichi de communi mathematica scientia liber}}.
\newblock B. G. Teubner, Lipsiae, 1891.

\bibitem{amazing}
J\"oran Friberg.
\newblock {\em Amazing Traces of a {B}abylonian Origin in {G}reek Mathematics}.
\newblock World Scientific, 2007.

\bibitem{guthrie}
W.~K.~C. Guthrie.
\newblock {\em {Plato, Protagoras and Meno}}.
\newblock Penguin Classics. Penguin Books, 1956.

\bibitem{aristotle}
Thomas~L. Heath.
\newblock {\em Mathematics in {A}ristotle}.
\newblock Claredon Press, Oxford, 1949.

\bibitem{euclidI}
Thomas~L. Heath.
\newblock {\em The thirteen books of {E}uclid's {E}lements, volume {I}: {B}ooks
  {I}-{II}}.
\newblock Dover Publications, second edition, 1956.

\bibitem{euclidisV}
J.~L. Heiberg, editor.
\newblock {\em {Euclidis Opera omnia, vol. V}}.
\newblock B. G. Teubner, Lipsiae, 1888.

\bibitem{pappus7a}
Alexander Jones, editor.
\newblock {\em {Pappus of Alexandria. Book 7 of the {\em Collection}. Part 1.
  Introduction, Text, and Translation}}, volume~8 of {\em Sources in the
  History of Mathematics and Physical Sciences}.
\newblock Springer, 1986.

\bibitem{klein}
Jacob Klein.
\newblock {\em {A Commentary on Plato's {\em Meno}}}.
\newblock University of Chicago Press, 1989.

\bibitem{knorr}
Wilbur~Richard Knorr.
\newblock {\em The Ancient Tradition of Geometric Problems}.
\newblock Dover Publications, 2012.

\bibitem{gerard}
Anthony Lo~Bello.
\newblock {\em {Gerard of Cremona's translation of the commentary of al-Nayrizi
  on Book I of Euclid's {\em Elements of Geometry}}}, volume~2 of {\em Ancient
  Mediterranean and Medieval Texts and Contexts}.
\newblock Brill, 2003.

\bibitem{alnayriziI}
Anthony Lo~Bello.
\newblock {\em {The Commentary of al-Nayrizi on Book I of Euclid's {\em
  Elements of Geometry}, with an Introduction on the Transmission of Euclid's
  {\em Elements} in the Middle Ages}}, volume~1 of {\em Ancient Mediterranean
  and Medieval Texts and Contexts}.
\newblock Brill, 2003.

\bibitem{albertus}
Anthony Lo~Bello.
\newblock {\em {The Commentary of Albertus Magnus on Book I of Euclid's {\em
  Elements of Geometry}}}, volume~3 of {\em Ancient Mediterranean and Medieval
  Texts and Contexts}.
\newblock Brill, 2003.

\bibitem{proclus}
Glenn~R. Morrow.
\newblock {\em {Proclus: A Commentary on the First Book of Euclid's {\em
  Elements}}}.
\newblock Princeton University Press, 1992.

\bibitem{mueller}
Ian Mueller.
\newblock {\em Philosophy of Mathematics and Deductive Structure in {Euclid's
  {\em Elements}}}.
\newblock Dover Publications, 2006.

\bibitem{mugler}
Charles Mugler.
\newblock {\em {Dictionnaire historique de la terminologie g\'eom\'etrique des
  Grecs}}, volume XXVIII of {\em \'Etudes et commentaires}.
\newblock Librairie C. Klincksieck, Paris, 1958.

\bibitem{netz}
Reviel Netz.
\newblock {\em The shaping of deduction in {G}reek mathematics: a study in
  cognitive history}, volume~51 of {\em Ideas in Context}.
\newblock Cambridge University Press, 1999.

\bibitem{szabo}
\'Arp\'ad Szab\'o.
\newblock {\em The beginnings of {G}reek mathematics}, volume~17 of {\em
  Sythese Historical Library}.
\newblock D. Reidel Publishing Company, Dordrecht, 1978.
\newblock Translated from the German by A. M. Ungar.

\bibitem{taisbak2003}
Christian~Marinus Taisbak.
\newblock Exceeding and falling short: Elliptical and hyperbolical application
  of areas.
\newblock {\em Science in Context}, 16(3):299--318, 2003.

\bibitem{data}
Christian~Marinus Taisbak.
\newblock {\em {Euclid's {\em Data}}, or The Importance of Being Given},
  volume~45 of {\em Acta Historica Scientiarum Naturalium et Medicinalium}.
\newblock Museum Tusculanum Press, Copenhagen, 2010.

\bibitem{LCL335}
Ivor Thomas, editor.
\newblock {\em Greek Mathematics, volume {I}}.
\newblock Number 362 in Loeb Classical Library. Harvard University Press, 1991.

\bibitem{vitracI}
Bernard Vitrac and Maurice Caveing.
\newblock {\em {Euclide d'Alexandrie. Les \'El\'ements. Volume 1: Introduction
  g\'en\'erale et Livres I--IV}}.
\newblock Biblioth\`eque d'histoire des sciences. Presses Universitaires de
  France, Paris, 1990.

\bibitem{hippias}
Robin Waterfield.
\newblock {Hippias Major and Hippias Minor}.
\newblock In Trevor~J. Saunders, editor, {\em {Plato, Early Socratic
  Dialogues}}, Penguin Classics, pages 211--293. Penguin Books, 2005.

\end{thebibliography}

\end{document}